\pdfoutput=1

\documentclass[titlepage,twoside,openright,10pt]{book}

\usepackage{amsthm,mathrsfs,bbm,amssymb,here,amsmath,cite,textcomp,stmaryrd}
\usepackage[a5paper, top=.8in, bottom=1.3in, left=.5in, right=.5in]{geometry}
\usepackage{fancyhdr}
\usepackage{makeidx}
\usepackage{graphicx}
\usepackage{xr}

\makeindex

\newcommand{\be}{\begin{equation}}
\newcommand{\ee}{\end{equation}}
\newcommand{\bea}{\begin{eqnarray}}
\newcommand{\eea}{\end{eqnarray}}
\newcommand{\ba}{\begin{array}}
\newcommand{\ea}{\end{array}}
\newcommand{\mc}{\mathcal}
\newcommand{\tp}{\mathcal T}

\newcommand{\eopsymbol}{$\blacksquare$}
\newcommand{\naturalnumbers}{{\mathbb N}}
\newcommand{\bigtimes}{\raisebox{-.2em}{\huge $\times$}}

\newcommand{\midvspace}{\vspace{.25in}}

\newtheoremstyle{alltheorem}
	{0in}	
	{0in}	
	{}		
	{}		
	{\bf}	
	{}		
	{\newline}	
	{}		

\theoremstyle{alltheorem}
\newtheorem{defi}{Definition}[chapter]
\newtheorem{prop}[defi]{Proposition}
\newtheorem{rema}[defi]{Remark}
\newtheorem{axio}[defi]{Axiom}
\newtheorem{lemm}[defi]{Lemma}
\newtheorem{theo}[defi]{Theorem}
\newtheorem{coro}[defi]{Corollary}
\newtheorem{exam}[defi]{Example}
\newtheorem{lede}[defi]{Lemma and Definition}

\newcommand{\alltheoremwidth}{\textwidth}

\newcommand{\bdefi}{
\begin{minipage}{\alltheoremwidth}
\begin{defi}}

\newcommand{\edefi}{
\hspace{\stretch{1}} \eopsymbol \end{defi} \end{minipage} \vspace{.25in}}

\newcommand{\bexam}{
\begin{minipage}{\alltheoremwidth}
\begin{exam}}

\newcommand{\eexam}{
\hspace{\stretch{1}} \eopsymbol \end{exam} \end{minipage} \vspace{.25in}}

\newcommand{\brema}{
\begin{minipage}{\alltheoremwidth}
\begin{rema}}

\newcommand{\erema}{
\hspace{\stretch{1}} \eopsymbol \end{rema} \end{minipage} \vspace{.25in}}

\newcommand{\baxio}{
\begin{minipage}{\alltheoremwidth}
\begin{axio}}

\newcommand{\eaxio}{
\hspace{\stretch{1}} \eopsymbol \end{axio} \end{minipage} \vspace{.25in}}

\newcommand{\blede}{
\begin{minipage}{\alltheoremwidth}
\begin{lede}}

\newcommand{\elede}{
\end{lede} \end{minipage} \vspace{.0in}}

\newcommand{\blemm}{
\begin{minipage}{\alltheoremwidth}
\begin{lemm}}

\newcommand{\elemm}{
\end{lemm} \end{minipage} \vspace{.0in}}

\newcommand{\btheo}{
\begin{minipage}{\alltheoremwidth}
\begin{theo}}

\newcommand{\etheo}{
\end{theo} \end{minipage} \vspace{.0in}}

\newcommand{\bprop}{
\begin{minipage}{\alltheoremwidth}
\begin{prop}}

\newcommand{\eprop}{
\end{prop} \end{minipage} \vspace{.0in}}

\newcommand{\bcoro}{
\begin{minipage}{\alltheoremwidth}
\begin{coro}}

\newcommand{\ecoro}{
\end{coro} \end{minipage} \vspace{.0in}}

\newcommand{\bproof}{
\vspace{-.05in} \begin{proof}}

\newcommand{\eproof}{
\end{proof} \vspace{.15in}}

\newcommand{\benum}{
\renewcommand{\theenumi}{\roman{enumi}}
\renewcommand{\labelenumi}{(\theenumi)}
\vspace{-.05in} \begin{enumerate}}

\newcommand{\eenum}{\end{enumerate}}

\usepackage[pdftex]{hyperref}

\fancypagestyle{plain}{
\lhead[\thepage]{}
\chead{}
\rhead[]{\thepage}
\lfoot{\ \\ \small \textcopyright \ 2013 \, Felix Nagel --- Set theory and topology, Part II: Topology -- Fundamental notions}
\cfoot{}
\rfoot{}
}

\begin{document}

\frontmatter


\thispagestyle{empty}

\vspace{\baselineskip}

\begin{center}
\textbf{\huge Set theory and topology\\[.3em]
\large An introduction to the foundations of analysis
\footnote{For remarks and suggestions please contact: stt.info@t-online.de}\\[3em]
\large Part II: \; Topology -- Fundamental notions}\\
\vspace{7\baselineskip}
{\sc Felix~Nagel}\\
\vspace{5\baselineskip}
\textbf{Abstract}\\
\vspace{1\baselineskip}
\parbox{0.9\textwidth}{We provide a formal introduction into the classic theorems of general topology and its axiomatic foundations in set theory. In this second part we introduce the fundamental concepts of topological spaces, convergence, and continuity, as well as their applications to real numbers. Various methods to construct topological spaces are presented.}
\end{center}
\vspace{\baselineskip}

\pagebreak

\thispagestyle{empty}

{\small The author received his doctoral degree from the University of Heidelberg for his thesis in electroweak gauge theory. He worked for several years as a financial engineer in the financial industry. His fields of interest are probability theory, foundations of analysis, finance, and mathematical physics. He lives in Wales and Lower Saxony.}


\chapter{Remark}
\label{sec-preface}

\setcounter{page}{1}

\thispagestyle{plain}

This is the second part of a series of articles on the foundations of analysis, cf.~\cite{Nagel}. For the Preface and Chapters~\ref{axiomatic foundation} to~\ref{numbers ii} see Part~\ref{part sets relations numbers}.

\pagebreak

\thispagestyle{plain}


\cleardoublepage
\phantomsection
\thispagestyle{plain}
\addcontentsline{toc}{chapter}{Table of contents}
\tableofcontents
\thispagestyle{plain}

\mainmatter


\setcounter{part}{1}
\setcounter{chapter}{4}

\part{Topology -- Fundamental notions}

\chapter{Topologies and filters}
\label{topologies}
\setcounter{equation}{0}

\pagebreak

\section{Set systems}
\label{set systems}

In this Section we introduce three basic functions on set systems, that are used in many places subsequently. Each of these functions is defined with respect to a given set~$X$ as a function that maps every subsystem of~${\mc P}(X)$ on a---generally larger---subsystem of~${\mc P}(X)$. A fourth function is introduced at the end of this Section, which is used in the context of neighborhood system in Section~\ref{neighborhood}.

The introduction of the following new symbol turns out to be useful.

\midvspace

\bdefi
Given a set $X$, we write $A \sqsubset X$ if $A \subset X$ and $A$ is finite.
\edefi

\bdefi
Given a set $X$, we define the following functions:
\begin{center}
\begin{tabular}{ll}
$\Psi_X : {\mc P}^2(X) \longrightarrow {\mc P}^2(X)$, \quad &
$\Psi_X ({\mc A}) = \left\{ \, \bigcap {\mc B} \, : \, {\mc B} \sqsubset {\mc A}, \; {\mc B} \neq \O \right\}$\,;\\[.8em]
$\Theta_X : {\mc P}^2(X) \longrightarrow {\mc P}^2(X)$, \quad &
$\Theta_X ({\mc A}) = \left\{ \, \bigcup {\mc B} \, : \, {\mc B} \subset {\mc A}, \; {\mc B} \neq \O \right\}$\,;\\[.8em]
$\Phi_X : {\mc P}^2(X) \longrightarrow {\mc P}^2(X)$, \quad &
$\Phi_X ({\mc A}) = \left\{ B \subset X \, : \, \exists A \in {\mc A} \;\; A \subset B \right\}$\\
\end{tabular}
\end{center}

When the set $X$ we refer to is evident from the context, we also use the short notations $\Psi$, $\Theta$, and $\Phi$, respectively.
\edefi

That is, for a system $\mc A$ of subsets of~$X$, $\Psi ({\mc A})$ is the system of all finite intersections of members of~$\mc A$, $\Theta ({\mc A})$ is the system of all unions of members of~$\mc A$, and $\Phi ({\mc A})$ is the system of all subsets of~$X$ that contain some member of~$\mc A$.

\midvspace

\brema
\label{rema set relations}
Given a set $X$, the following equations hold:
\benum
\item \label{rema set relations 1} $\Psi (\O) = \O \, , \quad \Theta (\O) = \O \, , \quad \Phi (\O) = \O$
\item \label{rema set relations 2} $\Psi \circ \Psi = \Psi \, , \quad \Theta \circ \Theta = \Theta \, , \quad \Phi \circ \Phi = \Phi$
\item \label{rema set relations 3} $\Psi (\left\{ \O \right\}) = \left\{ \O \right\} \, , \quad \Theta (\left\{ \O \right\}) = \left\{ \O \right\} \, , \quad \Phi (\left\{ \O \right\}) = {\mc P}(X)$
\item \label{rema set relations 4} $\Psi (\left\{ X \right\}) = \left\{ X \right\} \, , \quad \Theta (\left\{ X \right\}) = \left\{ X \right\} \, , \quad \Phi (\left\{ X \right\}) = \left\{ X \right\}$
\item \label{rema set relations 5} $\Psi (\left\{ \O, X \right\}) = \left\{ \O, X \right\} \, , \quad \Theta (\left\{ \O, X \right\}) = \left\{ \O, X \right\} \, , \quad \Phi (\left\{ \O, X \right\}) = {\mc P}(X)$
\eenum

\erema

The identities in Remark~\ref{rema set relations}~(\ref{rema set relations 2}) say that $\Psi$, $\Theta$, and $\Phi$ are projective.

The composition of two of the functions is not commutative. However, we have the following result.

\midvspace

\blemm
\label{lemm set relations}
Given a set $X$, we have for every ${\mc A} \subset {\mc P}(X)$:
\[
\Psi \, \Theta \, ({\mc A}) \, \subset \, \Theta \, \Psi \, ({\mc A}) \, , \quad \quad
\Psi \, \Phi \, ({\mc A}) \, \subset \, \Phi \, \Psi \, ({\mc A})
\]

The maps $(\Theta \, \Psi)$ and $(\Phi \, \Psi)$ are projective.
\elemm

\bproof
In order to prove the first claim, let $A = \bigcap_{i=1}^n \bigcup \left\{ A_{ij} \, : \, j \in J_i \right\}$ where $n \in \mathbb{N}$, $n > 0$, and, for every $i \in \naturalnumbers$, $1 \leq i \leq n$, $J_i$ is an index set and $A_{ij} \in {\mc A}$ ($j \in J_i$). We have
\[
A \, = \, \bigcup \left\{ \, {\bigcap}_{k = 1}^n A_{k j(k)} \; : \; j \in \bigtimes_{\!\! i = 1}^{\!\! n} \, J_i \right\} \, \in \, \Theta \, \Psi \, ({\mc A})
\]

To show the second claim, let $n \in \naturalnumbers$, $n > 0$, and for every $i \in \naturalnumbers$, $1 \leq i \leq n$, let $A_i \in {\mc A}$ and $B_i$ be a set with $A_i \subset B_i$. Further let $B = \bigcap_{i = 1}^n B_i$. It follows that $B \supset \bigcap_{i = 1}^n A_i$, and thus $B \, \in \, \Phi \, \Psi \, ({\mc A})$.

Now the last claim clearly follows.
\eproof

\blemm
\label{lemm equivalences set operations}
Given a set $X$, $({\mc P}^2(X),\subset)$ is an ordered space in the sense of~"$\leq$". In particular, $\subset$ is a reflexive pre-ordering. The maps $\Psi$, $\Theta$, and $\Phi$ as well as their compositions are $\subset$-increasing. For every ${\mc A}, {\mc B} \subset {\mc P}(X)$ we have:

\benum
\item \label{lemm equivalences set operations 1}
${\mc A} \subset_{\stackrel{}{\Psi}} {\mc B} \quad \Longleftrightarrow \quad \Psi ({\mc A}) \subset \Psi ({\mc B}) \quad \Longleftrightarrow \quad {\mc A} \subset \Psi ({\mc B})\\[.5em]
{} \quad \Longleftrightarrow \quad \forall A \in {\mc A} \quad \exists {\mc G} \sqsubset {\mc B} \quad {\mc G} \neq \O, \; A = \bigcap {\mc G}$

\vspace{.5em}

\item \label{lemm equivalences set operations 2}
${\mc A} \subset_{\stackrel{}{\Theta}} {\mc B} \quad \Longleftrightarrow \quad \Theta ({\mc A}) \subset \Theta ({\mc B}) \quad \Longleftrightarrow \quad {\mc A} \subset \Theta ({\mc B})\\[.5em]
{} \quad \Longleftrightarrow \quad \forall A \in {\mc A} \quad \exists {\mc H} \subset {\mc B} \quad {\mc H} \neq \O, \; A = \bigcup {\mc H}$

\vspace{.5em}

\item \label{lemm equivalences set operations 3}
${\mc A} \subset_{\stackrel{}{\Phi}} {\mc B} \quad \Longleftrightarrow \quad \Phi ({\mc A}) \subset \Phi ({\mc B}) \quad \Longleftrightarrow \quad {\mc A} \subset \Phi ({\mc B})\\[.5em]
{} \quad \Longleftrightarrow \quad \forall A \in {\mc A} \quad \exists B \in {\mc B} \quad A \supset B$

\vspace{.5em}

\item \label{lemm equivalences set operations 4}
${\mc A} \subset_{\stackrel{}{\Theta \Psi}} {\mc B} \quad \Longleftrightarrow \quad \Theta \, \Psi \, ({\mc A}) \, \subset \, \Theta \, \Psi \, ({\mc B}) \quad \Longleftrightarrow \quad {\mc A} \subset \, \Theta \, \Psi \, ({\mc B})\\[.5em]
{} \quad \Longleftrightarrow \quad \forall A \in {\mc A} \quad \exists {\mc H} \subset \Psi ({\mc B}) \quad {\mc H} \neq \O, \; A = \bigcup {\mc H}$

\vspace{.5em}

\item \label{lemm equivalences set operations 5}
${\mc A} \subset_{\stackrel{}{\Phi \Psi}} {\mc B} \quad \Longleftrightarrow \quad \Phi \, \Psi \, ({\mc A}) \, \subset \, \Phi \, \Psi \, ({\mc B}) \quad \Longleftrightarrow \quad {\mc A} \subset \, \Phi \, \Psi \, ({\mc B})\\[.5em]
{} \quad \Longleftrightarrow \quad \forall A \in {\mc A} \quad \exists {\mc G} \sqsubset {\mc B} \quad {\mc G} \neq \O, \; A \supset \bigcap {\mc G}$
\eenum

\elemm

\bproof
In each case the first equivalence is true by definition of the respective relation, cf.\ Definition~\ref{defi gen relation}.

To see the second equivalence notice that $\Psi$, $\Theta$, and $\Phi$ are projective by Remark~\ref{rema set relations}. The compositions $(\Theta \, \Psi)$ and $(\Phi \, \Psi)$ are projective by Lemma~\ref{lemm set relations}. The second equivalence follows by Lemma~\ref{lemm increasing projection}.

The third equivalence in each case is a consequence of the definition of the maps.
\eproof

\bdefi
\label{defi phi prime}
Given a set $X$, the function $\Phi'_X$ is defined by
\begin{center}
\begin{tabular}{l}
$\Phi'_X : {\mc P}\big(X \!\times {\mc P}(X)\big) \longrightarrow {\mc P}\big(X \!\times {\mc P}(X)\big)$ ,\\[.8em]
$\Phi'_X (R) = \big\{ (x,B) \in X \!\times {\mc P}(X) \, : \, \exists A \subset X \;\; (x,A) \in R, \, A \subset B \big\}$
\end{tabular}
\end{center}

When $X$ is evident from the context, we also use the short notation $\Phi'$ for~$\Phi'_X$.
\edefi

\blemm
\label{lemm phi prime relations}
Given a set $X$, $x \in X$, $A \subset X$, and $R \subset X \!\times {\mc P}(X)$, we have:
\benum
\item \label{lemm phi prime relations 1} $\left( \Phi'(R) \right) \left\{ x \right\} = \Phi ( R \left\{ x \right\} )$
\item \label{lemm phi prime relations 2} $\left( \Phi'(R) \right)  \left[ A \right] = \Phi ( R \left[ A \right] )$
\item \label{lemm phi prime relations 3} $\left( \Phi'(R) \right) \langle A \rangle \supset \Phi ( R \langle A \rangle )$
\eenum

\elemm

\bproof
(\ref{lemm phi prime relations 1}) and~(\ref{lemm phi prime relations 2}) clearly hold. To prove~(\ref{lemm phi prime relations 3}) notice that
\begin{center}
\begin{tabular}{l}
$\left( \Phi'(R) \right) \langle A \rangle \, = \, {\displaystyle \bigcap}_{x \in A} \left( \Phi'(R) \right) \left\{ x \right\} \, = \, {\displaystyle \bigcap}_{x \in A} \Phi ( R \left\{ x \right\} )$\\[.8em]
${} \quad  = \, \big\{ B \subset X \, : \, \forall x \in A \quad \exists B_x \subset B \quad (x,B_x) \in R \big\}$\\[.8em]
${} \quad  \, \supset \big\{ B \subset X \, : \, \exists C \subset B \quad \forall x \in A \quad (x,C) \in R \big\} \, = \, \Phi ( R \langle A \rangle )$
\end{tabular}
\end{center}
\eproof

\blemm
Given a set $X$, the pair $\left({\mc P}(X \!\times {\mc P}(X)),\subset \right)$ is an ordered space in the sense of~"$\leq$". In particular, $\subset$ is a reflexive pre-ordering on ${\mc P}(X \!\times {\mc P}(X))$. The map $\Phi'$ is $\subset$-increasing and projective. For every $R, S \subset X \!\times {\mc P}(X)$ we have:
\begin{center}
\begin{tabular}{l}
$R \subset_{\stackrel{}{\Phi'}} S \quad \Longleftrightarrow \quad \Phi' (R) \subset \Phi' (S) \quad \Longleftrightarrow \quad R \subset \Phi' (S)$\\[.8em] 
${} \quad \Longleftrightarrow \quad \forall (x,A) \in R \quad \exists (y,B) \in S \quad x = y, \, A \supset B$\\[.8em]
${} \quad \Longleftrightarrow \quad \forall x \in X \quad \left( \Phi' (R) \right) \left\{ x \right\} \, \subset \, \left( \Phi' (S) \right) \left\{ x \right\}$\\[.8em]
${} \quad \Longleftrightarrow \quad \forall x \in X \quad \Phi (R \left\{ x \right\}) \, \subset \, \Phi (S \left\{ x \right\})$\\[.8em]
${} \quad \Longleftrightarrow \quad \forall x \in X \quad R \left\{ x \right\} \subset_{\stackrel{}{\Phi}} S \left\{ x \right\}$
\end{tabular}
\end{center}
\elemm

\bproof
Exercise.
\eproof

\section{Topologies, bases, subbases}

We start with the definition of a topological space.

\midvspace

\bdefi
\label{defi top}
\index{Topology}
\index{Topological space}
\index{Open!set}
\index{Closed!set}
Given a set $X$, a system $\tp \subset {\mc P}(X)$ is called {\bf topology on}~$X$ if it has all of the following properties:
\benum
\item \label{defi top 1} $\O, X \in \tp$
\item \label{defi top 2} $\forall \, {\mc G} \subset \tp \quad {\mc G} \neq \O \; \Longrightarrow \; \bigcup {\mc G} \in \tp$
\item \label{defi top 3} $\forall A, B \in \tp \quad A \cap B \in \tp$
\eenum

The pair $\xi = (X, \tp)$ is called {\bf topological space}. The members of $\tp$ are called $\xi$-{\bf open} or $\tp$-{\bf open}. They are also called {\bf open} if the topology is evident from the context. A set $B \subset X$ is called $\xi$-{\bf closed} if $X \!\setminus\! B$ is $\xi$-open. If the set $X$ is evident from the context, we also say that $B$ is $\tp$-{\bf closed}. If the set $X$ and the topology $\tp$ are both evident from the context, we also say that $B$ is~{\bf closed}.
\edefi

Notice that property~(\ref{defi top 2}) in Definition~\ref{defi top} is equivalent to $\Theta (\tp) = \tp$. By property (\ref{defi top 3}) it follows that $\bigcap {\mc H} \in \tp$ for every ${\mc H} \sqsubset \tp$ with ${\mc H} \neq \O$ by the Induction principle. Therefore property~(\ref{defi top 3}) is equivalent to $\Psi (\tp) = \tp$. Hence the system of topologies on $X$ contains precisely the fixed points of $\Theta$ and $\Psi$ that additionally satisfy property~(\ref{defi top 1}).

We now define several simple topologies that serve as examples throughout the text.

\midvspace

\blede 
\label{lede examples topologies}
\index{Topology!discrete}
\index{Topology!indiscrete}
\index{Topology!cofinite}
\index{Topology!cocountable}
Given a set $X$, each of the following systems of subsets is a topology on~$X$:
\benum
\item \label{distop} $\tp_{\mathrm{ dis}} = {\mc P}(X)$ is called {\bf discrete topology}.
\item \label{intop} $\tp_{\mathrm{in}} = \left\{ \O, X \right\}$ is called {\bf indiscrete topology}.
\item \label{coftop} $\tp_{\mathrm{cf}} = \left\{ A \subset X \, : \, A^c \mathrm{\; is \; finite}\right\} \cup \left\{ \O \right\}$ is called {\bf cofinite topology}.
\item \label{cctop} $\tp_{\mathrm{cc}} = \left\{ A \subset X \, : \, A^c \mathrm{\; is \; countable}\right\} \cup \left\{ \O \right\}$ is called {\bf cocountable topology}.
\eenum

\elede

\bproof
Exercise.
\eproof

\bexam
Let $X$ be a set and $A_n \subset X$ ($n \in {\mathbb N}$) such that $A_n \subset A_{n+1}$ ($n \in {\mathbb N}$) and $\bigcup_{n \in {\mathbb N}} A_n = X$. Then ${\mc A} = \left\{ A_n : n \in {\mathbb N} \right\} \cup \left\{ \O, X \right\}$ is a topology on~$X$.
\eexam

The analogues of properties (\ref{defi top 1}) to (\ref{defi top 3}) in Definition~\ref{defi top} hold for the system of closed sets as follows.

\midvspace

\blemm
\label{systemclosed}
Let $\xi = (X,\tp)$ be a topological space. The system $\mc C$ of all $\xi$-closed sets has the following properties:
\benum
\item \label{systemclosed1} $\O, X \in {\mc C}$
\item \label{systemclosed2} $\forall \, {\mc G} \subset {\mc C} \quad {\mc G} \neq \O \; \Longrightarrow \; \bigcap {\mc G} \in {\mc C}$
\item \label{systemclosed3} $\forall A, B \in {\mc C} \quad A \cup B \in {\mc C}$
\eenum

\elemm

\bproof
Exercise.
\eproof

Clearly, property (\ref{systemclosed3}) in Lemma~\ref{systemclosed} implies that $\bigcup {\mc H} \in {\mc C}$ for every ${\mc H} \sqsubset {\mc C}$ with ${\mc H} \neq \O$ by the Induction principle (or by the analogue for open sets).

The following Lemma demonstrates that for a given system $\mc C$ of subsets of~$X$ one may first confirm that $\mc C$ is the system of all closed subsets for some topology on~$X$ and then construct the topology from~$\mc C$.

\midvspace

\blemm
\label{lemm closed sets}
Given a set $X$ and a system $\mc C$ of subsets of $X$ satisfying (\ref{systemclosed1}) to (\ref{systemclosed3}) in Lemma~\ref{systemclosed}, the system $\tp = \left\{ B^c : \, B \in {\mc C} \right\}$ is a topology on~$X$, and $\mc C$ is the system of all $\tp$-closed sets.
\elemm

\bproof
Exercise.
\eproof

We often encounter more than one topology on the same set. The notions in the following Definition are useful in this case.

\midvspace

\bdefi
\label{defi top comp}
\index{Topology!finer}
\index{Topology!coarser}
\index{Topology!strictly finer}
\index{Topology!strictly coarser}
\index{Topology!comparable}
Let $X$ be a set and $\tp_1$ and $\tp_2$ be two topologies on~$X$. If $\tp_2 \subset \tp_1$, then $\tp_1$ is called {\bf finer than}~$\tp_2$, and $\tp_2$ is called {\bf coarser than}~$\tp_1$. If $\tp_1 \subset \tp_2$ or $\tp_2 \subset \tp_1$, then $\tp_1$ and $\tp_2$ are called {\bf comparable}. If $\tp_2 \subset \tp_1$ and $\tp_1 \neq \tp_2$, then $\tp_1$ is called {\bf strictly finer than}~$\tp_2$, and $\tp_2$ is called {\bf strictly coarser than}~$\tp_1$.
\edefi

\brema
\label{top finer}
Let $(X,\tp)$ be a topological space. Then we have $\tp_{\mathrm{in}} \subset \tp \subset \tp_{\mathrm{dis}}$.
\erema

\blemm
\label{uncinc}
Given an uncountable set~$X$, the topologies defined in Lemma~\ref{lede examples topologies} (\ref{distop})--(\ref{cctop}) obey $\tp_{\mathrm{in}} \subset \tp_{\mathrm{cf}} \subset \tp_{\mathrm{cc}} \subset \tp_{\mathrm{dis}}$, and no two of them are identical.
\elemm

\bproof
Exercise.
\eproof

The following notation allows us a more formal treatment in the sequel.

\midvspace

\bdefi
Given a set~$X$, the system of all topologies on~$X$ is denoted by~$\mathscr{T}(X)$.
\edefi

\blede
\label{lede coarsest finest}
Let $X$ be a set, $\mathscr{A} \subset \mathscr{T}(X)$, and~$\tp \in \mathscr{A}$. The pair $(\mathscr{A}, \subset)$ is a space ordered in the sense of~"$\leq$". $\tp$~is called {\bf finest} ({\bf coarsest}) topology of~$\mathscr{A}$ if it is a maximum (minimum) of~$\mathscr{A}$. $\mathscr{A}$~has at most one finest and at most one coarsest topology.
\elede

\bproof
Exercise.
\eproof

\brema
\label{rema top bounds}
The discrete topology on~$X$ is the finest member of~$\mathscr{T}(X)$, and thus it is an upper bound of any subsystem $\mathscr{A} \subset \mathscr{T}(X)$. The indiscrete topology on~$X$ is the coarsest member of~$\mathscr{T}(X)$, and thus it is a lower bound of any subsystem $\mathscr{A} \subset \mathscr{T}(X)$.
\erema

It is proven below that for a given set~$X$ the supremum and the infimum of every system of topologies $\mathscr{A} \subset \mathscr{T}(X)$ exist and are unique. In regard to Remark~\ref{rema top bounds} this is equivalent to the least upper bound property of the ordering~$\subset$ on~$\mathscr{T}(X)$. The supremum of~$\mathscr{A}$, which is the coarsest topology that is finer than every $\tp \in \mathscr{A}$, is determined in Corollary~\ref{coro gen top single space}. The infimum of~$\mathscr{A}$, which is the finest topology that is coarser than every $\tp \in \mathscr{A}$, is determined in Corollary~\ref{coro im top same set}.

It is often convenient to think of a topology as a system of sets that is, in a sense, "generated" by a subsystem of open sets. The first step in this direction is to determine a subsystem of open sets such that every open set can be written as a union of members of such a generating system. Such a generating system is called a base for the topology. Below we consider a yet smaller subsystem, called subbase.

\midvspace

\bdefi
\label{defi base top}
\index{Base!for a topology}
Given a topological space $(X,\tp)$, a system ${\mc B} \subset \tp$ is called {\bf base for}~$\tp$ if $\tp = \Theta ({\mc B})$, i.e.\ if the system of all unions of members of $\mc B$ is identical to~$\tp$. We also say that ${\mc B}$ {\bf generates}~$\tp$.
\edefi

For a given topology there generally exists more than one base generating it.

\midvspace

\bexam
Given a set $X$, the system ${\mc B} = \left\{ \left\{ x \right\} \, : \, x \in X \right\}$ of all singletons is a base for the discrete topology $\tp_{\mathrm{dis}}$ on~$X$. Every other base for $\tp_{\mathrm{dis}}$ contains~$\mc B$. 
\eexam

\bexam
\label{baseIndiscrete}
Given a set $X$, the only base for the indiscrete topology $\tp_{\mathrm{in}}$ on $X$ is $\tp_{\mathrm{in}}$ itself.
\eexam

\blemm
Given a topological space $(X,\tp)$, a system ${\mc B} \subset \tp$ is a base for~$\tp$ iff for every $U \in \tp$ and $x \in U$ there exists $B \in {\mc B}$ such that $x \in B \subset U$.
\elemm

\bproof
Exercise.
\eproof

The following Lemma provides a characterization of a system of subsets of~$X$ to be a base for some topology.

\midvspace

\blede
\label{lemm char base}
\index{Base!topological}
Let $X$ be a set and $\O \neq {\mc B} \subset {\mc P}(X)$. $\mc B$~is a base for some topology on~$X$ iff it satisfies all of the following conditions:
\benum
\item \label{lemm char base 0} $\O \in {\mc B}$
\item \label{lemm char base 1} $X = \bigcup {\mc B}$
\item \label{lemm char base 2} $\forall A, B \in \mc B \quad \exists \, {\mc C} \subset \mc B \quad {\mc C} \neq \O \; \wedge \;  A \cap B = \bigcup \mc C$
\eenum

In this case, $\mc B$ is also called a {\bf topological base on}~$X$. The topology generated by $\mc B$ is unique.
\elede

\bproof
Assume that (\ref{lemm char base 0}) to (\ref{lemm char base 2}) hold, and let $\tp = \Theta ({\mc B})$. It follows by Definition~\ref{defi top} and Lemma~\ref{lemm set equalities great} that $\tp$ is a topology on~$X$. Clearly $\mc B$ generates~$\tp$. The converse implication and the uniqueness of the topology generated by $\mc B$ are obvious.
\eproof

Notice that we have excluded the case ${\mc B} = \O$ in Lemma and Definition~\ref{lemm char base} only because property~(\ref{lemm char base 1}) has to be well-defined. If a system $\mc B$ satisfies property~(\ref{lemm char base 0}), this obviously implies ${\mc B} \neq \O$.

We may compare two topological bases on the same set~$X$ with each other similarly as we compare two topologies. We even refer in our definitions to the corresponding notions for the comparison of two topologies.

\midvspace

\bdefi
\label{defi base comp}
\index{Base!finer}
\index{Base!coarser}
\index{Base!strictly finer}
\index{Base!strictly coarser}
\index{Base!comparable}
Let $\mc A$ and $\mc B$ be two topological bases on a set~$X$. $\mc B$ is called {\bf finer}, {\bf coarser}, {\bf strictly finer}, {\bf strictly coarser than} $\mc B$ if the generated topologies $\Theta ({\mc A})$ and $\Theta ({\mc B})$ have the respective property. $\mc A$ and $\mc B$ are called {\bf comparable} if $\Theta ({\mc A})$ and $\Theta ({\mc B})$ are comparable.
\edefi

Notice that every topology on a set~$X$ is a base for itself. If $\tp_1$ and $\tp_2$ are two topologies on a set~$X$, $\tp_1$~is finer than~$\tp_2$ in the sense of Definition~\ref{defi top comp} iff $\tp_1$ is finer than~$\tp_2$ in the sense of Definition~\ref{defi base comp}, etc.\ That is Definitions~\ref{defi top comp} and~\ref{defi base comp} are consistent when referring to two topologies in both cases.

\midvspace

\bdefi
Given a set~$X$, the system of all topological bases on~$X$ is denoted by~$\mathscr{T}_B(X)$. 
\edefi

\brema
Given a set~$X$, we have $\mathscr{T}(X) \subset \mathscr{B}(X)$ and $\mathscr{B}(X) = \Theta^{-1} \left[ \mathscr{T}(X) \right]$. For every $\tp \in \mathscr{T}(X)$, the system of all bases for~$\tp$ is given by~$\Theta^{-1} \left\{ \tp \right\}$.
\erema

When we compare two topologies on~$X$, we use the ordering $\subset$ on~$\mathscr{T}(X)$. The direct comparison of topological bases requires another pre-ordering on~${\mc P}^2(X)$ whose properties are analysed in Lemma~\ref{lemm equivalences set operations}.

\midvspace

\brema
\label{rema comp base conditions}
Given a set~$X$, the pair $\left(\mathscr{B}(X),\subset_{\stackrel{}{\Theta}}\right)$ is a pre-ordered space with a reflexive relation. Let ${\mc A}, {\mc B} \in \mathscr{B}(X)$. The following statements are true:
\benum 
\item \label{rema comp base conditions 1} $\mc A$ is finer than $\mc B$ \; $\Longleftrightarrow \;\; {\mc B} \subset_{\stackrel{}{\Theta}} {\mc A}$
\item \label{rema comp base conditions 2} $\mc A$ is strictly finer than $\mc B$ \; $\Longleftrightarrow \;\; \left({\mc B} \subset_{\stackrel{}{\Theta}} {\mc A}\right) \; \wedge \; \neg \left({\mc A} \subset_{\stackrel{}{\Theta}} {\mc B}\right)$
\item \label{rema comp base conditions 3} $\mc A$ and $\mc B$ are comparable \; $\Longleftrightarrow \;\; \left({\mc B} \subset_{\stackrel{}{\Theta}} {\mc A}\right) \; \vee \; \left({\mc A} \subset_{\stackrel{}{\Theta}} {\mc B}\right)$
\eenum

\erema

\brema
Let $X$ be a set and ${\mc A}, {\mc B} \subset {\mc P}(X)$. By Lemma~\ref{lemm equivalences set operations}~(\ref{lemm equivalences set operations 2}) we have
\[
\Theta ({\mc A}) = \Theta ({\mc B}) \quad \Longleftrightarrow \quad \left({\mc B} \subset_{\stackrel{}{\Theta}} {\mc A}\right) \; \wedge \; \left({\mc A} \subset_{\stackrel{}{\Theta}} {\mc B}\right)
\]

If one side is true (and hence both sides are true) and ${\mc A}$ is a topological base, then also ${\mc B}$ is a topological base.

Notice that $\Theta ({\mc A}) = \Theta ({\mc B})$ need not imply ${\mc A} = {\mc B}$.
\erema

The following is a counterexample.

\midvspace

\bexam
Let $X$ be an infinite set. We may define systems
\[
{\mc B}_n = \left\{ B_{nk} \subset X \, : \, k \in \naturalnumbers,\, 1 \leq k \leq 2^n \right\} \quad (n \in \mathbb{N})
\]

with the following properties:
\benum
\item $B_{nk} \cap B_{nl} = \O \quad \left(n, k, l \in \naturalnumbers \,;\;\; 1 \leq k, l \leq 2^n \,;\;\; k \neq l \right)$
\item $\bigcup {\mc B}_n = X \quad \left(n \in \naturalnumbers \right)$
\item $B_{nk} = B_{(n+1)(2k-1)} \cup B_{(n+1)(2k)} \quad \left(n, k \in \naturalnumbers \,;\;\; 1 \leq k \leq 2^n \right)$
\eenum

Then each of the systems
\[
{\mc C} = \bigcup \left\{ {\mc B}_n : \, \mbox{$n$ is odd} \right\} \cup \left\{ \O \right\}, \quad {\mc D} = \bigcup \left\{ {\mc B}_n : \, \mbox{$n$ is even} \right\} \cup \left\{ \O \right\}
\]

is a topological base on~$X$, and $\Theta ({\mc C}) = \Theta ({\mc D})$. However, we have ${\mc C} \neq {\mc D}$.
\eexam

\bdefi
\index{Second countable}
We say that a topological space $(X,\tp)$ or the topology~$\tp$ is {\bf second countable} if there exists a countable base for~$\tp$.
\edefi

We define below what we mean by "first countable" topological space or topology. This definition is based on the notion of neighborhood to be introduced in Section~\ref{neighborhood}.

\midvspace

\bexam
Given an uncountable set $X$, the discrete topology on $X$ is not second countable.
\eexam

In the same way as for the system of all open sets, there is a possibility to "generate" the system of all closed sets from an appropriate subsystem that we call "base for the closed sets". It is a system such that every closed set is an intersection of its members.

\midvspace

\bdefi
\index{Base!for the closed sets}
Let $\xi = (X,\tp)$ be a topological space and ${\mc C}$ the system of all $\xi$-closed sets. A system ${\mc D} \subset {\mc C}$ is called {\bf base for}~$\mc C$ or {\bf base for the} $\xi$-{\bf closed sets} if ${\mc C} = \left\{ \, \bigcap {\mc G} \, : \, {\mc G} \subset {\mc D}, \; {\mc G} \neq \O \right\}$. When the set~$X$ is evident from the context, $\mc C$ is also called {\bf base for the} $\tp$-{\bf closed sets}. When the set~$X$ as well as the topology~$\tp$ are evident from the context, $\mc C$ is also called {\bf base for the closed sets}.
\edefi

The base for a topology and the base for the closed sets are related by complementation as follows.

\midvspace

\blemm
\label{baseclosed}
Let $(X,\tp)$ be a topological space, $\mc C$ the system of all closed sets, and ${\mc B} \subset {\mc P}(X)$. $\mc B$~is a base for~$\tp$ iff $\left\{ B^c : \, B \in {\mc B} \right\}$ is a base for~$\mc C$.
\elemm

\bproof
Exercise.
\eproof

The analogue of Lemma~\ref{lemm char base} for the system of all closed sets is stated in the following Lemma.

\midvspace

\blemm
\label{lemm base closed}
Let $X$ be a set and $\O \neq {\mc D} \subset {\mc P}(X)$. $\mc D$~is a base for the $\tp$-closed sets where $\tp$ is some topology on~$X$ iff it satisfies all of the following conditions:

\benum
\item \label{lemm base closed 1} $X \in {\mc D}$
\item \label{lemm base closed 2} $\O = \bigcap \mc D$
\item \label{lemm base closed 3} $\forall A, B \in \mc D \quad \exists \, {\mc E} \subset \mc D \quad {\mc E} \neq \O \; \wedge \; A \cup B = \bigcap \mc E$
\eenum

\elemm

\bproof
This follows by Lemmas~\ref{lemm char base} and~\ref{baseclosed}.
\eproof

Again, notice that the case ${\mc D} = \O$ is excluded in Lemma~\ref{lemm base closed} in order for the intersection in property~(\ref{lemm base closed 2}) to be well-defined. If a system $\mc D$ satisfies property~(\ref{lemm base closed 1}), this obviously implies ${\mc D} \neq \O$.

As announced before we are often able to "generate" a topological base from an appropriate smaller subsystem, called a "subbase".

\midvspace

\bdefi
\label{defi subbase}
\index{Subbase!for a topology}
Given a topological space $(X,\tp)$, a system ${\mc S} \subset \tp$ is called {\bf subbase for}~$\tp$ if the system $\Psi ({\mc S})$ is a base for $\tp$. We also say that ${\mc S}$ {\bf generates} $\tp$.
\edefi

Generally, for a given topological base, there exist more than one subbase generating it.

There is a simple criterion to probe whether a given system of subsets of~$X$ is a subbase for some topology on~$X$ as shown in the following Lemma. To this end the notion of "finite intersection property" is introduced.

\midvspace

\bdefi
\index{Finite intersection property}
Let $X$ be a set and ${\mc C} \subset {\mc P}(X)$. We say that $\mc C$ {\bf has the finite intersection property} if $\O \notin \Psi ({\mc C})$, i.e.\ for every ${\mc H} \sqsubset \mc C$ with ${\mc H} \neq \O$ we have $\bigcap {\mc H} \neq \O$.
\edefi

\blemm
\label{lemm crit subbase}
\index{Subbase!topological}
Let $X$ be a set and $\O \neq {\mc S} \subset {\mc P}(X)$. $\mc S$~is a subbase for some topology $\tp$ on~$X$ iff both of the following statements are true:
\benum
\item \label{lemm crit subbase 1} $\mc S$ does not have the finite intersection property.
\item \label{lemm crit subbase 2} $X = \bigcup {\mc S}$
\eenum

In this case, $\mc S$ is also called a {\bf topological subbase on}~$X$, and the topology generated by $\mc S$ is unique. Furthermore, $\tp$~is the coarsest topology on~$X$ that contains~$\mc S$.
\elemm

\bproof
The first claim is easy to verify (exercise).

Now let $\tp_1$ and $\tp_2$ be two topologies on $X$ and $\mc S$ a subbase both for $\tp_1$ and for $\tp_2$. By Definition~\ref{defi subbase} $\Psi ({\mc S})$ is a base both for $\tp_1$ and for $\tp_2$. Hence, $\tp_1 = \tp_2$ by Lemma~\ref{lemm char base}.

To prove the last claim, let $\mc S$ be a topological subbase. Further let $\tp_1$ be a topology on $X$ with ${\mc S} \subset \tp_1$. It follows that $\Theta \, \Psi \, ({\mc S}) \subset \tp_1$.
\eproof

Notice again that we have excluded the case ${\mc S} = \O$ in order for property~(\ref{lemm crit subbase 2}) to be well-defined. Of course, if a system $\mc S$ satisfies property~(\ref{lemm crit subbase 1}), this implies ${\mc S} \neq \O$.

Notice that for any system ${\mc S} \subset {\mc P}(X)$ the system ${\mc S} \cup \left\{ \O, X \right\}$ is a topological subbase. Alternatively, one could use the conventions that an intersection of an empty system of sets is identical to~$X$ and that a union of an empty system of sets is identical to~$\O$. In such a framework any arbitrary system of subsets of $X$ would be a subbase for some topology on~$X$. In this account the intersection or union of an empty system of sets is not defined.

We now state a criterion for the second countability of a topology.

\midvspace

\blemm
Let $(X,\tp)$ be a topological space and $\mc S$ a topological subbase for~$\tp$. If $\mc S$ is countable, then $\tp$ is second countable. 
\elemm

\bproof
This follows by Lemma~\ref{lemm fin count}.
\eproof

Next we introduce some notions for the comparison of two subbases on the same set.

\midvspace

\bdefi
\label{defi subbase comp}
\index{Subbase!finer}
\index{Subbase!coarser}
\index{Subbase!strictly finer}
\index{Subbase!strictly coarser}
\index{Subbase!comparable}
Let ${\mc S}_1$ and ${\mc S}_2$ be two topological subbases on a set~$X$. ${\mc S}_1$~is called {\bf finer}, {\bf coarser}, {\bf strictly finer}, {\bf strictly coarser than}~${\mc S}_2$ if the topological bases $\Psi ({\mc S}_1)$ and $\Psi ({\mc S}_2)$ have the respective property. ${\mc S}_1$~and ${\mc S}_2$ are called {\bf comparable} if $\Psi ({\mc S}_1)$ and $\Psi ({\mc S}_2)$ are comparable.
\edefi

Note that the notions defined in Definition~\ref{defi subbase comp} are well-defined since every topological subbase generates a unique topological base. Furthermore, every topological base is a subbase itself. When comparing two bases, no confusion can arise whether this comparison is done in the sense of Definition~\ref{defi base comp} or Definition~\ref{defi subbase comp} since Definition~\ref{defi subbase comp} refers to the corresponding notions for bases.

\midvspace

\bdefi
Given a set~$X$, the system of all topological subbases on~$X$ is denoted by~$\mathscr{S}(X)$.
\edefi

\brema
Given a set~$X$, we have $\mathscr{B}(X) \subset \mathscr{S}(X)$ and $\mathscr{S}(X) = \Psi^{-1} \left[ \mathscr{B}(X) \right]$. For every $\tp \in \mathscr{T}(X)$, the system of all subbases for~$\tp$ is given by~$\left(\Theta \, \Psi \right)^{-1} \left\{ \tp \right\}$.
\erema

Similarly to the case of topological bases in Remark~\ref{rema comp base conditions}, the comparison between two topological subbases may be expressed in terms of a pre-ordering as stated by the following remark.

\midvspace

\brema
\label{rema comp subbase conditions}
Given a set~$X$, the pair $\left(\mathscr{S}(X),\subset_{\stackrel{}{\Theta \Psi}}\right)$ is a pre-ordered space with a reflexive relation. Let ${\mc A}, {\mc B} \in \mathscr{S}(X)$. The following statements are true by Lemma~\ref{lemm equivalences set operations}~(\ref{lemm equivalences set operations 4}):

\benum 
\item \label{rema comp subbase conditions 1} $\mc A$ is finer than~$\mc B$ \; $\Longleftrightarrow \;\; {\mc B} \subset_{\stackrel{}{\Theta \Psi}} {\mc A}$
\item \label{rema comp subbase conditions 2} $\mc A$ is strictly finer than~$\mc B$ \; $\Longleftrightarrow \;\; \left({\mc B} \subset_{\stackrel{}{\Theta \Psi}} {\mc A}\right) \; \wedge \; \neg \left({\mc A} \subset_{\stackrel{}{\Theta \Psi}} {\mc B}\right)$
\item \label{rema comp subbase conditions 3} $\mc A$ and $\mc B$ are comparable \; $\Longleftrightarrow \;\; \left({\mc B} \subset_{\stackrel{}{\Theta \Psi}} {\mc A}\right) \; \vee \; \left({\mc A} \subset_{\stackrel{}{\Theta \Psi}} {\mc B}\right)$
\eenum

\erema

\brema
Let $X$ be a set and ${\mc A}, {\mc B} \subset {\mc P}(X)$. By Lemma~\ref{lemm equivalences set operations}~(\ref{lemm equivalences set operations 4}) we have
\[
\Theta \, \Psi ({\mc A}) = \Theta \, \Psi ({\mc B}) \quad \Longleftrightarrow \quad \left({\mc B} \subset_{\stackrel{}{\Theta \Psi}} {\mc A}\right) \; \wedge \; \left({\mc A} \subset_{\stackrel{}{\Theta \Psi}} {\mc B}\right)
\]

In this case, if ${\mc A}$ is a topological subbase, then also ${\mc B}$ is a topological subbase and the topologies generated by $\mc A$ and $\mc B$ are the same. Furthermore, notice that $\Theta \, \Psi ({\mc A}) = \Theta \, \Psi ({\mc B})$ need not imply ${\mc A} = {\mc B}$.
\erema

\bdefi
\index{Subbase!for the closed sets}
Let $\xi = (X,\tp)$ be a topological space and $\mc C$ the system of all $\xi$-closed sets. A system ${\mc S} \subset {\mc C}$ is called {\bf subbase for}~$\mc C$ or {\bf subbase for the  $\xi$-closed sets} if $\left\{ \, \bigcup {\mc H} \, : \, {\mc H} \sqsubset {\mc S}, \; {\mc H} \neq \O \right\}$ is a base for~$\mc C$. If the set~$X$ is evident from the context, $\mc C$~is also called {\bf subbase for the $\tp$-closed sets}. If the set~$X$ and the topology~$\tp$ are evident, $\mc C$~is called {\bf subbase for the closed sets}.
\edefi

\blemm
\label{subbase closed}
Let $\xi = (X,\tp)$ be a topological space, $\mc C$~the system of all $\xi$-closed sets, and ${\mc S} \subset {\mc P}(X)$. Then $\mc S$ is a subbase for~$\tp$ iff the system $\left\{ S^c : \, S \in {\mc S} \right\}$ is a subbase for~$\mc C$.
\elemm

\bproof
Exercise.
\eproof

\blemm
Let $X$ be a set and $\O \neq {\mc S} \subset {\mc P}(X)$. $\mc S$~is subbase for the $\tp$-closed sets, where $\tp$ is some topology on $X$, iff both of the following conditions are satisfied:
\benum
\item There exists ${\mc H} \sqsubset {\mc S}$ such that ${\mc H} \neq \O$ and $X = \bigcup {\mc H}$.
\item $\O = \bigcap \mc S$
\eenum

\elemm

\bproof
This follows from Lemmas~\ref{lemm crit subbase} and~\ref{subbase closed}.
\eproof

\section{Filters}
\label{filters}

So far, in this chapter, we have considered systems of sets whose members are subsets of a common set~$X$ as well as properties between members of such systems. Thereby we have not referred to any specific points of~$X$. In much of the following we consider particular points, values of maps at specific points, systems of sets that contain specific points etc.\ In particular this is the case in all circumstances where convergence is examined, a notion that is introduced in Chapter~\ref{convergence}. To this end, in this section the notion of filter is introduced as a fundamental concept. Note that the definitions and claims in this section do not involve any topology on~$X$.

\midvspace

\bdefi
\label{defi filter}
\index{Filter}
Given a set $X$, a system ${\mc F} \subset {\mc P}(X)$ is called a {\bf filter on}~$X$ if it has all of the following properties:
\benum
\item \label{defi filter 1} $\O \notin \mc F$
\item \label{defi filter 2} $X \in \mc F$
\item \label{defi filter 3} $\forall A, B \subset X \quad A, B \in \mc F \; \Longrightarrow \; A \cap B \in \mc F$
\item \label{defi filter 4} $\forall A, B \subset X \quad \left(A \in {\mc F}\right) \, \wedge \, \left(A \subset B\right) \; \Longrightarrow \; B \in \mc F$
\eenum

\edefi

Properties (\ref{defi filter 1}) and (\ref{defi filter 3}) in Definition~\ref{defi filter} imply that $\mc F$ has the finite intersection property.

\midvspace

\brema
\label{rema filter characterization}
Let $X$ be a set and ${\mc F} \subset {\mc P}(X)$. $\mc F$~is a filter on~$X$ iff all of the following conditions are satisfied:

\benum
\item \label{rema filter characterization 1} ${\mc F} \neq \O$
\item \label{rema filter characterization 2} ${\mc F} \neq {\mc P}(X)$
\item \label{rema filter characterization 3} $\Psi ({\mc F}) = {\mc F}$
\item \label{rema filter characterization 4} $\Phi ({\mc F}) = {\mc F}$
\eenum

That is, the filters on~$X$ are precisely those systems that are fixed points of both $\Psi$ and~$\Phi$ and additionally are non-trivial in the sense of (\ref{rema filter characterization 1}) and~(\ref{rema filter characterization 2}).
\erema

When comparing two filters on the same set we use the same notions as for topologies.

\midvspace

\bdefi
\label{defi filter finer}
\index{Filter!finer}
\index{Filter!coarser}
\index{Filter!strictly finer}
\index{Filter!strictly coarser}
\index{Filter!comparable}
Let ${\mc F}_1$ and ${\mc F}_2$ be two filters on a set~$X$. If ${\mc F}_2 \subset {\mc F}_1$, then ${\mc F}_1$ is called {\bf finer than}~${\mc F}_2$ and ${\mc F}_2$ is called {\bf coarser than}~${\mc F}_1$. If ${\mc F}_2 \subset {\mc F}_1$ and ${\mc F}_1 \neq {\mc F}_2$, then ${\mc F}_1$ is called {\bf strictly finer than}~${\mc F}_2$ and ${\mc F}_2$ is called {\bf strictly coarser than}~${\mc F}_1$. If ${\mc F}_1 \subset {\mc F}_2$ or ${\mc F}_2 \subset {\mc F}_1$, the filters are called {\bf comparable}.
\edefi

\bdefi
Given a set~$X$, the system of all filters on~$X$ is denoted by~$\mathscr{F}(X)$.
\edefi

\blede
Let $X$ be a set, $\mathscr{A} \subset \mathscr{F}(X)$, and ${\mc F} \in \mathscr{A}$. The pair $\left(\mathscr{A}, \subset \right)$ is an ordered space in the sense~"$\leq$". $\mc F$ is called {\bf finest} ({\bf coarsest}) filter of~$\mathscr{A}$ if it is a maximum (minimum) of~$\mathscr{A}$. $\mathscr{A}$~has at most one finest and at most one coarsest filter.
\elede

\bproof
Exercise.
\eproof

\bdefi
\index{Ultrafilter}
Given a filter~$\mc F$ on a set~$X$, $\mc F$ is called {\bf ultrafilter} if it is a weak maximum of~$\mathscr{F}(X)$ with respect to the ordering~$\subset$, i.e.\ if there is no filter on~$X$ that is strictly finer than~$\mc F$.
\edefi

Given a specific filter it is often convenient to consider only a particular subsystem, called filter base, instead of the entire filter. A filter base is defined by the requirement that every member of the filter contains a member of the filter base. More formally we have the following definition.

\midvspace

\bdefi
\label{defi filter base}
\index{Filter base}
\index{Base!for a filter}
Given a filter~$\mc F$ on a set~$X$, a subsystem ${\mc B} \subset {\mc F}$ is called {\bf filter base for}~$\mc F$ if ${\mc F} = \Phi ({\mc B})$. We also say that $\mc B$ {\bf generates}~$\mc F$.
\edefi

Notice that this definition is similar to the definition of a base for a topology, see Definition~\ref{defi base top}. The two are related by the concept of neighborhood system, which is described below in Section~\ref{neighborhood}.

\midvspace

\blede
\label{lemm characterization filter base}
\index{Filter base}
Let $X$ be a set and ${\mc B} \subset {\mc P}(X)$. $\mc B$~is a filter base for some filter on~$X$ iff all of the following conditions are satisfied:
\benum
\item \label{lemm characterization filter base 1} $\O \notin \mc B$
\item \label{lemm characterization filter base 2} ${\mc B} \neq \O$
\item \label{lemm characterization filter base 3} $\forall A, B \in {\mc B} \quad \exists C \in {\mc B} \quad C \subset A \cap B$
\eenum

In this case, $\mc B$ is also called {\bf filter base on}~$X$. The filter generated by~$\mc B$ is unique.
\elede

\bproof
Assume that (\ref{lemm characterization filter base 1}) to (\ref{lemm characterization filter base 3}) hold. Then clearly ${\mc F} =\Phi ({\mc B})$ is a filter on~$X$. The converse implication and the uniqueness of the generated filter are obvious.
\eproof

Notice that property~(\ref{lemm characterization filter base 3}) in Lemma and Definition~\ref{lemm characterization filter base} is equivalent to $\Psi ({\mc B}) \subset \Phi({\mc B})$.

We now define appropriate notions for the comparison of two filter bases.

\midvspace

\bdefi
\label{defi filter base finer}
\index{Filter base!finer}
\index{Filter base!coarser}
\index{Filter base!strictly finer}
\index{Filter base!strictly coarser}
\index{Filter base!comparable}
Let $X$ be a set, and $\mc A$ and $\mc B$ two filter bases on~$X$. $\mc A$~is called {\bf finer}, {\bf coarser}, {\bf strictly finer}, {\bf strictly coarser than}~$\mc B$ if $\Phi ({\mc A})$ and $\Phi ({\mc B})$ have the respective property. $\mc A$~and $\mc B$ are called {\bf comparable} if $\Phi ({\mc A})$ and $\Phi ({\mc B})$ are comparable.
\edefi

This definition relies on the fact that the filter generated by a filter base is unique. Also the fact that every filter is a filter base for itself is taken into account by referring to the properties of filters in the definition.

\midvspace

\bdefi
Given a set~$X$, the system of all filter bases on~$X$ is denoted by~$\mathscr{F}_{\!\stackrel{}{B}}(X)$.
\edefi

\brema
Given a set~$X$, we have $\mathscr{F}(X) \subset \mathscr{F}_{\!\stackrel{}{B}}(X)$ and $\mathscr{F}_{\!\stackrel{}{B}}(X) = \Phi^{-1} \left[ \mathscr{F}(X) \right]$. For every ${\mc F} \in \mathscr{F}(X)$, the system of all filter bases for~$\mc F$ is given by~$\Phi^{-1} \left\{ {\mc F} \right\}$.
\erema

The comparison of two filter bases on a set~$X$ may be expressed by the pre-ordering $\subset_{\stackrel{}{\Phi}}$.

\midvspace

\brema
\label{coro filter bases comparison}
Given a set~$X$, the pair $\left(\mathscr{F}_{\!\stackrel{}{B}}(X),\subset_{\stackrel{}{\Phi}}\right)$ is a pre-ordered space with a reflexive relation. For every ${\mc A}, {\mc B} \in \mathscr{F}_{\!\stackrel{}{B}}(X)$, the following statements are true by Lemma~\ref{lemm equivalences set operations}~(\ref{lemm equivalences set operations 3})

\benum
\item \label{coro filter bases comparison 1} $\mc A$ is finer than $\mc B$ \; $\Longleftrightarrow \;\; {\mc B} \subset_{\stackrel{}{\Phi}} {\mc A}$
\item \label{coro filter bases comparison 2} $\mc A$ is strictly finer than $\mc B$ \; $\Longleftrightarrow \;\; \left({\mc B} \subset_{\stackrel{}{\Phi}} {\mc A}\right) \; \wedge \; \neg \left({\mc A} \subset_{\stackrel{}{\Phi}} {\mc B}\right)$
\item \label{coro filter bases comparison 3} $\mc A$ and $\mc B$ are comparable \; $\Longleftrightarrow \;\; \left({\mc B} \subset_{\stackrel{}{\Phi}} {\mc A}\right) \; \vee \; \left({\mc A} \subset_{\stackrel{}{\Phi}} {\mc B}\right)$
\eenum

\erema

\brema
Let $X$ be a set and ${\mc A}, {\mc B} \subset {\mc P}(X)$. By Lemma~\ref{lemm equivalences set operations}~(\ref{lemm equivalences set operations 3}) we have $\Phi ({\mc A}) = \Phi ({\mc B})$ iff ${\mc A} \subset_{\stackrel{}{\Phi}} {\mc B}$ and ${\mc B} \subset_{\stackrel{}{\Phi}} {\mc A}$. In this case, if $\mc A$ is a filter base on~$X$, then $\mc B$ is also a filter base on~$X$ and the two filters generated by~$\mc A$ and~$\mc B$ are the same. Notice however that $\Phi ({\mc A}) = \Phi ({\mc B})$ need not imply ${\mc A} = {\mc B}$.
\erema

\bdefi
\index{Cluster point!filter}
\index{Filter!fixed}
\index{Filter!free}
\index{Cluster point!filter base}
\index{Filter base!fixed}
\index{Filter!free}
Given a filter $\mc F$ on a set~$X$, a point $x \in X$ is called a {\bf cluster point of~$\mc F$} if $x \in \bigcap {\mc F}$. $\mc F$~is called {\bf free} if it has no cluster points, otherwise it is called {\bf fixed}. Similarly, for a filter base $\mc B$ on~$X$, a point $x \in X$ is called a {\bf cluster point of}~$\mc B$ if it is a cluster point of~$\Phi ({\mc B})$. $\mc B$~is called {\bf free} if it has no cluster points, otherwise it is called {\bf fixed}.
\edefi

\brema
Given a filter base $\mc B$ on a set~$X$ and a point $x \in X$, $x$~is a cluster point of~$\mc B$ iff $x \in \bigcap {\mc B}$.
\erema

\bexam
\label{exfrech}
Given a set $X$, the system ${\mc B} = \left\{ \; ]x,\infty[ \;\, : \, x \in {\mathbb R}\right\}$ is a free filter base on~$\mathbb R$.
\eexam

\bexam
For every $r \in \; ]0,\infty[ \;$, let $B_r = \left\{ (x,y) \in {\mathbb R}^2 \, : \, x^2 + y^2 \leq r \right\}$. The system ${\mc B} = \left\{ B_r \, : \, r \in \; ]0,\infty[ \; \right\}$ is a fixed filter base on~${\mathbb R}^2$.
\eexam

\bexam
\index{Filter!indiscrete}
Given a set $X$, the system ${\mc F} = \left\{ X \right\}$ is a filter. It is called the {\bf indiscrete filter on~$X$}.
\eexam

\blemm
\index{Filter!discrete}
Given a set $X$ and $A \subset X$ with $A \neq \O$, the system ${\mc F} = \left\{ F \subset X \, : \, A \subset F \right\}$ is a fixed filter on~$X$. The system $\left\{ A \right\}$ is a filter base for~$\mc F$. $\mc F$~is an ultrafilter iff $A = \left\{ x \right\}$ for some $x \in X$. In this case it is called the {\bf discrete filter at~$x$}.
\elemm

\bproof
Exercise.
\eproof

\blede
\label{lemm two filter bases intersection}
\index{Intersection!of two filter bases}
\index{Filter base!intersection}
Let $\mc A$ and $\mc B$ be two filter bases on a set~$X$, where $A \cap B \neq \O$ for every $A \in {\mc A}$, $B \in {\mc B}$. The system ${\mc C} = \left\{ A \cap B \, : \, A \in {\mc A},\, B \in {\mc B} \right\}$ is a filter base. It is called the {\bf filter-base intersection of} $\mc A$ {\bf and}~$\mc B$. $\mc C$~is a supremum of $\left\{ {\mc A}, {\mc B} \right\}$ in the pre-ordered space $\left(\mathscr{F}_{\!\stackrel{}{B}}(X),\subset_{\stackrel{}{\Phi}}\right)$. Equivalently, the filter $\Phi ({\mc C})$ is the supremum of $\left\{ \Phi ({\mc A}), \Phi ({\mc B}) \right\}$ in the ordered space $\left(\mathscr{F}(X),\subset\right)$, i.e.\ it is the coarsest filter on~$X$ that is finer than $\Phi ({\mc A})$ and~$\Phi ({\mc B})$.
\elede

\bproof
Let $A_1, A_2 \in {\mc A}$ and $B_1, B_2 \in {\mc B}$. We may choose $A_3 \in {\mc A}$ such that $A_3 \subset A_1 \cap A_2$, and $B_3 \in {\mc B}$ such that $B_3 \subset B_1 \cap B_2$. It follows that
\[
A_3 \cap B_3 \; \subset \; A_1 \cap B_1 \cap A_2 \cap B_2
\]

Thus $\mc C$ is a filter base. We clearly have ${\mc A} \subset_{\stackrel{}{\Phi}} {\mc C}$ and ${\mc B} \subset_{\stackrel{}{\Phi}} {\mc C}$. Assume that $\mc D$~is a filter base on~$X$ with ${\mc A} \subset_{\stackrel{}{\Phi}} {\mc D}$ and ${\mc B} \subset_{\stackrel{}{\Phi}} {\mc D}$. We show that ${\mc C} \subset_{\stackrel{}{\Phi}} {\mc D}$. Let $A \in {\mc A}$ and $B \in {\mc B}$. We may choose $A', B' \in {\mc D}$ such that $A' \subset A$ and $B' \subset B$. There exists $D \in {\mc D}$ such that $D \subset A' \cap B'$. Thus we have $D \subset A \cap B$. This shows that $\mc C$ is a supremum. It clearly follows that $\Phi ({\mc C})$ is a supremum, which is unique since $\subset$ is an ordering.
\eproof

\blemm
\label{lemm filter base finite intersection}
\index{Intersection!finite, of filter bases}
\index{Filter base!intersection}
Let $X$ be a set and ${\mc B}_i \in \mathscr{F}_{\!\stackrel{}{B}}(X)$ ($i \in I$) where $I$ an index set. Furthermore we define the system
\[
{\mc B} = \left\{ \, {\textstyle \bigcap}_{j \in J} \, B_j \, : \, J \sqsubset I, \; J \neq \O, \; B_j \in {\mc B}_j \; (j \in J) \right\}
\]

If $\O \notin {\mc B}$, then ${\mc B} \in \mathscr{F}_{\!\stackrel{}{B}}(X)$. In this case, $\mc B$ is a supremum of~$\left\{ {\mc B}_i \, : \, i \in I \right\}$ in the pre-ordered space $\left(\mathscr{F}_{\!\stackrel{}{B}}(X),\subset_{\stackrel{}{\Phi}}\right)$. Equivalently, the filter $\Phi ({\mc B})$ is the supremum of $\left\{ \Phi ({\mc B}_i) \, : \, i \in I \right\}$ in the ordered space $\left(\mathscr{F}(X),\subset\right)$, i.e.\ it is the coarsest filter on~$X$ that is finer than $\Phi ({\mc B}_i)$ for every~$i \in I$.
\elemm

\bproof
Assume that $\O \notin {\mc B}$. It follows by Lemma~\ref{lemm characterization filter base} that $\mc B$ is a filter base on $X$. Clearly, $\mc B$ is finer than ${\mc B}_i$ for every $i \in I$. It remains to show that every filter base on~$X$ that is finer than ${\mc B}_i$ for every $i \in I$ is also finer than~$\mc B$. Assume that $\mc A$ is a filter base on $X$ that is finer than ${\mc B}_i$ for every $i \in I$, and let $B \in {\mc B}$. We may choose $J \sqsubset I$ with $J \neq \O$ and $B_j \in {\mc B}_j$ ($j \in J$) such that $B = \bigcap_{j \in J} B_j$. For every $j \in J$ we may choose $A_j \in {\mc A}$ such that $A_j \subset B_j$ by assumption. There exists $A \in {\mc A}$ such that $A \subset \bigcap_{j \in J} A_j$. Hence $A \subset B$. This shows that $\mc B$ is a supremum of~$\left\{ {\mc B}_i \, : \, i \in I \right\}$. It clearly follows that $\Phi ({\mc B})$ is a supremum of $\left\{ \Phi ({\mc B}_i) \, : \, i \in I \right\}$, which is unique since $\subset$ is an ordering.
\eproof

Notice that the filter base~$\mc B$ defined in Lemma~\ref{lemm filter base finite intersection} is not a direct generalization of Lemma and Definition~\ref{lemm two filter bases intersection}, as the index set $J$ may be a singleton. However, it is easy to see that the generated filters are the same if $J \sim 2$.

\midvspace

\bdefi
\index{Ultrafilter base}
Let ${\mc B}$ be a filter base on a set~$X$. $\mc B$~is called {\bf ultrafilter base} if $\Phi ({\mc B})$ is an ultrafilter.
\edefi

The following are a few characterizations of an ultrafilter base.

\midvspace

\blemm
\label{lemm ultra base filter}
Let $\mc B$ be a filter base on a set~$X$. The following statements are equivalent.

\benum
\item \label{lemm ultra base filter 1} ${\mc B}$ is an ultrafilter base.
\item \label{lemm ultra base filter 2} $\mc B$ is a weak maximum of $\mathscr{F}_{\!\stackrel{}{B}}(X)$ with respect to the pre-ordering~$\subset_{\stackrel{}{\Phi}}$, i.e.\ there is no filter base on~$X$ that is strictly finer than~$\mc B$.
\item \label{lemm ultra base filter 3} $\forall A \subset X \quad A \in \Phi ({\mc B}) \; \vee \; A^c \in \Phi ({\mc B})$
\item \label{lemm ultra base filter 4} $\forall A \subset X \quad \left\{ A \right\} \subset_{\stackrel{}{\Phi}} {\mc B} \; \vee \ \left\{ A^c \right\} \subset_{\stackrel{}{\Phi}} {\mc B}$
\eenum

\elemm

\bproof
(\ref{lemm ultra base filter 1}) and (\ref{lemm ultra base filter 2}) are clearly equivalent.

The equivalence of (\ref{lemm ultra base filter 3}) and (\ref{lemm ultra base filter 4}) is obvious as well.

To show that (\ref{lemm ultra base filter 1}) implies (\ref{lemm ultra base filter 4}), let $A \subset X$ and assume that neither $\left\{ A \right\} \subset_{\stackrel{}{\Phi}} {\mc B}$ nor $\left\{ A^c \right\} \subset_{\stackrel{}{\Phi}} {\mc B}$ holds. It follows by Lemma and Definition~\ref{lemm two filter bases intersection} that each of the systems
\[
{\mc A} = \left\{ B \cap A \, : \, B \in {\mc B} \right\}, \quad \quad {\mc C} = \left\{ B \cap A^c \, : \, B \in {\mc B} \right\}
\]

is a filter base and that $\mc A$ and $\mc C$ are both finer than~$\mc B$. We clearly have $\Phi ({\mc A}) \neq \Phi ({\mc C})$. Hence at least one of the filter bases is strictly finer than~$\mc B$, which is a contradiction to the fact that $\mc B$ is an ultrafilter base.

To show that (\ref{lemm ultra base filter 4}) implies (\ref{lemm ultra base filter 1}), let $\mc A$ be a filter base on~$X$ that is finer than~$\mc B$, and $A \in {\mc A}$. Since $\left\{ A^c \right\} \subset_{\stackrel{}{\Phi}} {\mc B}$ clearly does not hold, we have $\left\{ A \right\} \subset_{\stackrel{}{\Phi}} {\mc B}$. Therefore $\mc B$ is finer than~$\mc A$, and hence $\mc B$ is an ultrafilter base.
\eproof

Note that Lemma~\ref{lemm ultra base filter} says that ultrafilters are precisely those filters that contain either $A$ or $A^c$ for every $A \subset X$.

\midvspace

\blemm
Let $X$ be a set. $\mc F$~is a fixed ultrafilter on~$X$ iff there is a point $x \in X$ such that ${\mc F} = \left\{ F \subset X \, : \, x \in F \right\}$.
\elemm

\bproof
If $\mc F$ is a fixed ultrafilter, it clearly has a unique cluster point by Lemma \ref{lemm ultra base filter} (\ref{lemm ultra base filter 3}), say~$x$. Let $F \subset X$ with $x \in F$. Since $F^c \notin {\mc F}$, we have $F \in {\mc F}$ by Lemma~\ref{lemm ultra base filter}~(\ref{lemm ultra base filter 3}). This shows that ${\mc F} = \left\{ F \subset X \, : \, x \in F \right\}$.

The converse is obvious.
\eproof

\btheo
\label{theo existence ultrafilter base}
\index{Ultrafilter base!existence}
For every filter base ${\mc B}$ on a set~$X$ there exists an ultrafilter base that is finer than~$\mc B$.
\etheo

\bproof
Let $M = \left\{ {\mc B}_i \, : \, i \in I \right\} \subset \mathscr{F}_{\!\stackrel{}{B}}(X)$ be the set of all filter bases on~$X$ that are finer than~$\mc B$ where $I$ is an index set. The relation $\subset_{\stackrel{}{\Phi}}$ is a pre-ordering on~$M$. Let $L = \left\{ {\mc B}_j \, : \, j \in J \right\}$ ($J \subset I$) be a chain and
\[
{\mc A} = \Big\{ \, {\textstyle \bigcap}_{k \in K} \, B_k \, : \, K \sqsubset J, \; K \neq \O, \; B_k \in {\mc B}_k \; (k \in K) \Big\}
\]

Since $\O \notin {\mc A}$, $\mc A$ is a filter base that is finer than ${\mc B}_j$ for every $j \in J$ by Lemma~\ref{lemm filter base finite intersection}. Thus $\mc A$ is an upper bound of~$L$. Let $\mc C$ be a weak maximum of~$M$ according to Theorem~\ref{theo Zorn}. To prove that $\mc C$ is an ultrafilter base, assume that $\mc D$ is a filter base on~$X$ with ${\mc C} \subset_{\stackrel{}{\Phi}} {\mc D}$. Since ${\mc B} \subset_{\stackrel{}{\Phi}} {\mc C}$, we have ${\mc B} \subset_{\stackrel{}{\Phi}} {\mc D}$, and thus ${\mc D} \in M$. Since $\mc C$ is a weak maximum of~$M$, ${\mc C} \subset_{\stackrel{}{\Phi}} {\mc D}$ implies ${\mc D} \subset_{\stackrel{}{\Phi}} {\mc C}$.
\eproof

\bcoro
\index{Ultrafilter!existence}
For every filter $\mc F$ on a set~$X$, there is an ultrafilter~$\mc G$ that is finer than~$\mc F$.
\ecoro

\bproof
By Theorem~\ref{theo existence ultrafilter base} we may choose an ultrafilter base $\mc B$ such that $\mc B$ is finer than~$\mc F$. Then $\Phi ({\mc B})$ is an ultrafilter that is finer than~$\mc F$.
\eproof

\section{Neighborhoods}
\label{neighborhood}

We briefly outline the major notions defined in this Section. Given a topological space $(X,\tp)$, a neighborhood $U$ of a point $x$ is defined as a---not necessarily open---set for which there exists an open set~$V$ such that $x \in V \subset U$. The system of all neighborhoods of a particular point~$x$ is called the neighborhood system of~$x$. One may also consider the ensemble of neighborhood systems for every $x \in X$, which may be called the neighborhood system of the topological space. The natural way to describe this is to define a structure relation on $X$ (cf.\ Definition~\ref{defi structure relation}), that contains the pair $(x,U)$ for every point~$x$ and every neighborhood~$U$ of~$x$. More formally this leads to Definitions~\ref{defi neighborhood system} and~\ref{defi neighborhood} below. In order to describe the neighborhood system it would be possible to define a function on~$X$ such that for each point $x \in X$ the value $f(x)$ is the neighborhood system of~$x$. The domain of such a function is $X$ and its range is a subset of~${\mc P}^2(X)$. However, this turns out to be notationally disadvantageous when choosing subsets of the neighborhood system of single points.

\midvspace

\bdefi
\label{defi neighborhood system}
\index{Neighborhood system}
\index{Neighborhood system!open}
\index{Neighborhood system!closed}
Let $\xi = (X,\tp)$ be a topological space and $\mc C$ the system of all $\xi$-closed sets. The structure relation ${\mc N}_{\xi}$, or short $\mc N$, defined by
\[
(x,U) \in {\mc N} \quad \Longleftrightarrow \quad \exists V \in \tp \quad x \in V \subset U
\]

is called the {\bf neighborhood system of}~$\xi$. The structure relation
\[
{\mc N}^{\, \mathrm{open}}_{\xi} = \, {\mc N}_{\xi} \, \cap \,  \left( X \!\times \tp \right)
\]

is called the {\bf open neighborhood system of}~$\xi$, also denoted by~${\mc N}^{\, \mathrm{open}}_{\phantom{1}}$. The structure relation
\[
{\mc N}^{\, \mathrm{closed}}_{\xi} = \, {\mc N}_{\xi} \, \cap \, \left( X \!\times {\mc C} \right)
\]

is called the {\bf closed neighborhood system of}~$\xi$, also denoted by~${\mc N}^{\, \mathrm{closed}}_{\phantom{1}}$.
\edefi

\brema
\label{rema properties nh system}
Given a topological space $(X,\tp)$, the following statements hold:

\benum
\item \label{rema properties nh system 1} $\forall x \in X \quad \forall U \subset X \quad (x,U) \in {\mc N}^{\, \mathrm{open}}_{\phantom{1}} \; \Longleftrightarrow \; U \in \tp(x)$
\item \label{rema properties nh system 2} $\Phi' \big({\mc N}^{\, \mathrm{open}}_{\phantom{1}}\big) = \Phi' ({\mc N}) = {\mc N}$
\eenum

Note that we use the notation of Definition~\ref{defi set system point} in~(\ref{rema properties nh system 1}) and the map $\Phi'$ defined in Definition~\ref{defi phi prime} in~(\ref{rema properties nh system 2}).
\erema

We now introduce the neighborhood system of an arbitrary set $A \subset X$. The special case in which $A$ is a singleton leads to the definition of the neighborhood system of a point. Remember that for a relation $R$ on~$X$ we have defined
\begin{eqnarray*}
R \left\{ x \right\} \!\!\! & = & \!\! \left\{ y \in X \, : \, (x,y) \in R \right\},\\
R \langle A \rangle \!\!\! & = & \!\! \big\{ y \in X \, : \, \forall x \in A \;\; (x,y) \in R \big\}
\end{eqnarray*}

\midvspace

\bdefi
\label{defi neighborhood}
\index{Neighborhood}
\index{Neighborhood!open}
\index{Neighborhood!closed}
Let $\mc N$ be the neighborhood system of a topological space~$(X,\tp)$.

\benum
\item For every $A \subset X$ the system ${\mc N} \langle A \rangle$ is called the {\bf neighborhood system of}~$A$, and every member is called a {\bf neighborhood of}~$A$.
\item For every $A \subset X$, the systems ${\mc N}^{\, \mathrm{open}}_{\phantom{1}} \langle A \rangle$ and ${\mc N}^{\, \mathrm{closed}}_{\phantom{1}} \langle A \rangle$ are called the {\bf open} and the {\bf closed neighborhood systems of}~$A$, respectively, and every member is called an {\bf open (closed) neighborhood of}~$A$.
\item For every $x \in X$, the system ${\mc N} \! \left\{ x \right\} = {\mc N} \langle \left\{ x \right\} \rangle$ is called the {\bf neighborhood system of}~$x$. A member $U \in {\mc N} \! \left\{ x \right\}$ is called {\bf neighborhood of}~$x$.
\item For every $x \in X$, the systems
\[
{\mc N}^{\, \mathrm{open}}_{\phantom{1}} \! \left\{ x \right\} = {\mc N}^{\, \mathrm{open}}_{\phantom{1}} \langle \left\{ x \right\} \rangle, \quad \quad {\mc N}^{\, \mathrm{closed}}_{\phantom{1}} \! \left\{ x \right\} = {\mc N}^{\, \mathrm{closed}}_{\phantom{1}} \langle \left\{ x \right\} \rangle
\]

are called the {\bf open} and {\bf closed neighborhood systems of}~$x$, respectively. Their members are respectively called the {\bf open} and {\bf closed neighborhoods of}~$x$.
\eenum

\edefi

Notice that for every open $A \subset X$, $A$ is a neighborhoood of itself. Moreover, we have ${\mc N} \langle \O \rangle = {\mc P}(X)$, i.e.\ every subset of~$X$ is a neighborhood of the empty set. We state a few more consequences of the preceding definitions in the following Lemma.

\midvspace

\blemm
\label{lemm characterization neighborhood}
Let $(X,\tp)$ be a topological space and $U, V \subset X$. The following statements hold:

\benum
\item \label{lemm characterization neighborhood 1} $U \in \tp \quad \Longleftrightarrow \quad \forall x \in U \quad U \in {\mc N} \! \left\{ x \right\}$
\item \label{lemm characterization neighborhood 2} $U \in \tp \quad \Longleftrightarrow \quad \forall A \subset U \quad U \in {\mc N} \langle A \rangle$
\item \label{lemm characterization neighborhood 3} $V \in {\mc N} \langle U \rangle \quad \Longleftrightarrow \quad \exists W \in \tp \quad U \subset W \subset V$
\eenum

\elemm

\bproof
To prove~(\ref{lemm characterization neighborhood 1}) assume $U \in {\mc N} \! \left\{ x \right\}$ holds for every $x \in U$. For each $x \in U$ we may choose $V_x \in \tp$ such that $x \in V_x \subset U$. Hence $U = \bigcup_{x \in U} V_x \in \tp$. The converse is obvious.

(\ref{lemm characterization neighborhood 2}) is a direct consequence of (\ref{lemm characterization neighborhood 1}).

To prove (\ref{lemm characterization neighborhood 3}) assume $V \in {\mc N} \langle U \rangle$. Then $V \in {\mc N} \! \left\{ x \right\}$ for every $x \in U$. For each $x \in U$, we may choose $W_x \in \tp$ such that $x \in W_x \subset V$. Therefore $U \subset W \subset V$ where $W = \bigcup_{x \in U} W_x$. The converse is obvious. 
\eproof

Our definitions of neighborhood and neighborhood system are based on the notion of topological space. In the following Lemma we list properties of the neighborhood system. We then show that for a given structure relation with these properties there is a unique topology such that the neighborhood system is this structure relation.

\midvspace

\blemm
\label{nhsystem}
The neighborhood system $\mc N$ of a topological space $(X,\tp)$ has the following properties:

\benum
\item \label{nhsystem0} $\forall \, x \in X \quad {\mc N} \! \left\{ x \right\} \neq \O$
\item \label{nhsystem1} $\forall \, x \in X \quad \forall \, U \!\in {\mc N} \! \left\{ x \right\} \quad x \in U$
\item \label{nhsystem2} $\forall \, x \in X \quad \forall \, U \!\in {\mc N} \! \left\{ x \right\} \quad \forall \, V \!\supset U \quad V \!\in {\mc N} \! \left\{ x \right\}$
\item \label{nhsystem3} $\forall \, x \in X \quad \forall \, U, V \!\in {\mc N} \! \left\{ x \right\} \quad U \cap V \in {\mc N} \! \left\{ x \right\}$
\item \label{nhsystem4} $\forall \, x \in X \quad \forall \, U \!\in {\mc N} \! \left\{ x \right\} \quad \exists \, V \!\in {\mc N} \! \left\{ x \right\} \quad \forall \, y \in V \quad U \!\in {\mc N} \! \left\{ y \right\}$
\eenum

\elemm

\bproof
(\ref{nhsystem0}) to (\ref{nhsystem3}) are obvious. To show (\ref{nhsystem4}), let $x \in X$ and $U \!\in {\mc N} \! \left\{ x \right\}$. We may choose $V \!\in \tp$ such that $x \in V \!\subset U$. Then, for every $y \in V$, we have $U \!\in {\mc N} \! \left\{ y \right\}$.
\eproof

\blemm
\label{lemma top from neighborhood}
Let $X$ be a set and ${\mc N} \subset X \times {\mc P}(X)$ such that properties (\ref{nhsystem0}) to~(\ref{nhsystem4}) in Lemma~\ref{nhsystem} are satisfied. The system $\tp = \big\{ U \subset X \, : \, \forall x \in U \;\; U \!\in {\mc N} \! \left\{ x \right\} \!\big\}$ is the unique topology on~$X$ such that $\mc N$ is the neighborhood system of~$(X,\tp)$.
\elemm

\bproof
Notice that $\tp$ is a topology on~$X$ by properties (\ref{nhsystem0}), (\ref{nhsystem2}), and~(\ref{nhsystem3}). Let ${\mc N}_{\xi}$ be the neighborhood system of $\xi = (X,\tp)$. We show that ${\mc N}_{\xi} = {\mc N}$. Fix $x \in X$. 

First assume that $U \!\in {\mc N}_{\xi} \!\left\{ x \right\}$. There exists $V \!\in \tp$ such that $x \in V \!\subset U$. Therefore $V \!\in {\mc N} \! \left\{ x \right\}$ by definition of~$\tp$, and thus $U \!\in {\mc N} \! \left\{ x \right\}$.

Now assume $U \!\in {\mc N} \! \left\{ x \right\}$. We may define $V = \big\{ y \in X \, : \, U \in {\mc N} \! \left\{ y \right\} \!\big\}$. Fix $y \in V$. By property~(\ref{nhsystem4}) there is $W \!\in {\mc N} \! \left\{ y \right\}$ such that $U \!\in {\mc N} \! \left\{ z \right\}$ for every $z \in W$. Thus $W \!\subset V$, and therefore $V \!\in {\mc N} \! \left\{ y \right\}$ by property~(\ref{nhsystem2}). It follows that $V \!\in \tp$. Since $x \in V \!\subset U$ by property~(\ref{nhsystem1}), we have $U \!\in {\mc N}_{\xi} \!\left\{ x \right\}$.

The uniqueness of $\tp$ follows by Lemma~\ref{lemm characterization neighborhood}~(\ref{lemm characterization neighborhood 1}).
\eproof

In the following Lemma we list several statements that are satisfied by the neighborhoods of subsets of~$X$---in contrast to the neighborhoods of single points of~$X$ considered in Lemma~\ref{nhsystem}.

\midvspace

\blemm
\label{nhsystem sets}
The neighborhood system $\mc N$ of a topological space $(X,\tp)$ has the following properties:

\benum
\item \label{nhsystem sets 1} $\forall A \subset X \quad \forall \, U \!\in {\mc N} \langle A \rangle \quad A \subset U$
\item \label{nhsystem sets 2} $\forall A \subset X \quad \forall \, U \in {\mc N} \langle A \rangle \quad \forall \, V \supset U \quad V \!\in {\mc N} \langle A \rangle$
\item \label{nhsystem sets 3} $\forall A \subset X \quad \forall \, U, V \!\in {\mc N} \langle A \rangle \quad U \cap V \in {\mc N} \langle A \rangle$
\item \label{nhsystem sets 4} $\forall A \subset X \quad \forall \, U \!\in {\mc N} \langle A \rangle \quad \exists \, V \!\in {\mc N} \langle A \rangle \quad U \!\in {\mc N} \langle V \rangle$
\item \label{nhsystem sets 5} For every index set~$I$ and every $A_i \subset X$ ($i \in I$), we have
\[
{\mc N} \big\langle \, {\textstyle \bigcup}_{i \in I} A_i \big\rangle \, = \, {\textstyle \bigcap}_{i \in I} \, {\mc N} \langle A_i \rangle
\]

\eenum

\elemm

\bproof
(\ref{nhsystem sets 1}) to (\ref{nhsystem sets 3}) and~(\ref{nhsystem sets 5}), follow by Lemma~\ref{lemm characterization neighborhood}. To show~(\ref{nhsystem sets 4}), let $A \subset X$ and $U \!\in {\mc N} \langle A \rangle$. We may choose $V \!\in \tp$ such that $A \subset V \!\subset U$ by Lemma~\ref{lemm characterization neighborhood}~(\ref{lemm characterization neighborhood 3}). Then $V \!\in {\mc N} \langle V \rangle$ by Lemma~\ref{lemm characterization neighborhood}~(\ref{lemm characterization neighborhood 2}), and therefore $U \!\in {\mc N} \langle V \rangle$.
\eproof

\blemm
\label{lemm neighborhood fixed filter}
Let $(X,\tp)$ be a topological space and $A \subset X$ with $A \neq \O$. The system ${\mc N} \langle A \rangle$ is a fixed filter on~$X$.
\elemm

\bproof
It follows by Lemma~\ref{nhsystem sets} (\ref{nhsystem sets 1}) to~(\ref{nhsystem sets 3}) that ${\mc N} \langle A \rangle$ is a filter. This filter is obviously fixed.
\eproof

Notice that in the following Lemma the statements (\ref{lemma top from neighborhood sets 1}) to~(\ref{lemma top from neighborhood sets 5}) correspond to the statements (\ref{nhsystem sets 1}) to~(\ref{nhsystem sets 5}) in Lemma~\ref{nhsystem sets}, the statement~(\ref{lemma top from neighborhood sets 6}) guarantees that $\mc M$ can be derived from a structure relation, and the statement~(\ref{lemma top from neighborhood sets 7}) excludes that the system is empty for any $A \subset X$.

\midvspace

\blemm
\label{lemma top from neighborhood sets}
Let $X$ be a set and ${\mc M} : {\mc P}(X) \longrightarrow {\mc P}^2(X)$ a map such that the following statements hold:

\benum
\item \label{lemma top from neighborhood sets 1}
$\forall A \subset X \quad \forall \, U \!\in {\mc M}(A) \quad A \subset U$

\item \label{lemma top from neighborhood sets 2}
$\forall A \subset X \quad \forall \, U \!\in {\mc M}(A) \quad \forall \, V \!\supset U \quad V \!\in {\mc M}(A)$

\item \label{lemma top from neighborhood sets 3}
$\forall A \subset X \quad \forall \, U, V \!\in {\mc M}(A) \quad U \cap V \in {\mc M}(A)$

\item \label{lemma top from neighborhood sets 4}
$\forall A \subset X \quad \forall \, U \!\in {\mc M}(A) \quad \exists \, V \!\in {\mc M}(A) \quad U \!\in {\mc M}(V)$

\item \label{lemma top from neighborhood sets 5}
For every index set~$I$ and every $A_i \subset X$ ($i \in I$), we have
\[
{\mc M} \! \left( {\textstyle \bigcup}_{i \in I} A_i \right) \, = \, {\textstyle \bigcap}_{i \in I} \, {\mc M} (A_i)
\]

\item \label{lemma top from neighborhood sets 6}
${\mc M} (\O) = {\mc P}(X)$

\item \label{lemma top from neighborhood sets 7}
$\forall A \subset X \quad {\mc M}(A) \neq \O$
\eenum

There is a unique topology $\tp$ on~$X$, such that ${\mc M}(A) = {\mc N} \langle A \rangle$ for every $A \subset X$, where ${\mc N}$ is the neighborhood system of~$(X,\tp)$.
\elemm

\bproof
We give two proofs. In the first one we use Lemma~\ref{lemma top from neighborhood}, while in the second one we do not.\\

{\em First proof:}\\
We may define a relation ${\mc N} \subset X \!\times {\mc P}(X)$ by ${\mc N} \! \left\{ x \right\} = {\mc M} (\left\{ x \right\})$ for every $x \in X$. Then ${\mc N}$ has properties (\ref{nhsystem0}) to~(\ref{nhsystem4}) in Lemma~\ref{nhsystem}. There is a unique topology $\tp$ on~$X$ such that $\mc N$ is the neighborhood system of $(X,\tp)$ by Lemma~\ref{lemma top from neighborhood}. Fix $A \subset X$ with $A \neq \O$. Then we have
\[
{\mc N} \langle A \rangle = \, {\textstyle \bigcap}_{x \in A} \, {\mc N} \! \left\{ x \right\} \, = \, {\textstyle \bigcap}_{x \in A} \, {\mc M} (\left\{ x \right\}) \, = \, {\mc M}(A)
\]

Moreover ${\mc N} \langle \O \rangle = {\mc P}(X)$.\\

{\em Second proof:}\\
We define $\tp = \left\{ U \subset X \, : \, \forall A \subset U \;\; U \!\in {\mc M}(A) \right\}$. First we show that $\tp$ is a topology on~$X$. Let $U_i \in \tp$ ($i \in I$), where $I$ is an index set, $U = \bigcup_{i \in I} U_i$, and $A \subset U$. Further let $A_i = A \cap U_i$ ($i \in I$). It follows that $A = \bigcup_{i \in I} A_i$. We have
\[
{\mc M}(A) = {\mc M} \left( {\textstyle \bigcup}_{i \in I} A_i \right) = \, {\textstyle \bigcap}_{i \in I} \, {\mc M}(A_i)
\]

by~(\ref{lemma top from neighborhood sets 5}). Moreover, for every $i \in I$, $A_i \subset U_i$ implies $U_i \in {\mc M}(A_i)$, and thus $U \in {\mc M}(A_i)$ by~(\ref{lemma top from neighborhood sets 2}). Hence we obtain $U \!\in {\mc M}(A)$. This proves that $\tp$ is a topology.

Let $\mc N$ be the neighborhood system of $\xi = (X,\tp)$. We show that ${\mc N} \langle A \rangle = {\mc M}(A)$ for every $A \subset X$. Clearly, ${\mc N} \langle \O \rangle = {\mc M}(\O)$. Fix $A \subset X$ with $A \neq \O$.

First assume that $U \!\in {\mc N} \langle A \rangle$. There exists $V \!\in \tp$ such that $A \subset V \subset U$. Therefore $V \!\in {\mc M}(A)$ by definition of~$\tp$, and thus $U \!\in {\mc M}(A)$.

Now assume $U \!\in {\mc M}(A)$, and define $V \! = \bigcup \left\{ B \subset X \, : \, U \in {\mc M}(B) \right\}$. Then we have
\begin{eqnarray*}
{\mc M}(V) \!\!\! & = & \!\! {\mc M}\left( \, \bigcup \big\{ B \subset X \, : \, U \!\in {\mc M}(B) \big\} \right)\\[.2em]
 & = & \!\! \bigcap \big\{ {\mc M}(B) \, : \, B \subset X, \; U \!\in {\mc M}(B) \big\}
\end{eqnarray*}

It follows that $U \!\in {\mc M}(V)$. Hence $U \!\in {\mc M}(B)$ for every $B \subset V$ by~(\ref{lemma top from neighborhood sets 5}). Fix $B \subset V$. Then there is $W \!\in {\mc M}(B)$ such that $U \!\in {\mc M}(W)$ by~(\ref{lemma top from neighborhood sets 4}). Hence $W \!\subset V$, and thus $V \!\in {\mc M}(B)$ by~(\ref{lemma top from neighborhood sets 2}). It follows that $V \!\in \tp$. Since $V \!\subset U$ by~(\ref{lemma top from neighborhood sets 1}), we have $U \!\in {\mc N} \langle A \rangle$.

The uniqueness follows by Lemma~\ref{lemm characterization neighborhood}~(\ref{lemm characterization neighborhood 2}).
\eproof

Lemma~\ref{lemm neighborhood fixed filter} says that the neighboorhood system of a non-empty set $A \subset X$ is a fixed filter. This suggests the idea to define the "neighborhood base" of~$A$ as a filter base for that filter. For a point $x \in X$, a neighborhood base is a subset of the neighborhood system of~$x$ such that each neighboorhood contains a member of the base. We begin with the definition of neighborhood base, continue with the definition of neighborhood base of a subset of~$X$, and then consider the neighborhood base of a point of~$X$ as a special case. It is then shown below that a neigborhood base of a subset is a filter base for the neighborhood system of that subset.

\midvspace

\bdefi
\index{Neighborhood base}
Let $\mc N$ be the neighborhood system of a topological space $\xi = (X,\tp)$. A system ${\mc B} \subset {\mc N}$ is called {\bf neighborhood base of}~$\xi$ if ${\mc N} = \Phi' ({\mc B})$ where $\Phi'$ is the map defined in Definition~\ref{defi phi prime}. We also say that $\mc B$ {\bf generates}~$\mc N$.
\edefi

In general, for a given topology there exists more than one neighborhood base.

\midvspace

\brema
Given a topological space $\xi$, the system ${\mc N}^{\, \mathrm{open}}_{\phantom{1}}$ is a neighborhood base of~$\xi$ by Lemma~\ref{rema properties nh system}~(\ref{rema properties nh system 2}).
\erema

\bdefi
Let $\xi = (X,\tp)$ be a topological space and $A \subset X$. A system ${\mc B} \subset {\mc N} \langle A \rangle$ is called {\bf neighborhood base of}~$A$ if $\Phi ({\mc B}) = {\mc N} \langle A \rangle$. If, in addition, $A = \left\{ x \right\}$ for some $x \in X$, then $\mc B$ is called {\bf neighborhood base of}~$x$.
\edefi

Again, for a given topology and a given set $A \subset X$, in general, more than one neighborhood base of $A$ exists.

\midvspace

\blemm
Let $\xi = (X,\tp)$ be a topological space, $\mc B$ a neighboorhood base of~$\xi$, and $x \in X$. Then ${\mc B} \left\{ x \right\}$ is a neighborhood base of~$x$.
\elemm

\bproof
This follows by Lemma~\ref{lemm phi prime relations}~(\ref{lemm phi prime relations 1}).
\eproof

Notice that, given a neighborhood base $\mc B$ of a topological space $(X,\tp)$ and $A \subset X$, ${\mc B} \langle A \rangle$ need not be a neighborhood base of $A$, because generally equality in Lemma~\ref{lemm phi prime relations}~(\ref{lemm phi prime relations 3}) may not hold.

\midvspace

\brema
\label{rema neighborhood filter}
Let $(X,\tp)$ be a topological space and $A \subset X$. The following statements hold:

\benum
\item \label{rema neighborhood filter 1} Every neighborhood base of~$A$ is a filter base for ${\mc N} \langle A \rangle$.
\item \label{rema neighborhood filter 2} The system ${\mc N}^{\, \mathrm{open}}_{\phantom{1}} \langle A \rangle$ is a neighborhood base of~$A$ by Lemma~\ref{lemm characterization neighborhood}~(\ref{lemm characterization neighborhood 3}). 
\eenum

Notice that ${\mc N}^{\, \mathrm{closed}}_{\phantom{1}} \langle A \rangle$ need not be a neighborhood base of~$A$.
\erema

There is a straightforward way to obtain a neighborhood base of a topological space from a topological base demonstrated in the following Lemma.

\midvspace

\blemm
\label{lemm top neigborhood base}
Given a topological space $\xi = (X,\tp)$ and a base $\mc A$ for~$\tp$. Then the structure relation $\mc B$ defined by
\[
(x,U) \in {\mc B} \;\; \Longleftrightarrow \;\; U \!\in {\mc A}(x)
\]

is a neighborhood base of~$\xi$. For every $x \in X$, ${\mc A}(x)$ is a neighborhood base of~$x$ that has solely open member sets.
\elemm

\bproof
Let $x \in X$ and $V \!\in {\mc N} \! \left\{ x \right\}$. We may choose $W \!\in \tp$ such that $x \in W \!\subset V$, and $U \!\in {\mc A}$ such that $x \in U \!\subset W$. Then $(x,U) \in {\mc B}$. The second claim is obvious.
\eproof

We now address the following two questions: For a given topological space, can a topological base be constructed from a neighborhood base? And: Without specifying a topology beforehand, when is a given structure relation a neighborhood base of some topology?

\midvspace

\blemm
Let $(X,\tp)$ be a topological space and $\mc B$ a neighborhood base that has only open member sets. The system ${\mc B} \left[ X \right] = \bigcup_{x \in X} {\mc B} \left\{ x \right\}$ is a topological base. 
\elemm

\bproof
Exercise.
\eproof

Concerning the second question we proceed similarly as in the case of neighborhood systems and derive a list of properties of a neighborhood base and then show below that if a structure relation on a set~$X$ has these properties, then there is a topology~$\tp$ on~$X$ such that the structure relation is a neighborhood base of~$(X,\tp)$.

\midvspace

\blemm
\label{nhbase}
Let $\mc B$ be a neighborhood base of a topological space. $\mc B$~has the following properties:
\benum
\item \label{nhbase0} $\forall x \in X \quad {\mc B} \left\{ x \right\} \neq 0$
\item \label{nhbase1} $\forall x \in X \quad \forall \, U \!\in {\mc B} \left\{ x \right\} \quad x \in U$
\item \label{nhbase3} $\forall x \in X \quad \forall \, U, V \!\in {\mc B} \left\{ x \right\} \quad \exists \, W \!\in {\mc B} \left\{ x \right\} \quad W \!\subset U \cap V$
\item \label{nhbase4} $\forall x \in X \quad \forall \, U \!\in {\mc B} \left\{ x \right\} \quad \exists \, V \!\in {\mc B} \left\{ x \right\} \quad \forall y \in V \quad \exists \, W \!\in {\mc B} \left\{ y \right\} \quad W \!\subset U$
\eenum

\elemm

\bproof
This follows by Lemma~\ref{nhsystem}.
\eproof

\midvspace

\blemm
\label{lemma top from neighborhood base}
Let $X$ be a set and ${\mc B} \subset X \times {\mc P}(X)$ such that properties (\ref{nhbase0}) to (\ref{nhbase4}) in Lemma~\ref{nhbase} are satisfied. There is a unique topology $\tp$ on~$X$ such that $\mc B$ is a neighborhood base of~$(X,\tp)$.
\elemm

\bproof
The structure relation ${\mc N} = \Phi' ({\mc B})$ satisfies conditions (\ref{nhsystem0}) to (\ref{nhsystem4}) in Lemma \ref{nhsystem}. Hence there exists a topology~$\tp$ on $X$ by Lemma~\ref{lemma top from neighborhood} such that $\mc N$ is neighborhood system of $(X,\tp)$, and thus $\mc B$ is a neighborhood base of~$(X,\tp)$.
\eproof

Recall that a topological space that has a countable base is called second countable. The following definition refers to the cardinality of the neighborhood base of every point of~$X$.

\midvspace

\bdefi
\index{First countable}
\index{Topology!first countable}
\index{Topological space!first countable}
Let $\xi = (X,\tp)$ be a topological space. If there is a neighborhood base $\mc B$ of~$\xi$ such that ${\mc B} \left\{ x \right\}$ is countable for every $x \in X$, then the space~$\xi$ or the topology $\tp$ is called {\bf first countable}.
\edefi

\bexam
Given a set $X$, the discrete topology $\tp_{\mathrm{dis}}$ on~$X$, and a point $x \in X$, the neighborhood system of~$x$ is $\tp_{\mathrm{dis}}(x)$. $\tp_{\mathrm{dis}}$ is clearly first countable.
\eexam

\bexam
Given a set $X$, the indiscrete topology $\tp_{\mathrm{in}}$ on~$X$, and a point $x \in X$, the neighborhood system of $x$ is~$\left\{ X \right\}$.
\eexam

\blemm
\index{Second countable}
\index{Topology!second countable}
\index{Topological space!second countable}
A second countable topological space $(X,\tp)$ is also first countable.
\elemm

\bproof
This follows by Lemma~\ref{lemm top neigborhood base}.
\eproof

\blemm
\label{lemm countable int base}
Let $(X,\tp)$ be a first countable topological space and $x \in X$. There is a neighborhood base $\left\{ B_m \, : \, m \in \naturalnumbers \right\}$ of~$x$ such that all $B_m$ ($m \in \naturalnumbers$) are open and
\[
\forall m, n \in \naturalnumbers \quad m < n \;\; \Longrightarrow \;\; B_n \subset B_m
\]

\elemm

\bproof
By Lemma~\ref{lemm top neigborhood base}, we may choose a neighborhood base $\left\{ C_m \, : \, m \in \naturalnumbers \right\}$ for~$x$ such that all $C_m$ ($m \in \naturalnumbers$) are open. For every $n \in \naturalnumbers$ we define $B_n = \bigcap \left\{ C_m \, : \, m \in \naturalnumbers,\, m \leq n \right\}$.
\eproof

We conclude this section with a criterion for the comparison of two topologies on the same set in terms of neighborhoods and neighborhood bases.

\midvspace

\blemm
\label{base finer neighborhood}
Let $\tp_1$ and $\tp_2$ be two topologies on a set~$X$, ${\mc N}_1$ and ${\mc N}_2$ their respective neighborhood systems, and ${\mc B}_1$ and ${\mc B}_2$ two neighborhood bases generating ${\mc N}_1$ and ${\mc N}_2$, respectively. The following three statements are equivalent:

\benum
\item \label{base finer neighborhood 1} $\tp_2 \subset \tp_1$
\item \label{base finer neighborhood 2} $\forall x \in X \quad {\mc N}_2 \!\left\{ x \right\} \subset \, {\mc N}_1 \!\left\{ x \right\}$
\item \label{base finer neighborhood 3} $\forall x \in X \quad {\mc B}_2 \!\left\{ x \right\} \subset_{\stackrel{}{\Phi}} {\mc B}_1 \!\left\{ x \right\}$
\eenum

\elemm

\bproof
Since ${\mc B}_1 \! \left\{ x \right\}$ is a filter base for~${\mc N}_1 \!\left\{ x \right\}$, and ${\mc B}_2 \!\left\{ x \right\}$ is a filter base for~${\mc N}_2 \!\left\{ x \right\}$ by Remark~\ref{rema neighborhood filter}~(\ref{rema neighborhood filter 1}), statements (\ref{base finer neighborhood 2}) and~(\ref{base finer neighborhood 3}) are equivalent by Lemma~\ref{coro filter bases comparison}~(\ref{coro filter bases comparison 1}). The fact that (\ref{base finer neighborhood 1}) implies (\ref{base finer neighborhood 2}) follows by the definition of neighborhood system. The converse implication follows by Lemma~\ref{lemm characterization neighborhood}~(\ref{lemm characterization neighborhood 1}).
\eproof

\section{Interval topology}
\label{interval topology}

In this and in the next Section we discuss two specific kinds of topologies: those associated to pre-orderings with full field, and those generated by pseudo-metrics.

\midvspace

\blede
\label{lede interval top}
\index{Interval topology}
\index{Topology!interval}
Given a pre-ordered space $(X,R)$ where $R$ has full field, the system 
\[
{\mc S} \, = \, \big\{ \left] -\infty, x \right[\, , \;\left] x, \infty \right[\; : \, x \in X \big\} \cup \left\{ \O \right\}
\]

is a topological subbase on~$X$. The generated topology is called $R$-{\bf interval topology}, or short {\bf interval topology}, and also written~$\tau(R)$. If $Y \subset X$ is order dense, then also the system
\[
{\mc R} \, = \, \big\{ \left] -\infty, y \right[\, , \;\left] y, \infty \right[\; : \, y \in Y \big\} \cup \left\{ \O \right\}
\]

is a subbase for the interval topology.
\elede

\bproof
We first show that $\mc S$ is a topological subbase by Lemma~\ref{lemm crit subbase}. Since $R$ has full field, there is, for each $x \in X$, a point $y \in X$ such that $x \in \left] -\infty, y \right[\,$ or $x \in \left] y, \infty \right[\,$. It follows that $\bigcup {\mc S} = X$, and clearly ${\mc S} \neq \O$.

If $Y$ is order dense, then ${\mc S} \subset \Theta ({\mc R})$ by Remark~\ref{rema ordered dense}, and therefore $\bigcup {\mc R} = X$. Thus $\mc R$ is a topological subbase. Moreover, $\Psi ({\mc S}) \subset \Psi \, \Theta \, ({\mc R}) \subset \Theta \, \Psi \, ({\mc R})$ by Lemma~\ref{lemm set relations}, and hence $\Theta \, \Psi \, ({\mc S}) \subset \Theta \, \Psi \, ({\mc R})$. Therefore $\mc R$ and $\mc S$ generate the same topology.
\eproof

Notice that for some pre-ordered spaces $(X,\prec)$ the system
\[
{\mc A} = \big\{ \left] -\infty, x \right[\, , \;\left] x, \infty \right[\; : \, x \in X \big\}
\]

may contain $\O$ or may not have the finite intersection property. In these two cases $\mc A$ alone (without explicitly including~$\O$) is a topological subbase.

Under additional assumptions on the pre-ordered space, canonical bases for the interval topology can be specified.

\midvspace

\blemm
\label{lemm interval top base}
Let $(X,\prec)$ be a pre-ordered space where $\prec$ has full field and the systems
\[
{\mc S}_- = \big\{ \left] -\infty, x \right[\; : \, x \in X \big\} \cup \left\{ \O \right\},\quad \quad
{\mc S}_+ = \big\{ \left] x, \infty \right[\; : \, x \in X \big\} \cup \left\{ \O \right\}
\]

satisfy $\Psi ({\mc S}_-) = {\mc S}_-$ and $\Psi ({\mc S}_+) = {\mc S}_+$, that is ${\mc S}_-$ and ${\mc S}_+$ are fixed points of~$\Psi$. Further let
\[
{\mc S} \, = \, {\mc S}_+ \cup \, {\mc S}_-\,, \quad \quad {\mc A} = \big\{ \left] x, y \right[\; : \, x, y \in X,\; x \prec y \big\} \cup \left\{ \O \right\}
\]

The system ${\mc A} \cup {\mc S}$ is a base for the interval topology. If $\prec$ has full domain and full range, then $\mc A$ alone is a base for the interval topology.

Furthermore, let $Y \subset X$ be an order dense set such that the systems
\[
{\mc R}_- = \big\{ \left] -\infty, x \right[\; : \, x \in Y \big\} \cup \left\{ \O \right\}, \quad \quad {\mc R}_+ = \big\{ \left] x, \infty \right[\; : \, x \in Y \big\} \cup \left\{ \O \right\}
\]

satisfy $\Psi ({\mc R}_-) = {\mc R}_-$ and $\Psi ({\mc R}_+) = {\mc R}_+$. Moreover, let
\[
{\mc R} \, = \, {\mc R}_- \cup \, {\mc R}_+ \,, \quad \quad {\mc B} = \big\{ \left] x, y \right[\; : \, x, y \in Y,\; x \prec y \big\} \cup \left\{ \O \right\}
\]

Then ${\mc B} \cup {\mc R}$ is a base for the interval topology. If $\prec$ has full domain and full range, then $\mc B$ alone is a base for the interval topology.
\elemm

\bproof
To prove that ${\mc A} \cup {\mc S}$ is a base for the interval topology, we show that $\Psi ({\mc S}) = {\mc A} \cup {\mc S}$. Let $\O \neq A \in \Psi ({\mc S})$. We have
\[
\left( A = A_E \right) \; \vee \; \left( A = A_F \right) \; \vee \; \left( A = A_E \cap A_F \right)
\]

where
\[
A_E \, = \; \bigcap \big\{ \left] x, \infty \right[ \; : \, x \in E \big\}, \quad A_F \, = \; \bigcap \big\{ \left] -\infty, y \right[ \; : \, y \in F \big\}
\]

and $E, F \sqsubset X$ are index sets. Since ${\mc S}_-$ and ${\mc S}_+$ are fixed points of~$\Psi$, there are points $x, y \in X$ such that $A_E = \; \left] x, \infty \right[ \;$ and $A_F = \; \left] -\infty, y \right[\;$.

To prove the second claim we show that $\Theta \, \Psi ({\mc S}) = \Theta ({\mc A})$ under the given conditions. We may choose $A$, $A_E$, and $A_F$ as above. Since $\prec$ has full domain, we have
\[
A_E \, = \; \left] x, \infty \right[ \; = \; \bigcup \big\{ \left] x, z \right[ \; : \, z \in X,\; x \prec z \big\}
\]

by Remark~\ref{rema order union}. Thus, if $A = A_E$, then $A \in \Theta ({\mc A})$. Similarly, since $\prec$ has full range, we have
\[
A_F \, = \; \left] -\infty, y \right[\; = \; \bigcup \big\{ \left] z, y \right[ \; : \, z \in X,\;  z \prec y \big\}
\]

Hence, if $A = A_F$, then $A \in \Theta ({\mc A})$. If $A = A_E \cap A_F$, we have $A = \, \left] x, y \right[ \,$ and therefore $A \in \Theta ({\mc A})$.

To see that ${\mc B} \cup {\mc R}$ is a base for the interval topology, we show that $\Psi ({\mc R}) = {\mc B} \cup {\mc R}$. Let $\O \neq B \in \Psi ({\mc R})$. Since ${\mc R}_-$ and ${\mc R}_+$ are fixed points of~$\Psi$, we have either $B = \; \left] u, \infty \right[ \;$ or $B = \; \left] -\infty, v \right[\;$ or $B = \; \left] u, v \right[ \;$ where $u, v \in Y$.

The last claim follows by Remarks~\ref{rema order union} and~\ref{rema ordered dense}.
\eproof

Note that in the particular case in which $\prec$ is a connective pre-ordering each of the systems ${\mc S}_-$, ${\mc S}_+$, ${\mc R}_-$, and ${\mc R}_+$ is a fixed point of~$\Psi$ by Lemma~\ref{lemm tot ord subset min max}.

In order to handle conditions such as those of Lemma~\ref{lemm interval top base} in a stringent manner, we define the following notions.

\midvspace

\bdefi
Let $(X,\prec)$ be a pre-ordered space. If the systems
\[
{\mc S}_- = \big\{ \left] -\infty, x \right[\; : \, x \in X \big\} \cup \left\{ \O \right\}, \quad \quad {\mc S}_+ = \big\{ \left] x, \infty \right[\; : \, x \in X \big\} \cup \left\{ \O \right\}
\]

are fixed points of~$\Psi$, we say that $\prec$ {\bf has the interval intersection property}. If, in addition, $\prec$ has full domain and full range, then $\prec$ is called {\bf interval relation}.
\edefi

\bdefi
\index{Standard topology}
\index{Topology!standard}
The interval topology of $({\mathbb R},<)$, where $<$ is the standard ordering in the sense of~"$<$" on~${\mathbb R}$, is called {\bf standard topology on}~${\mathbb R}$. The interval topology on $({\mathbb R}_+,<)$ is called {\bf standard topology on}~${\mathbb R}_+$. 
\edefi

\brema
\label{rema r interval rel}
The system  $\big\{ \left] -\infty, x \right[\, , \, \left] x, \infty \right[\, : \, x \in {\mathbb R} \big\}$ is a subbase for the standard topology on~$\mathbb R$. The ordering~$<$ is an interval relation on~${\mathbb R}$. Thus the system
\[
\big\{ \left] x, y \right[\; : \, x, y \in {\mathbb R},\; x < y \big\} \cup \left\{ \O \right\}
\]

is a base for the standard topology by Lemma~\ref{lemm interval top base}. Moreover, the system
\[
\big\{ \left] x, y \right[\; : \, x, y \in {\mathbb D},\; x < y \big\} \cup \left\{ \O \right\}
\]

is a base for the same topology by Lemmas~\ref{lemm D dense} and~\ref{lemm interval top base}. Hence the standard topology on~${\mathbb R}$ is second countable by Remark~\ref{rema dyadic countable} and Lemma~\ref{lemm countable sets}.
\erema

\brema
\label{rema r+ interval rel}
The system
\[
{\mc S} = \big\{ \left] -\infty, x \right[, \,\left] x, \infty \right[\; : \, x \in {\mathbb R}_+ \big\}
\]

is a subbase for the standard topology on~${\mathbb R}_+$. Let
\[
{\mc A} = \big\{ \left] x, y \right[\; : \, x, y \in {\mathbb R}_+,\; x < y \big\} \cup \left\{ \O \right\}
\]

The system ${\mc S} \cup {\mc A}$ is a base for the standard topology by Lemma~\ref{lemm interval top base}. The relation~$<$ has the interval intersection property. However it is not an interval relation because it has a minimum. Further let
\[
{\mc R} = \big\{ \left] -\infty, x \right[, \,\left] x, \infty \right[\; : \, x \in {\mathbb D}_+ \big\},\quad {\mc B} = \big\{ \left] x, y \right[\; : \, x, y \in {\mathbb D}_+,\; x < y \big\} \cup \left\{ \O \right\}
\]

The system ${\mc R} \cup {\mc B}$ is a base for the same topology by Lemmas~\ref{lemm D0 dense} and~\ref{lemm interval top base}. Hence the standard topology on~${\mathbb R}_+$ is second countable by Corollary~\ref{coro dyadic countable} and Lemmas~\ref{lemm countable sets} and~\ref{lemm count times count} .
\erema

We now generalize the concept of interval topology to the case of a system of pre-orderings on the same set~$X$.

\midvspace

\blede
\label{defi r int top}
Let $X$ be a set and $\mc R$ a system of pre-orderings on~$X$ where $\mc R$ has full field. Then the system
\[
{\mc S} = \big\{ \left] -\infty, x \right[_{\, R}, \;\left] x, \infty \right[_{\, R} \, : \, x \in X,\, R \in {\mc R} \big\} \cup \left\{ \O \right\}
\]

is a topological subbase. The topology generated by~$\mc S$ is called {\bf ${\mc R}$-interval topology}, and written~$\tau({\mc R})$. Further let $Y \subset X$ be $\mc R$-dense in~$X$. Then
\[
{\mc Q} = \big\{ \left] -\infty, y \right[_{\, R}, \;\left] y, \infty \right[_{\, R} \, : \, y \in Y,\, R \in {\mc R} \big\} \cup \left\{ \O \right\}
\]

is a subbase for the same topology.
\elede

\bproof
To see that $\mc S$ is a topological subbase, note that for every $x \in X$ there is $R \in {\mc R}$ and $y \in X$ such that $x \in \, \left] -\infty, y \right[_{\, R}$ or $x \in \, \left] y, \infty \right[_{\, R} \, $. It follows that $\bigcup {\mc S} = X$, and clearly ${\mc S} \neq \O$.

To see the second claim, assume that $Y$ is $\mc R$-dense. Then ${\mc S} \subset \Theta ({\mc Q})$ by Remark~\ref{rema ordered dense}, and therefore $\bigcup {\mc Q} = X$. Thus $\mc Q$ is a topological subbase. Moreover, $\Psi ({\mc S}) \subset \Psi \, \Theta \, ({\mc Q}) \subset \Theta \, \Psi \, ({\mc Q})$ by Lemma~\ref{lemm set relations}, and hence $\Theta \, \Psi \, ({\mc S}) \subset \Theta \, \Psi \, ({\mc Q})$. Therefore $\mc Q$ and $\mc S$ generate the same topology.
\eproof

\blemm
\label{lemm pre-order equal}
Let $X$ be a set and $\mc R$ a system of pre-orderings on~$X$ where $\mc R$ is independent and has full field. Further let $S = \bigcap {\mc R}$. Then we have $\tau({\mc R}) \subset \tau(S)$. If $\mc R$ is finite, then $\tau(S) \subset \tau({\mc R})$
\elemm

\bproof
This follows by Lemma~\ref{lemm independent induced rel}.
\eproof

Example~\ref{exam rn ordering} is an important example of the construction defined in Definition~\ref{defi r int top}. We introduce the following notion.

\midvspace

\bdefi
\label{defi int top rn}
\index{Standard topology}
\index{Topology!standard}
Let $n \in \naturalnumbers$, $n \geq 1$, and $(X_i,R_i) = ({\mathbb R},<)$ for every $i \in \naturalnumbers$, $1 \leq i \leq n$. Further let $X = \bigtimes_{\!\! i = 1}^{\!\! n} \, X_i$ and ${\mc R} = \left\{ p_i^{-1} \left[ R_i \right] \, : \, i \in \naturalnumbers, \, 1 \leq i \leq n \right\}$. The $\mc R$-interval topology on ${\mathbb R}^n$ is called {\bf standard topology on}~${\mathbb R}^n$.
\edefi
 
\brema
\label{rema rn interval rel}
With definitions as in Definition~\ref{defi int top rn}, the system
\[
{\mc S} = \big\{ \left] -\infty, x \right[_{\,R}, \;\left] x, \infty \right[_{\,R} \; : \, x \in {\mathbb D}^n,\, R \in {\mc R} \big\}
\]

is a subbase for the standard topology on~${\mathbb R}^n$ by Lemma~\ref{defi r int top}. Since $\Psi ({\mc S})$ is a base, this topology is second countable by Remark~\ref{rema dyadic countable} and Lemmas~\ref{lemm countable sets}, \ref{lemm count times count}, and~\ref{lemm fin count}. Let $S = \bigcap {\mc R}$. Since $\mc R$ is independent, the $S$-interval topology is the standard topology on~${\mathbb R}^n$ by Lemma~\ref{lemm pre-order equal}. Since $S$ is an interval relation, the system
\[
\big\{ \left] x, y \right[ \; : \, x, y \in {\mathbb R}^n,\, (x,y) \in S \big\} \cup \left\{ \O \right\},
\]

where the interval refers to the ordering~$S$, is a base for this topology by Lemma~\ref{lemm interval top base}. Similarly, it is clear that the system
\[
\big\{ \left] x, y \right[ \; : \, x, y \in {\mathbb D}^n,\, (x,y) \in S \big\} \cup \left\{ \O \right\}
\]

is a base for the same topology.
\erema

Similarly, a standard topology on any finite product of positive reals can be defined as follows.

\midvspace

\bdefi
\label{defi int top r+n}
\index{Standard topology}
\index{Topology!standard}
Let $n \in \naturalnumbers$, $n \geq 1$, and $(X_i,R_i) = ({\mathbb R}_+,<)$ for every $i \in \naturalnumbers$, $1 \leq i \leq n$. Further let ${\mc R} = \left\{ p_i^{-1} \left[ R_i \right] \, : \, i \in \naturalnumbers,\, 1 \leq i \leq n \right\}$. The $\mc R$-interval topology on ${\mathbb R}_+^n$ is called {\bf standard topology on}~${\mathbb R}_+^n$. 
\edefi
 
\brema
\label{rema r+n interval rel}
With definitions as in Definition~\ref{defi int top r+n}, the system
\[
{\mc S} = \big\{ \left] -\infty, x \right[_{\,R}, \;\left] x, \infty \right[_{\,R} \, : \, x \in {\mathbb D}_+^n,\, R \in {\mc R} \big\}
\]

is a subbase for the standard topology on~${\mathbb R}_+^n$ by Lemma~\ref{defi r int top}. Since $\Psi ({\mc S})$ is a base, this topology is second countable by Corollary~\ref{coro dyadic countable} and Lemmas~\ref{lemm countable sets}, \ref{lemm count times count}, and~\ref{lemm fin count}. $\mc R$~is upwards independent, but not downwards independent. 
\erema

It is possible to construct other topologies on a pre-ordered space. While using all improper intervals and the empty set as a subbase in Lemma and Definition~\ref{lede interval top}, we may use only the lower or only the upper segments, respectively supplemented by the empty set, as subbase as done in the following Lemma.

\midvspace

\blemm
\label{lemm left right interval top}
Let $(X,\prec)$ be a pre-ordered space. If $\prec$ has full domain, then the system
\[
{\mc S}_- = \big\{ \left] -\infty, x \right[\; : \, x \in X \big\} \cup \left\{ \O \right\}
\]

is a topological subbase on~$X$. If $\prec$ has full range, then the system
\[
{\mc S}_+ = \big\{ \left] x, \infty \right[\; : \, x \in X \big\} \cup \left\{ \O \right\}
\]

is a topological subbase on~$X$.

Let $Y \subset X$ be order dense. If $\prec$ has full domain, then the system
\[
{\mc R}_- = \big\{ \left] -\infty, y \right[\; : \, y \in Y \big\} \cup \left\{ \O \right\}
\]

is a topological subbase on~$X$. If $\prec$ has full range, then the system
\[
{\mc R}_+ = \big\{ \left] y, \infty \right[\; : \, y \in Y \big\} \cup \left\{ \O \right\}
\]

is a topological subbase on~$X$.
\elemm

\bproof
Exercise.
\eproof

Clearly, if a topological subbase is a fixed point of~$\Psi$, it is a base for its generated topology. In case of the interval topology it is, under certain conditions, even a topology if we only include $X$ to the set system.

\midvspace

\blemm
Let $(X,\prec)$ be a totally ordered space in the sense of "$<$". We define the systems
\begin{eqnarray*}
{\mc S}_- \!\!\!\! & = & \!\!\! \big\{ \left] -\infty, x \right[\; : \, x \in X \big\} \cup \left\{ \O, X \right\},\\[.2em]
{\mc S}_+ \!\!\!\! & = & \!\!\! \big\{ \left] x, \infty \right[\; : \, x \in X \big\} \cup \left\{ \O, X \right\}
\end{eqnarray*}

If $\prec$ has the least upper bound property, then $\Theta ({\mc S}_-) = {\mc S}_-$ and $\Theta ({\mc S}_+) = {\mc S}_+$. In this case, ${\mc S}_-$ and ${\mc S}_+$ are topologies on~$X$.
\elemm

\bproof
Let $I$ be an index set, $x_i \in X$ ($i \in I$), and $A = \left\{ x_i \, : \, i \in I \right\}$. If $A$ has no upper bound, then we have $X = \bigcup_{i \in I} \; \left] -\infty, x_i \right[ \;$. If $A$ has an upper bound, let $x$ be the supremum of~$A$. Then we have $\, \left] -\infty, x \right[ \, = \bigcup_{i \in I} \; \left] -\infty, x_i \right[ \;$. It follows that $\Theta ({\mc S}_-) = {\mc S}_-$. The second equation follows similarly by Theorem~\ref{theo least upper greatest lower}. Since ${\mc S}_-$ and ${\mc S}_+$ are fixed points of~$\Psi$, the last claim follows.
\eproof

\bexam
\label{rless}
Let $<$ denote the standard ordering in the sense of "$<$" on~$\mathbb R$. The system
\[
\tp_< = \big\{ \left]-\infty, x \right[ \; : \, x \in {\mathbb R}\big\} \cup \left\{ \O, X \right\}
\]

is a topology on~$\mathbb R$. The system
\[
\big\{ \left] -\infty, x \right[ \; : \, x \in {\mathbb D} \big\} \cup \left\{ \O \right\}
\]

is a base for this topology.
\eexam

\bexam
\label{rleq}
Let $\leq$ denote the standard ordering in the sense of "$\leq$" on~$\mathbb R$. Notice that we have $x \in \, \left]-\infty, x \right[ \,$ for every $x \in {\mathbb R}$. The system
\[
\big\{ \left]-\infty, x \right[ \, : \, x \in {\mathbb R} \big\} \cup \left\{ \O \right\}
\]

is a topological base. However, the system
\[
\big\{ \left]-\infty, x \right[ \, : \, x \in {\mathbb D} \big\} \cup \left\{ \O \right\}
\]

is not a base for the generated topology. The conditions of Lemma~\ref{lemm left right interval top} are not satisfied as ${\mathbb D}$ is not order dense. 
\eexam

\section{Pseudo-metrics}

In this Section we consider topologies generated by pseudo-metrics and metrics. Some of the most important examples of metric spaces are the real numbers $\mathbb R$ and their $n$-fold Cartesian product $\mathbb R^n$ together with the Euclidean metric. This metric is only introduced in Chapter~\ref{functions on Rn} as it requires the definition of the square root function on the reals. The concept of pseudo-metric space can be generalized to that of uniform space, which is a set together with a system of---generally more than one---pseudo-metrics.

We start with the basic definitions.

\midvspace

\bdefi
\label{metric}
\index{Pseudo-metric}
\index{Pseudo-metric space}
\index{Metric}
\index{Metric space}
Given a set $X$, a {\bf pseudo-metric on} $X$ is a map $d: X \times X \longrightarrow {\mathbb R}_+$ with the following properties:
\benum
\item \label{metric1} $\forall x \in X \quad d(x, x) = 0$
\item \label{metric2} $\forall x, y \in X \quad d(y,x) = d(x,y)$ \quad (symmetry)
\item \label{metric3} $\forall x, y, z \in X \quad d(x,z) \leq d(x,y) + d(y,z)$ \quad (triangle inequality)
\eenum

The pair $(X,d)$ is called {\bf pseudo-metric space}. The map $d$ is called {\bf metric on}~$X$ if the following statement holds instead of~(\ref{metric1}):
\benum
\item[(i)$^\prime$] $\forall x, y \in X \quad x = y \quad \Longleftrightarrow \quad d(x, y) = 0$
\eenum
In this case $(X,d)$ is called {\bf metric space}.
\edefi

\bdefi
\index{Bounded}
\index{Pseudo-metric!bounded}
Given a pseudo-metric space $(X,d)$, $d$ is called {\bf bounded} if there is $r \in {\mathbb R}_+$ such that
\[
\forall x, y \in X \quad d(x,y) < r
\]

\edefi

\brema
Given a set $X$ and a bounded pseudo-metric~$d$ on~$X$, we have
\[
\sup \big\{ d(x,y) \, : \, x, y \in X \big\} < \infty
\]

by Lemma~\ref{lemm least upper bound reals}.
\erema 

\bexam
Let $X$ be a set and $(Y,d)$ a bounded pseudo-metric space. The map
\[
D : Y^X \!\!\times Y^X \longrightarrow {\mathbb R}_+\,, \quad \quad D(f,g) = \, {\textstyle \sup_{x \in X}} \, d \,\big(f(x), g(x)\big)
\]

is a bounded pseudo-metric on~$Y^X$. If $d$ is a metric, then $D$ is a metric.
\eexam

\bdefi
\index{Isometry}
Let $(X,d_X)$ and $(Y,d_Y)$ be pseudo-metric spaces and $f : X \longrightarrow Y$ a surjective map. $f$~is called {\bf isometry} if \,\!\! $d_Y (f(x),f(y)) = d_X(x,y)$ for every $x, y \in X$.
\edefi

In the following it is shown how a topology can be related to a given pseudo-metric on a set~$X$. In this sense, a pseudo-metric space is a special case of a topological space.

\midvspace

\blede
\label{lede pseudo-m topology}
\index{Topology!pseudo-metric}
\index{Topology!metric}
\index{Open sphere}
\index{Closed sphere}
Let $(X,d)$ be a pseudo-metric space. We define the function
\[
B : X \!\times \left] 0, \infty \right[\; \longrightarrow  {\mc P}(X), \quad B(x,r) = \left\{ y \in X \, : \, d(x,y) < r \right\}
\]

Further let $R \subset \, \left] 0, \infty \right[ \,$ such that for every $K \in \, \left] 0, \infty \right[ \,$ there exists $r \in R$ with $r < K$. The system
\[
{\mc B} \, = \, \big\{ B(x,r) \, : \, x \in X,\; r \in R \big\} \cup \left\{ \O \right\}
\]

is a topological base on~$X$. The topology~$\Theta ({\mc B})$ is called {\bf pseudo-metric topology of}~$(X,d)$. This topology is also denoted by~$\tau(d)$. The definition of~$\tau(d)$ is independent of the choice of~$R$. We say that $d$ {\bf generates}~$\tau(d)$.

If $d$ is a metric, then $\tau(d)$ is also called {\bf metric topology}. For every $r \in \, \left] 0, \infty \right[ \,$ and $x \in X$, the set $B(x,r)$ is $\tau(d)$-open. It is called {\bf open sphere about $x$ with $d$-radius $r$}. For every $r \in \, \left] 0, \infty \right[ \,$ and $x \in X$, the set $\left\{ y \in X \, : \, d(x,y) \leq r \right\}$ is $\tau(d)$-closed. It is called {\bf closed sphere about~$x$ with $d$-radius~$r$}. Moreover, the system 
\[
{\mc C} = \big\{ \big( x, B(x,r) \big) \, : \, x \in X,\, r \in R \big\}
\]

is a neighborhood base of~$(X, \tau(d))$.
\elede

\bproof
The claims that $\mc B$ is a topological base on~$X$ and that closed spheres are closed sets follow by Lemma and Definition~\ref{lemm char base} and the triangle inequality (exercise).

Now it is clear that $\mc C$ is a neighborhood base.

To see that $\tau(d)$ is independent of the choice of the radii, let, for $i \in \left\{ 1, 2 \right\}$, $R_i \subset \, \left] 0, \infty \right[ \,$ such that for every $K \in \, \left] 0, \infty \right[ \,$ there exists $r \in R_i$ with $r < K$. We denote by $\tp_i$ ($i \in \left\{ 1, 2 \right\}$) the respective pseudo-metric topologies. For $i \in \left\{ 1, 2 \right\}$, the system
\[
{\mc C}_i = \big\{ \big( x, B(x,r) \big) \, : \, x \in X,\, r \in R_i \big\}
\]

is a neighborhood base of~$(X, \tp_i)$. For every $x \in X$ we have ${\mc C}_2 \! \left\{ x \right\} \subset_{\stackrel{}{\Phi}} {\mc C}_1 \! \left\{ x \right\}$ and ${\mc C}_1 \! \left\{ x \right\} \subset_{\stackrel{}{\Phi}} {\mc C}_2 \! \left\{ x \right\}$. It follows by Lemma~\ref{base finer neighborhood} that $\tp_1 = \tp_2$.
\eproof

Note that we may choose $R = \, \left] 0, \infty \right[ \;$ in Lemma and Definition~\ref{lede pseudo-m topology}.

\midvspace

\bexam
Let $X$ be a set and $d : X \times X \longrightarrow {\mathbb R}_+$ a map with $d(x,x) = 0$ for every $x \in X$, and $d(x,y) = 1$ for every $x, y \in X$ with $x \neq y$. Then $d$ is a bounded metric, and $\tau(d)$ is the discrete topology on~$X$.
\eexam

\bexam
Let $X$ be a set and $d : X \times X \longrightarrow {\mathbb R}_+$ a map with $d(x,y) = 0$ for every $x,y \in X$. Then $d$ is a bounded pseudo-metric, and $\tau(d)$ is the indiscrete topology on~$X$.
\eexam

\blemm
\label{lemm b metric}
The function $d : {\mathbb R} \times {\mathbb R} \longrightarrow {\mathbb R}_+$, $d(x,y) = |x - y|$, is a metric. The generated topology $\tau(d)$ is the standard topology on~${\mathbb R}$.
\elemm

\bproof
Exercise.
\eproof

\blede
\label{lemm sup metric}
Let $n \in \naturalnumbers$ with $n \geq 1$. The function
\[
d : {\mathbb R}^n \!\times {\mathbb R}^n \longrightarrow {\mathbb R}_+\,, \quad d(x,y) = \mathrm{max} \, \big\{ |x_k - y_k| \, : \, k \in \naturalnumbers,\, 1 \leq k \leq n \big\}
\]

is a metric. It is called {\bf maximum metric on}~${\mathbb R}^n$.
\elede

\bproof
Exercise.
\eproof

It is shown in Chapter~\ref{functions on Rn} that the topology generated by the maximum metric is the standard topology on~${\mathbb R}^n$.

Given a topological space one may raise the question whether there exists a metric or pseudo-metric that generates it.

\midvspace

\bdefi
\index{Topological space!pseudo-metrizable}
\index{Topological space!metrizable}
Given a topological space $\xi = (X,\tp)$, $\xi$ is called {\bf pseudo-metrizable} if there exists a pseudo-metric $d$ on $X$ such that $\tau(d) = \tp$. Similarly, $\xi$ is called {\bf metrizable} if there exists a metric $d$ on $X$ such that $\tau(d) = \tp$.
\edefi

Given the fact that pseudo-metric and metric spaces are very important mathematical concepts, pseudo-metrizability and metrizability are central features a topological space may or may not have.

The following Lemma provides a tool to compare the topologies generated by two different pseudo-metrics on the same set $X$. In particular, it can be used to prove that two specific different pseudo-metrics generate the same topology.

\midvspace

\blemm
\label{lemm pseudo-m comp}
Let $X$ be a set. For $i \in \left\{ 1, 2 \right\}$ let $d_i$ be a pseudo-metric on~$X$ and define the function
\[
B_i : X \!\times \left] 0, \infty \right[\; \longrightarrow  {\mc P}(X), \quad B_i(x,r) = \left\{ y \in X \, : \, d_i(x,y) < r \right\}
\]

That is, for $i \in \left\{ 1, 2 \right\}$, $x \in X$, and $r \in \; ]0,\infty[ \, $, the set $B_i(x,r)$ is the open sphere about $x$ with $d_i\,$-radius~$r$. $\tau (d_1)$~is finer than $\tau (d_2)$ iff for every $x \in X$ and $r \in \; ]0,\infty[ \, $, there is $s \in \; ]0,\infty[ \,$ such that $B_1(x,s) \subset B_2(x,r)$.
\elemm

\bproof
This follows by Lemma~\ref{base finer neighborhood}.
\eproof

The following Corollary is an application of Lemma~\ref{lemm pseudo-m comp}. It shows that for every pseudo-metric space there is a pseudo-metric bounded by $1$ that generates the same topology as the original pseudo-metric.

\midvspace

\bcoro
\label{coro bounded pseudo-m}
Let $(X,d)$ be a pseudo-metric space. We define the functions $e : X \times X \longrightarrow {\mathbb R}_+$ where $e(x,y)$ is the minimum of $d(x,y)$ and~$1$, and
\[
f : X \times X \longrightarrow {\mathbb R}_+, \quad f(x,y) = \frac{d(x,y)}{1 + d(x,y)}
\]

Both $e$ and $f$ are pseudo-metrics on~$X$. All three pseudo-metrics generate the same topology. Further, if $d$ is a metric, then each of the maps $e$ and $f$ is a metric, too.
\ecoro

\bproof
Exercise.
\eproof

We finally obtain the following result about pseudo-metric topologies.

\midvspace

\btheo
\label{theo pseudo-m first countable}
For every pseudo-metric space the generated topology is first countable.
\etheo

\bproof
We may choose $R = {\mathbb D}_+ \!\!\setminus\! \left\{ 0 \right\}$ in Lemma and Definition~\ref{lede pseudo-m topology}. Then $R$~is countable by Corollary~\ref{coro dyadic countable}. Hence there is a countable neighborhood base of~$x$ for every $x \in X$ by Lemma and Definition~\ref{lede pseudo-m topology}.
\eproof


\chapter{Convergence and continuity}
\label{convergence}
\setcounter{equation}{0}

\pagebreak

This Chapter is devoted to the notion of convergence, and to the related topics of limit points, adherence points, continuity of functions, as well as the closure, interior, and boundary of sets. There are three different concepts of convergence, which are interrelated: sequences, nets, and filters. Each of these concepts is treated in one of the following Sections.

\section{Sequences}
\label{seqconv}

The elementary concept of convergence uses sequences as defined in the following Definition.

\midvspace

\bdefi
\label{defseq}
\index{Sequence}
Given a set $X$, a {\bf sequence in} $X$ is a function $x: \naturalnumbers \longrightarrow X$. It is denoted by $(x_n : n \in \naturalnumbers)$, or, if it is known from the context that $x$ is a sequence, by the short notations $(x_n)$ or~$x_n$. The value~$x(n)$ of $x$ at a point $n \in \naturalnumbers$ is also denoted by~$x_n$. 
\edefi

In Chapter~\ref{topologies} we have introduced topological spaces and pseudo-metric spaces. In both cases we may define, using sequences, the notions of convergence, limit points, and adherence points. It is then demonstrated below that the respective definitions agree for a pseudo-metric space and the generated pseudo-metric topology.

\midvspace

\bdefi
\label{defi ev fre sequence}
\index{Eventually}
\index{Frequently}
Let $X$ be a set, $(x_n)$ a sequence in $X$, and $\varphi(x)$ a formula. We say that $\varphi(x_n)$ is true {\bf eventually}, or short $\varphi(x_n)$ {\bf eventually}, if there is $N \!\in \naturalnumbers$ such that $\varphi(x_n)$ is true for every $n \in \naturalnumbers$ with $N \leq n$. We say that $\varphi(x_n)$ is true {\bf frequently}, or short $\varphi(x_n)$ {\bf frequently}, if, for every $N \!\in \naturalnumbers$, there is $n \in \naturalnumbers$, $N \leq n$ such that $\varphi(x_n)$ is true. 
\edefi

Notice that Definition~\ref{defi ev fre sequence} is actually not a single definition but provides one definition for each formula~$\varphi(x)$. In the remainder of the text the notions "eventually" and "frequently" are, however, always used together with specific formulae, not with generic formula variables. Thus Definition~\ref{defi ev fre sequence} should be regarded as an abbreviated form of writing down a finite number of Definitions. This issue is similar to the one mentioned in the context of the Separation schema~\ref{axio sep schema} and the Replacement schema~\ref{axio repl schema}, see the discussion below Axiom~\ref{axio sep schema}.

\midvspace

\bdefi
\label{def seq top conv}
\index{Limit point!sequence}
\index{Sequence!limit point}
\index{Convergent!sequence}
\index{Sequence!convergent}
\index{Adherence point!sequence}
\index{Sequence!adherence point}
Given a topological space $(X,\tp)$, a subset $A \subset X$, and a sequence $(x_n)$ in~$A$, a point $x \in X$ is called {\bf limit point of}~$(x_n)$ if $x_n \!\in U$ eventually for every $U \!\in {\mc N} \!\left\{ x \right\}$. In this case we say that $(x_n)$ {\bf converges} to~$x$, or write $x_n \rightarrow x$. The set of all limit points of~$(x_n)$ is denoted by $\lim_n x_n$. If $(x_n)$ has a limit point, it is called {\bf convergent}. If $(x_n)$ has a unique limit point, say~$x$, we have $\lim_n x_n = \left\{ x \right\}$, and we also write $\lim_n x_n = x$. Further, we say that a point $x \in X$ is an {\bf adherence point of}~$(x_n)$ if $x_n \!\in U$ frequently for every $U \!\in {\mc N} \!\left\{ x \right\}$. The set of all adherence points of~$(x_n)$ is denoted by $\mathrm{adh}_n \, x_n$. If $(x_n)$ has a unique adherence point, say~$x$, we have $\mathrm{adh}_n \, x_n = \left\{ x \right\}$, and we also write $\mathrm{adh}_n \, x_n = x$.
\edefi

\bdefi
\label{def seq pseudo-m conv}
Given a pseudo-metric space $(X,d)$, a subset $A \subset X$, and a sequence $(x_n)$ in~$A$, a point $x \in X$ is called {\bf limit point of} $(x_n)$ if $d(x,x_n) < \varepsilon$ eventually for every $\varepsilon \in {\mathbb R}$, $\varepsilon > 0$. In this case we say that $(x_n)$ {\bf converges} to~$x$, or write $x_n \rightarrow x$.  Further we say that a point $x \in X$ is an {\bf adherence point of}~$(x_n)$ if $d(x,x_n) < \varepsilon$ frequently for every $\varepsilon \in {\mathbb R}$, $\varepsilon > 0$. For the limit points and for the adherence points of~$(x_n)$ we use the same notations as for limit points and adherence points of sequences in a topological space.
\edefi

\brema
Let $(X,d)$ be a pseudo-metric space, $(x_n)$~a sequence in~$X$, and $x \in X$. Since ${\mathbb D}_+$ is dense in~${\mathbb R}_+$, we have:
\benum
\item $x \in \lim_n x_n \;\; \Longleftrightarrow \;\; \big(\, d(x,x_n) < \varepsilon \;\; \mbox{eventually for every} \;\; \varepsilon \in {\mathbb D}_+ \!\!\setminus\! \left\{ 0 \right\} \big)$
\item $x \in \mathrm{adh}_n \, x_n \;\; \Longleftrightarrow \;\; \big(\, d(x,x_n) < \varepsilon \;\; \mbox{frequently for every} \;\; \varepsilon \in {\mathbb D}_+ \!\!\setminus\! \left\{ 0 \right\} \big)$
\eenum

\erema

Obviously a limit point of a sequence in a topological space or in a pseudo-metric space is also an adherence point of the sequence. In general, a sequence $(x_n)$ can have more than one adherence point. A sequence can also have more than one limit point.

\midvspace

\blemm
\label{lemm base conv seq}
Given a topological space $(X,\tp)$ and a sequence $(x_n)$ in~$X$, a point $x \in X$ is a limit point of~$(x_n)$ iff there is a neighborhood base $\mc B$ of~$x$ such that $x_n \!\in B$ eventually for every $B \in {\mc B}$.
\elemm

\bproof
Exercise.
\eproof

The following Lemma ensures that the repeated definition of the same notions in Definitions \ref{def seq top conv} and~\ref{def seq pseudo-m conv} is meaningful.

\midvspace

\blemm
Let $(X,d)$ be a pseudo-metric space, $(x_n)$~a sequence in~$X$, and $x \in X$. Then $(x_n)$ converges to~$x$ with respect to~$d$ in the sense of Definition~\ref{def seq pseudo-m conv} iff it converges to~$x$ with respect to~$\tau(d)$ in the sense of Definition~\ref{def seq top conv}.
\elemm

\bproof
Exercise.
\eproof

\bdefi
\label{defi subsequence}
\index{Subsequence}
Let $X$ and $Y$ be two sets, $(x_n)$~a sequence in~$X$, and $(y_n)$ a sequence in~$Y$. If there exists a strictly increasing map $f : \naturalnumbers \longrightarrow \naturalnumbers$ such that $y_n = x_{f(n)}$ for every $n \in \naturalnumbers$, then $(y_n)$ is called {\bf subsequence of}~$(x_n)$, or short $(y_n) \subset (x_n)$.
\edefi

\btheo
\label{seq adh sub}
Given a pseudo-metric space~$(X,d)$ and a sequence $(x_n)$ in~$X$, $x$~is adherence point of~$(x_n)$ iff there exists a subsequence $(y_n) \subset (x_n)$ such that $y_n \rightarrow x$.
\etheo

\bproof
Assume $x \in \mathrm{adh}_n \, x_n$. We may choose a bijection $f : \naturalnumbers \longrightarrow {\mathbb D}_+ \!\!\setminus\! \left\{ 0 \right\}$ by Corollary~\ref{coro dyadic countable}. We define a map $g : \naturalnumbers \longrightarrow \naturalnumbers$ where $g(m)$ is the minimum of
\[
\big\{ k \in \naturalnumbers \, : \, k \geq m, \, d(x, x_k) < e \big\}
\]

and $e$ is the minimum of the finite set $f \left[ \sigma(m) \right]$. $g$~is clearly unbounded. We further define a map $h : \naturalnumbers \longrightarrow \naturalnumbers$ by Recursion as follows:
\benum
\item $h(0) = g(0)$
\item $h(\sigma(m))$ is the minimum of $\left\{ k \in \mathrm{ran} \, g \, : \, k > h(m) \right\}$ for $m \in \naturalnumbers$
\eenum

We have $\mathrm{ran} \, h = \mathrm{ran} \, g$ by the Induction principle, and $x_{g(m)} \rightarrow x$. It follows that $x_{h(m)} \rightarrow x$. Moreover, $h$~is strictly increasing by definition.

The converse is clear.
\eproof

The following result says that an increasing (decreasing) sequence that is bounded converges, provided certain conditions on the ordering are satisfied.

\midvspace

\blemm
\label{lemm increasing convergent}
Let $(X,<)$ be a totally ordered space where $<$ is an interval relation and has the least upper bound property, $(x_n)$ a bounded sequence in~$X$, and $A = \left\{ x_n \, : \, n \in \naturalnumbers \right\}$. The following statements hold:
\benum
\item \label{lemm increasing convergent 1} If $(x_n)$ is increasing, then $x_n \rightarrow \sup A\,$ with respect to the interval topology.
\item \label{lemm increasing convergent 2} If $(x_n)$ is decreasing, then $x_n \rightarrow \inf A\,$ with respect to the interval topology.
\eenum

\elemm

\bproof
In order to show~(\ref{lemm increasing convergent 1}), assume that $(x_n)$ is increasing. Since $A$ has an upper bound, it has a supremum by assumption. This supremum is unique by Lemma and Definition~\ref{lede ordered sup unique}. We define $x = \sup A$.\\ By Lemmas~\ref{lemm interval top base} and~\ref{lemm top neigborhood base} the system
\[
{\mc A} = \big\{ \left] y, z \right[\; : \, y, z \in X,\; y < x < z \big\}
\]

is a neighborhood base of~$x$. Now let $y, z \in X$ with $y < x < z$. Assume there is no $m \in \naturalnumbers$ such that $x_m \in \,\left] y, z \right[\,$. Then $y$ is an upper bound of~$A$, which is a contradiction. Hence there is $n \in \naturalnumbers$ such that $x_m \in \,\left] y, z \right[\,$ for every $m \in \naturalnumbers$ with $m \geq n$ since $(x_n)$ is increasing.

To see~(\ref{lemm increasing convergent 2}), notice that $A$ has a greatest lower bound if it has a lower bound, by Theorem~\ref{theo least upper greatest lower}. The remainder of the proof is similar to that of~(\ref{lemm increasing convergent 1}).
\eproof

\brema
\label{rema r increasing conv}
Notice that $({\mathbb R},<)$ satisfies the conditions of Lemma~\ref{lemm increasing convergent} by Lemmas~\ref{rema r interval rel} and~\ref{lemm least upper bound reals}.
\erema

\blemm
\label{lemm sequ inv conv}
Let $(x_n)$ be a sequence in ${\mathbb R}_+ \!\!\setminus\! \left\{ 0 \right\}$. If $(x_n)$ is increasing and unbounded, then $x_n^{-1} \rightarrow 0$ with respect to the standard topology on~${\mathbb R}$.
\elemm

\bproof
Let $y \in {\mathbb R}_+  \!\!\setminus\! \left\{ 0 \right\}$. Under the stated conditions, there is $m \in \naturalnumbers$ such that $x_n > y^{-1}$ for every $n \in \naturalnumbers$ with $n > m$. Hence $x_n^{-1} < y$ for $n \in \naturalnumbers$, $n > m$ by Corollary~\ref{coro inequality inv} and Remark~\ref{rema r interval rel}.
\eproof

We may compare the convergence properties of a sequence with respect to two different topologies on the same set. To this end we introduce the following notions.

\midvspace

\bdefi
\index{Sequentially stronger}
\index{Topology!sequentially stronger}
\index{Sequentially equivalent}
\index{Topology!sequentially equivalent}
Let $X$ be a set and $\tp_1$ and $\tp_2$ be two topologies on~$X$. Then $\tp_1$ is called {\bf sequentially stronger than} $\tp_2$ if $x_n \rightarrow x$ with respect to~$\tp_1$ implies $x_n \rightarrow x$ with respect to~$\tp_2$ for every sequence $(x_n)$ in~$X$ and every $x \in X$. $\tp_1$~and $\tp_2$ are called {\bf sequentially equivalent} if $\tp_1$ is sequentially stronger than $\tp_2$ and vice versa.
\edefi

The following is an intuitive result that relates the comparison of two topologies on the same set to the comparison of sequence convergence.

\midvspace

\blemm
\label{lemm finer seq stronger}
Let $X$ be a set and $\tp_1$ and $\tp_2$ be two topologies on~$X$. If $\tp_1$ is finer than $\tp_2$, then it is also sequentially stronger.
\elemm

\bproof
Assume that $\tp_1$ is finer than $\tp_2$ and let $(x_n)$ be a sequence in $X$, $x \in X$, and $x_n \rightarrow x$ with respect to~$\tp_1$. For $i \in \left\{ 1, 2 \right\}$ let ${\mc N}_i$ be the neighborhood system of~$(X,\tp_i)$. We have ${\mc N}_2 \! \left\{ x \right\} \subset {\mc N}_1 \! \left\{ x \right\}$ by Lemma~\ref{base finer neighborhood}. Thus $x_n \rightarrow x$ with respect to~$\tp_2$.
\eproof

The condition under which the converse is true is proven in Section~\ref{closure} because a result about the closure of sets is required that is not derived before.

\midvspace

\blede
Let $X$ and $Y$ be two sets, $f : X \longrightarrow Y$ a map, $A \subset X$, and $(x_n)$ a sequence in~$A$. Then $(f(x_n) : n \in \naturalnumbers)$ is a sequence in~$Y$, also denoted by~$(f(x_n))$ or~$f(x_n)$.
\elede

\bproof
This is clear.
\eproof

\section{Nets}
\label{nets}

In this Section the concept of sequences is generalized by the introduction of nets.

\midvspace

\bdefi
\label{defi net}
\index{Net}
Given a set $X$ and a directed space $(D,\leq)$, a function $x: D \longrightarrow X$ is called a {\bf net in}~$X$. It is denoted by $(x_n : n \in D)$, or, if it is known from the context that $x$ is a net with some domain~$D$, by the short notations $(x_n)$ or~$x_n$. The value~$x(n)$ of~$x$ at a point $n \in D$ is also denoted by~$x_n$. 
\edefi

Note that $(\naturalnumbers,\leq)$ is a directed space. Therefore, given a set~$X$, the nets in~$X$ that are of the form $(x_n : n \in \naturalnumbers)$, i.e.\ those where the directed set is~$\naturalnumbers$, are precisely the sequences in~$X$. In this case the definitions of the notations $(x_n : n \in \naturalnumbers)$, $(x_n)$, and $x_n$ for nets in Definition~\ref{defi net} agree with the respective definitions for sequences in Definition~\ref{defseq}. Many of the illustrative properties of sequences generalize to similar properties of nets.

\midvspace

\bdefi
\label{defi ev fre nets}
\index{Eventually}
\index{Frequently}
Let $X$ be a set, $(x_n : n \in D)$ a net in~$X$, and $\varphi(x)$ a formula. We say that $\varphi(x_n)$ is true {\bf eventually}, or short $\varphi(x_n)$ {\bf eventually}, if there is $N \!\in \naturalnumbers$ such that $\varphi(x_n)$ is true for every $n \in \naturalnumbers$ with $N \leq n$. We say that $\varphi(x_n)$ is true {\bf frequently}, or short $\varphi(x_n)$ {\bf frequently}, if, for every $N \!\in \naturalnumbers$, there is $n \in \naturalnumbers$ with $N \leq n$ such that $\varphi(x_n)$ is true. 
\edefi

Notice that in the case $D = {\mathbb N}$ Definition~\ref{defi ev fre nets} agrees with Definition~\ref{defi ev fre sequence} where the same notions are defined for sequences. Regarding the usage of formula variables in Definition~\ref{defi ev fre nets} the same remarks apply as with respect to Definition~\ref{defi ev fre sequence}.

\midvspace

\bdefi
\label{defi conv nets}
\index{Convergent!net}
\index{Limit point!net}
\index{Net!limit point}
\index{Adherence point!net}
\index{Net!adherence point}
Given a topological space $(X,\tp)$, a subset $A \subset X$, and a net $(x_n : n \in D)$ in~$A$, a point $x \in X$ is called {\bf limit point of}~$(x_n)$ if $x_n \!\in U$ eventually for every $U \!\in {\mc N} \! \left\{ x \right\}$. In this case we say that $(x_n)$ {\bf converges to}~$x$, or write $x_n \rightarrow x$. The set of all limit points of $(x_n)$ is denoted by $\lim_n x_n$. If $(x_n)$ has a limit point, it is called {\bf convergent}. If $(x_n)$ has a unique limit point, say~$x$, we also write $\lim_n x_n = x$. Further, we say that a point $x \in X$ is an {\bf adherence point of} $(x_n)$ if $x_n \!\in U$ frequently for every $U \!\in {\mc N} \! \left\{ x \right\}$. The set of all adherence points of~$(x_n)$ is denoted by $\mathrm{adh}_n \, x_n$. If $(x_n)$ has a unique adherence point, say~$x$, we also write $\mathrm{adh}_n \, x_n = x$.
\edefi

Notice that every limit point of a net $(x_n)$ is also an adherence point. A net can have more than one limit point. We remark that in the case $D = {\mathbb N}$ Definition~\ref{defi conv nets} agrees with the Definition~\ref{def seq top conv}.

\midvspace

\bexam
\label{exneighb}
Let $(X,\tp)$ be a topological space, $x \in X$, and $\leq$ the relation on ${\mc N} \!\left\{ x \right\}$ that is defined by
\[
U \leq V \quad \Longleftrightarrow \quad V \subset U
\]

Then $\left({\mc N} \!\left\{ x \right\}, \leq \right)$ is a directed space (which justifies the notation~$\leq$). We may choose a point $x_U \in X$ for each $U \!\in {\mc N} \!\left\{ x \right\}$. Then $(x_U : U \!\in {\mc N} \!\left\{ x \right\})$ is a net in~$X$ and $x_U \rightarrow x$.

Similarly, if $\mc B$ is a neighborhood base of~$x$ and the relation $\leq$ on~$\mc B$ is defined as above for $U, V \in {\mc B}$, then $\left({\mc B}, \leq \right)$ is a directed space. We may choose a point $x_B \in X$ for each $B \in {\mc B}$. Then $(x_B : B \in {\mc B})$ is a net in~$X$ and $x_B \rightarrow x$.
\eexam

The analogue of Lemma~\ref{lemm base conv seq} also holds for nets.

\midvspace

\blemm
Given a topological space $(X,\tp)$ and a net $(x_n : n \in D)$ in~$X$, a point $x \in X$ is a limit point of $(x_n)$ iff there is a neighborhood base $\mc B$ of~$x$ such that $x_n \in B$ eventually for every $B \in {\mc B}$.
\elemm

\bproof
Exercise.
\eproof

In the case of a pseudo-metric topology we have the following characterization of convergence.

\midvspace

\blemm
\label{lemm pseudom conv}
Let $(X,d)$ be a pseudo-metric space, $(x_n : n \in D)$ a net in $X$, and $x \in X$. $(x_n)$~converges to~$x$ with respect to~$\tau(d)$ iff the net $\left( d(x_n,x) : n \in D \right)$ converges to~$0$ with respect to the standard topology on~${\mathbb R}_+$.
\elemm

\bproof
Exercise.
\eproof

\brema
\label{rema equiv net r}
Let $(x_n)$ be a net in~${\mathbb R}$ and $x \in {\mathbb R}$. Then we have:
\[
x_n \rightarrow x \quad \Longleftrightarrow \quad |x_n - x| \rightarrow 0
\]

where the limit on the left-hand side is with respect to the standard topology on~$\mathbb R$ and the limit on the right-hand side is with respect to the standard topology on~${\mathbb R}_+$.
\erema

We now introduce subnets, which are the analogue of subsequences as defined in Definition~\ref{defi subsequence}.

\midvspace

\bdefi
\index{Subnet}
Given two sets $X$ and $Y$, and a net $(x_n : n \in D)$ in~$X$, a net $(y_m : m \in E)$ in~$Y$ is called {\bf subnet of}~$(x_n)$ if there is a map $f: E \longrightarrow D$ such that
\benum
\item $\forall m \in E \quad y_m = x_{f(m)}$
\item $\forall n \in D \quad \exists m \in E \quad \forall p \in E \quad p \geq m \;\; \Longrightarrow \;\; f(p) \geq n$
\eenum

We also use the short notation $(y_m) \subset (x_n)$.
\edefi

\blemm
\label{lemm lim adh subnet}
Given a topological space $(X,\tp)$, a net $(x_n)$ in~$X$, and a subnet $(y_m) \subset (x_n)$, the following statements hold:
\benum
\item $\lim_n x_n \subset \, \lim_m y_m$
\item $\mathrm{adh}_m \, y_m \subset \, \mathrm{adh}_n \, x_n$
\eenum

\elemm

\bproof
Exercise.
\eproof

Lemma and Definition~\ref{product directed set} can be used to derive the analogue of Theorem~\ref{seq adh sub} for nets.

\midvspace

\btheo
\label{net adh sub}
Given a topological space $(X,\tp)$ and a net $(x_n : n \in D)$ in~$X$, a point $x \in X$ is adherence point of~$(x_n)$ iff there is a subnet $(y_r : r \in E) \subset (x_n)$ such that $y_r \rightarrow x$.
\etheo

\bproof
Assume that $x \in \mathrm{adh}_n \, x_n$. Let the relation $\leq$ on ${\mc N} \!\left\{ x \right\}$ be defined as in Example~\ref{exneighb}. We define
\[
E = \big\{ (n,U) \in D \times {\mc N} \!\left\{ x \right\} \, : \, x_n \in U \big\}
\]

and the relation $\leq$ on~$E$ as the restriction to~$E$ of the relation introduced in Lemma and Definition~\ref{product directed set}. Then $(E,\leq)$ is a directed space.
\\

\hspace{0.05\textwidth}
\parbox{0.95\textwidth}
{[The relation $\leq$ is clearly transitive and reflexive. Let $(m,U), (n,V) \in E$, and $W = U \cap V$. Then we have $W \!\in {\mc N} \!\left\{ x \right\}$, and there is $k \in \naturalnumbers$ such that $k \geq m, n$ and $x_k \in W$. It follows that $\leq$ is a directive relation on~$E$.]}
\\

Further let $p : D \times {\mc N} \!\left\{ x \right\} \longrightarrow D$ be the projection on the first component, and $(y_r : r \in E)$ the net in~$X$ such that $y = x \circ p$. Then $(y_r)$ clearly is a subnet of~$(x_n)$. Furthermore, $y_r \rightarrow x$.

The converse is clear.
\eproof

The following Theorem states the existence of a diagonal net that has among its limit points the limit points of iterated limits.

\midvspace

\btheo[Iterated limits]
\label{theo iterated limits}
\index{Iterated limits}
Let $(X,\tp)$ be a topological space and $(D,\leq)$ a directed space. Further let $(E_m,\leq)$ be a directed space for every $m \in D$ and $F = \bigcup_{m \in D} \left\{ m \right\} \times E_m$. Then there exists a net $(T_r : r \in R)$ in~$F$ such that for every map $S: F \longrightarrow X$ and every map $S' : D \longrightarrow X$ with $S'(m) \in \lim_n S(m,n)$ ($m \in D$) we have $\lim_m S'(m) \subset \lim_r S \circ T (r)$. If ${\mc N}^{\, \mathrm{closed}}$ is a neighborhood base, then equality holds instead of "$\subset$".
\etheo

\bproof
We define the product directed space $(R,\leq)$ where $R = D \!\times\! P$ and $P = \bigtimes_{\!\! m \in D}\, E_m$, the projections $p_m : P \longrightarrow E_m$ ($m \in D$), and the net $T : R \longrightarrow F$, $T(m,e) = \left(m, p_m(e)\right)$ for $m \in D$ and $e \in P$. Now let $S$ and $S'$ be two functions satisfying the stated conditions.

If $S'' \!\in \lim_m S'(m)$, and $U \!\in {\mc N}^{\, \mathrm{open}} \left\{ S'' \right\}$, then there is $M \!\in D$ such that $S'(m) \in U$ for every $m \in D$ with $m \geq M$. We may choose $e \in P$ such that $S(m,n) \in U$ for every $m \in D$ with $m \geq M$ and every $n \in E_m$ with $n \geq p_m(e)$. Hence, $S \circ T (m,e') \in U$ for every $(m,e') \in R$ with $(m,e') \geq (M,e)$.

Conversely, assume that ${\mc N}^{\, \mathrm{closed}}$ is a neighborhood base and let $S'' \!\in \lim_r S \circ T (r)$. Further let $U \!\in {\mc N}^{\, \mathrm{closed}} \left\{ S'' \right\}$. There are $M \!\in D$, $e \in P$ such that $S(m,p_m(e')) \in U$ for every $(m,e') \in R$ with $(m,e') \geq (M,e)$. Then, for every $m \in D$ with $m \geq M$, we have $S(m,n) \in U$ for $n \in E_m$, $n \geq p_m(e)$. It follows that $\lim_n S(m,n) \subset U$ for every $m \in D$, $m \geq M$. Thus $S'' \!\in \lim_m S'(m)$.
\eproof

\blede
Let $X$ and $Y$ be sets, $f : X \longrightarrow Y$ a map, $A \subset X$, and $(x_n : n \in D)$ a net in~$A$. Then $(f(x_n) : n \in D)$ is a net in~$Y$, and also denoted by $(f(x_n))$ or $f(x_n)$.

\elede

\bproof
This is clear.
\eproof

The following construction is useful in the context of product spaces.

\midvspace

\bprop
\label{prop net prod dir}
Let $(D_i,R_i)$ ($i \in I$) be directed sets where $I$ is an index set, $(D,R)$~the product directed set, and $p_i : D \longrightarrow D_i$ ($i \in I$) the projections. Further let $(X,\tp)$ be a topological space, $A \subset X$, $j \in I$, and $x : D_j \longrightarrow A$ a net in~$A$. The net $y : D \longrightarrow A$, $y(r) = x \big( p_j(r) \big)$, is a subnet of~$(x_n)$. Moreover we have
\[
\mathrm{ran} \, y = \mathrm{ran} \, x\,, \quad \quad {\textstyle \lim_n} \, x_n = \, {\textstyle \lim_r} \, y_r\,, \quad \quad \mathrm{adh}_n \, x_n = \, \mathrm{adh}_r \, y_r
\]

\eprop

\bproof
Exercise.
\eproof

\bdefi
\label{defi product net}
\index{Product net}
\index{Net!product}
Let $I$ be an index set. For every $i \in I$ let $(D_i,R_i)$ be a directed set, $(X_i,\tp_i)$ be a topological space, $A_i \subset X_i$\,, and $(x^i_n : n \in D_i)$ a net in~$A_i$\,. Let $(D,R)$ be the product directed set and $p_i : D \longrightarrow D_i$ ($i \in I$) the projections. Further let $X = \bigtimes_{\!\! i \in I}\, X_i$ and $q_i : X \longrightarrow X_i$ ($i \in I$) be the corresponding projections. Moreover let $A = \bigtimes_{\!\! i \in I}\, A_i$\,. The net $x : D \longrightarrow A$ defined by $q_i \big( x(r) \big) = x^i \big( p_i(r) \big)$ for every $i \in I$ and $r \in D$ is called {\bf product net} and denoted by~$\prod_{i \in I}\, (x^i_n)$. If $I = \sigma(n) \!\setminus\! m$ for some $m, n \in \naturalnumbers$ with $m < n$, then we also write $\prod_{k = m}^n (x^i_k)$ for the product net. If $I = \naturalnumbers \!\setminus\! m$ for some $m \in \naturalnumbers$, then we also write $\prod_{k = m}^{\infty} (x^i_k)$ for the product net.
\edefi

\blemm
\label{lemm prod net conv}
With definitions as in Definition~\ref{defi product net} we have for every $i \in I$:
\[
{\textstyle \lim_r} \, q_i(x_r) = \, {\textstyle \lim_n} \, x^i_n\,, \quad \quad \mathrm{adh}_r \, q_i(x_r) = \, \mathrm{adh}_n \, x^i_n
\]

\elemm

\bproof
This follows by Propositon~\ref{prop net prod dir}.
\eproof

\section{Filters}
\label{filters conv}

We have introduced the notion of filter in Section~\ref{filters}. The concept of filter can be widely used as an alternative to nets when probing convergence and related properties. Although in many cases the usage of nets seems more illustrative, there are cases where filters are advantageous.

\midvspace

\bdefi
\index{Convergent!filter}
\index{Limit point!filter}
\index{Filter!limit point}
\index{Adherence point!filter}
\index{Filter!adherence point}
Given a topological space $(X,\tp)$ and a filter $\mc F$ on~$X$, a point $x \in X$ is called {\bf limit point of}~$\mc F$ if ${\mc N} \!\left\{ x \right\} \subset {\mc F}$. In this case we say that $\mc F$ {\bf converges to}~$x$, written ${\mc F} \rightarrow x$. If $\mc F$ has a limit point, it is called {\bf convergent}. The set of all limit points of~$\mc F$ is denoted by $\lim {\mc F}$. Further, a point $x \in X$ is called an {\bf adherence point of}~$\mc F$ if $U \cap F \neq \O$ for every $U \!\in {\mc N} \!\left\{ x \right\}$ and every $F \!\in {\mc F}$. The set of all adherence points of~$\mc F$ is denoted by $\mathrm{adh} \, {\mc F}$. Let ${\mc B} \subset {\mc F}$ be a filter base for~$\mc F$. A point $x \in X$ is called {\bf limit point of}~$\mc B$ if it is a limit point of~$\mc F$. In this case we say that $\mc B$ {\bf converges to}~$x$, written ${\mc B} \rightarrow x$. A point $x \in X$~is called {\bf adherence point of}~$\mc B$ if it is an adherence point of~$\mc F$. The set of all limit (adherence) points of~$\mc B$ is denoted by $\lim {\mc B}$ ($\mathrm{adh} \, {\mc B}$). 
\edefi

\bexam
Given a topological space $(X,\tp)$ and a point $x \in X$, we obviously have \mbox{${\mc N} \!\left\{ x \right\} \rightarrow x$}.
\eexam

\bexam
The filter in Example~\ref{exfrech} has no adherence points.
\eexam

We now establish a connection between sequences and filters, similarly to that between sequences and nets in Section~\ref{nets}. The agreement of the respective limit and adherence points is shown.

\midvspace

\blede
\label{lede filter generated}
Given a set $X$ and a sequence $(x_n)$ in~$X$, the system
\[
{\mc F} = \big\{F \subset X \, : \, x_n \in F \;\, \mathrm{eventually} \big\}
\]

is a filter. We say that the sequence {\bf generates}~$\mc F$. The system ${\mc B} = \left\{ B_n \, : \, n \in \naturalnumbers \right\}$ where $B_n = \left\{ x_k \, : \, k \geq n \right\}$ ($n \in \naturalnumbers$) is a filter base for~$\mc F$. We also say that the sequence {\bf generates}~$\mc B$.
\elede

\bproof
$\mc B$ is a filter base by Lemma and Definition~\ref{lemm characterization filter base}. Clearly, ${\mc F} = \Phi ({\mc B})$.
\eproof

\blemm
Let $(X,\tp)$ be a topological space, $(x_n)$ a sequence in $X$, and $\mc F$ the filter generated by~$(x_n)$. We have \,\! $\lim_n x_n = \lim {\mc F}$ \,\! and \,\! $\mathrm{adh}_n \, x_n = \mathrm{adh} \, {\mc F}$.
\elemm

\bproof
Exercise.
\eproof

For a given set~$X$, Lemma and Definition~\ref{lede filter generated} introduces an injective map from the system of all sequences in~$X$ to the system of all filters on~$X$, which is clearly not bijective in general.

Next, two types of connections between nets and filters on a set~$X$ are established and the relationship between their limit and adherence points is investigated. In the first case a net, the so-called associated net, is defined for a given filter. However, the limit points and the adherence points of the filter and the associated net are generally not the same.

\midvspace

\blede
\label{filternet}
\index{Associated}
\index{Net!associated}
Let $X$ be a set and denote by $\leq$ the relation $\supset$ on~${\mc P}(X)$. Further let $\mc F$ be a filter on~$X$ and $\mc B$ a filter base for~$\mc F$. Then both $({\mc B},\leq)$ and $({\mc F},\leq)$ are directed spaces. We may choose a point $x_F \in F$ for each $F \in {\mc F}$. Then $(x_F : F \in {\mc F})$ is a net in~$X$. Similarly, we may choose a point $x_B \in B$ for each $B \in {\mc B}$. Then $(x_B : B \in {\mc B})$ is a net in~$X$. If a net in~$X$ can be constructed in this way, it is called {\bf associated with} $\mc F$ or $\mc B$, respectively.
\elede

\bproof
This is obvious.
\eproof

\blemm
\label{conv assoc net}
Let $(X,\tp)$ be a topological space, $\mc F$~a filter on~$X$, $\mc B$~a filter base for~$\mc F$, and $(x_F : F \in {\mc F})$ and $(x_B : B \in {\mc B})$ nets associated with $\mc F$ and $\mc B$, respectively. Then the following two statements hold:

\benum
\item \label{conv assoc net 1} $\lim {\mc F} \subset \lim_F x_F$
\item \label{conv assoc net 2} $\lim {\mc F} \subset \lim_B x_B$
\item \label{conv assoc net 3} $\mathrm{adh}_F \, x_F \subset \mathrm{adh} \, {\mc F}$
\item \label{conv assoc net 4} $\mathrm{adh}_B \, x_B \subset \mathrm{adh} \, {\mc F}$
\eenum

\elemm

\bproof
To see~(\ref{conv assoc net 1}), let $x \in \lim {\mc F}$ and $U \!\in {\mc N} \!\left\{ x \right\}$. Since $U \!\in {\mc F}$, it follows that $x_V \in U$ for $V \!\in {\mc F}$, $V \geq U$.

To show~(\ref{conv assoc net 2}), let $x \in \lim {\mc F}$ and $U \!\in {\mc N} \!\left\{ x \right\}$. Since $U \!\in {\mc F}$, there is $B \in {\mc B}$ such that $B \subset U$. It follows that $x_A \in U$ for $A \in {\mc B}$, $A \geq B$.

In order to prove~(\ref{conv assoc net 3}), let $x \in \mathrm{adh}_F \, x_F$, $F \in {\mc F}$, and $U \!\in {\mc N} \!\left\{ x \right\}$. Then there is $x_G$ with $G \geq F$ such that $x_G \in U$. Since we also have $x_G \in F$, it follows that $F \cap U \neq \O$.

To see~(\ref{conv assoc net 4}), let $x \in \mathrm{adh}_B \, x_B$, $F \in {\mc F}$, and $U \!\in {\mc N} \!\left\{ x \right\}$. We may choose $B \in {\mc B}$ such that $B \subset F$. Then there is $x_A$ with $A \geq B$ such that $x_A \in U$. Since we also have $x_A \in F$, it follows that $F \cap U \neq \O$.
\eproof

We now show how a filter on~$X$ can be defined for a given net in~$X$, and a net in~$X$ can be defined for a given filter on~$X$, such that when starting with any filter and performing both steps the original filter is obtained. Note that generally there does not exist any bijection between all nets in~$X$ and all filters on~$X$. In fact, there exists no set that contains every directed space, and therefore there is no set that contains every net in~$X$. However, limit points and adherence points agree between corresponding nets and filters in this approach, so a true equivalence of the two concepts is demonstrated.
\midvspace

\blede
\label{lede gen filter}
\index{Filter!generated}
\index{Filter base!generated}
Given a set $X$ and a net $(x_n : n \in D)$ in $X$, the system
\[
{\mc F} = \big\{ F \subset X \, : \, x_n \in F \;\, \mathrm{eventually} \big\}
\]

is a filter. We say that $(x_n)$ {\bf generates}~$\mc F$. The system ${\mc B} = \left\{ B_n \, : \, n \in D \right\}$ where $B_n = \left\{ x_k \, : \, k \geq n \right\}$ ($n \in D$) is a filter base for~$\mc F$. We also say that $(x_n)$ {\bf generates}~$\mc B$.
\elede

\bproof
$\mc B$ is a filter base by Lemma and Definition~\ref{lemm characterization filter base}. Clearly, we have ${\mc F} = \Phi ({\mc B})$.
\eproof

Notice that Lemma and Definition~\ref{lede gen filter} does not refer to any topology on~$X$.

\midvspace

\blemm
\label{lemm filter gen conv}
Let $(X,\tp)$ be a topological space, $(x_n)$~a net in~$X$, $\mc F$~the filter generated by~$(x_n)$, and $x \in X$. We have \,\! $\lim_n x_n = \lim {\mc F}$ \,\! and \,\! $\mathrm{adh}_n \, x_n = \mathrm{adh} \, {\mc F}$.
\elemm

\bproof
Exercise.
\eproof

\blede
\label{lede net gen}
\index{Net!generated}
Let $X$ be a set, $\mc F$~a filter on~$X$, $\mc B$~a filter base for~$\mc F$, and
\[
D = \big\{ (x,B) \, : \, x \in B \in {\mc B} \big\}
\]

Let the relation $\leq$ on~$D$ be defined as follows:
\[
(x,B) \leq (y,C) \quad \Longleftrightarrow \quad C \subset B
\]

Then $(D,\leq)$ is a directed space (which justifies the notation). The net \mbox{$(x_n : n \in D)$} in~$X$ where $x_{(x,B)} = x$ for every $(x,B) \in D$ is called {\bf generated by}~${\mc B}$.

The filter generated by $(x_n)$ is~$\mc F$.
\elede

\bproof
$(D,\leq)$ clearly is a directed space.

Let $\mc G$ be the filter generated by~$(x_n)$. 

To see that ${\mc G} \subset {\mc F}$, let $G \in {\mc G}$. Then there is $n \in D$ such that $x_k \in G$ for $k \in D,\, k \geq n$. Moreover there is $B \in {\mc B}$ such that $n = (x_n, B)$. Then we have $x_k \in B$ for $k \in D,\, k \geq n$. Now let $x \in B$. Then $(x,B) \geq (x_n,B)$, and therefore $x = x_{(x,B)} \in G$. Thus we have $B \subset G$. It follows that $G \in {\mc F}$.

Conversely, let $B \in {\mc B}$. We may choose $x \in B$ and define $n = (x, B)$. It follows that $x_k \in B$ for $k \in D,\, k \geq n$. Thus we obtain $B \in {\mc G}$.
\eproof

Notice that also Lemma and Definition~\ref{lede net gen} does not refer to any topology on~$X$.

The following is the analogue to Theorems~\ref{seq adh sub} and~\ref{net adh sub} where similar results for sequences and nets are proven, respectively.

\midvspace

\btheo
Let $(X,\tp)$ be a topological space and $\mc F$ a filter on~$X$. A point $x \in X$ is an adherence point of~$\mc F$ iff there is a filter $\mc G$ finer than $\mc F$ such that ${\mc G} \rightarrow x$.
\etheo

\bproof
First let $x \in \mathrm{adh} \,{\mc F}$. Then the system
\[
{\mc G} = \big\{ F \cap U \, : \, F \in {\mc F},\; U \in {\mc N} \!\left\{ x \right\} \! \big\}
\]

is a filter. Moreover, we have ${\mc F} \subset {\mc G}$ and ${\mc G} \rightarrow x$.

Conversely, if there exists a filter $\mc G$ finer than $\mc F$ with ${\mc G} \rightarrow x$, then ${\mc N} \!\left\{ x \right\} \subset {\mc G}$. It follows that $F \cap U \neq \O$ for every $F \in {\mc F}$ and $U \in {\mc N} \!\left\{ x \right\}$.
\eproof

We now briefly analyse how filter bases behave under mappings und thereby introduce the notions of image filter and inverse image filter.

\midvspace

\blede
\label{lede image filter}
\index{Image filter}
\index{Filter!image}
Let $X$ and $Y$ be two sets, ${\mc B}$ a filter base on~$X$, and $f: X \longrightarrow Y$ a map. Then $f \, \llbracket {\mc B} \, \rrbracket$ is a filter base on~$Y$. The generated filter $\Phi \left( f \, \llbracket {\mc B} \, \rrbracket \right)$ is called {\bf image filter of} ${\mc B}$ {\bf under}~$f$.
\elede

\bproof
Exercise.
\eproof

Notice that, even if $\mc B$ in Lemma and Definition~\ref{lede image filter} is a filter on~$X$, the image $f \, \llbracket {\mc B} \, \rrbracket$ need not be a filter on~$Y$.

\midvspace

\blemm
\label{lemm image filter base id}
With definitions as in Lemma and Definition~\ref{lede image filter} we have $\Phi \left( f \, \llbracket {\mc B} \, \rrbracket \right) = \Phi \left( f \, \llbracket \Phi ({\mc B}) \rrbracket \right)$, that is the image filter of a filter base and of its generated filter are the same.
\elemm

\bproof
Exercise.
\eproof

\blede
\index{Inverse image filter}
\index{Filter!inverse image}
Let $X$ and $Y$ be two sets, $\mc B$~a filter base on~$Y$, and $f : X \longrightarrow Y$ a surjective map. Then $f^{-1} \, \llbracket {\mc B} \, \rrbracket$ is a filter base for a filter on~$X$. The generated filter
$\Phi \left( f^{-1} \, \llbracket {\mc B} \, \rrbracket \right)$ is called {\bf inverse image filter of}~$\mc B$ {\bf under}~$f$.
\elede

\bproof
Exercise.
\eproof

Again, $f^{-1} \, \llbracket {\mc B} \, \rrbracket$ need not be a filter on~$X$, even if $\mc B$ is a filter on~$Y$.

\section{Continuous functions}
\label{continuous functions}

The topic of this Section are continuous functions. We introduce three types of continuity and then show how they are interrelated. The first one does not refer to any topology or pseudo-metric. It is based on two filters, one on the domain of the function and the other on its range. Moreover, it is a local definition, i.e.\ it refers to a single point of the domain. The second definition is based on two topologies, one on the domain and the other on the range. There we introduce both continuity in a point and continuity of the whole function. The third type of continuity is introduced in the context of pseudo-metric spaces and a priori does not refer to any topology. Also in this case a local and global form of continuity is defined. It is then shown that continuity with respect to pseudo-metrics is equivalent to continuity with respect to the generated topologies. Finally we analyse the special case of interval topologies.

\midvspace

\blede
\label{defcontfilter}
\index{Continuous!filter}
Let $X$, $Y$ be two sets, $f : X \longrightarrow Y$ a function, $x \in X$, $y = f(x)$, $\mc F_x$ a filter on~$X$ such that $x$ is a cluster point of~$\mc F_x$\,, and $\mc F_y$ a filter on~$Y$ such that $y$ is a cluster point of~$\mc F_y$\,. $f$~is called ${\mc F}_x\,$-${\mc F}_y\,$-{\bf continuous in}~$x$ if \,\! ${\mc F_y} \subset_{\stackrel{}{\Phi}} f \, \llbracket {\mc F_x} \rrbracket$\,.

Further, let $\mc B_x$ and $\mc B_y$ be filter bases for $\mc F_x$ and $\mc F_y$, respectively. Then $f$ is ${\mc F}_x\,$-${\mc F}_y\,$-continuous in $x$ iff \,\! ${\mc B_y} \subset_{\stackrel{}{\Phi}} f \, \llbracket {\mc B_x} \rrbracket$\,.
\elede

\bproof
Exercise.
\eproof

\bdefi
\index{Continuous!topological space}
Given two topological spaces $(X,\tp_X)$, $(Y,\tp_Y)$, a map $f : X \longrightarrow Y$, and a point $x \in X$, $f$ is called $\tp_X\,$-$\tp_Y\,$-{\bf continuous in} $x$, or short {\bf continuous in}~$x$, if $f^{-1} \, \llbracket \, \tp_Y \!\left(f(x)\right) \, \rrbracket \subset \tp_X$\,.
\edefi

The following Theorem shows how the two types of continuity are related and provides various characterizations of continuity of a function in a point.

\midvspace

\btheo
\label{theo local cont}
Given two topological spaces $(X,\tp_X)$ and~$(Y,\tp_Y)$, a map $f: X \longrightarrow Y$, a subbase ${\mc S}_Y$ for $\tp_Y$, a point $x \in X$, and a neighborhood base ${\mc E}$ of~$x$, the following statements are equivalent:
\benum
\item \label{theo local cont 1} $f$ is continuous in~$x$.
\item \label{theo local cont 2} $f^{-1} \, \llbracket \, {\mc S}_Y \!\left( f(x) \right) \, \rrbracket \, \subset \, \tp_X$
\item \label{theo local cont 2a} ${\mc S}_Y \!\left( f(x) \right) \, \subset_{\stackrel{}{\Phi}} f \, \llbracket {\mc E} \, \rrbracket$
\item \label{theo local cont 3} $f^{-1} \, \llbracket \, {\mc N} \!\left\{ f(x) \right\} \rrbracket \, \subset \, {\mc N} \!\left\{ x \right\}$
\item \label{theo local cont 4} ${\mc N} \!\left\{ f(x) \right\} \, \subset_{\stackrel{}{\Phi}} f \, \llbracket \, {\mc N} \!\left\{ x \right\} \rrbracket$\,, i.e.\ $f$ is ${\mc N} \!\left\{ x \right\}$-${\mc N} \!\left\{ f(x) \right\}$-continuous in~$x$.
\item \label{theo local cont 5} For every filter base $\mc B$ on~$X$, ${\mc B} \rightarrow x$ implies $f \, \llbracket {\mc B} \, \rrbracket \rightarrow f(x)$.
\item \label{theo local cont 6} For every net $(x_n)$ in~$X$, $x_n \rightarrow x$ implies $f(x_n) \rightarrow f(x)$.
\eenum

\etheo

\bproof
We show the equivalence of (\ref{theo local cont 1}) to~(\ref{theo local cont 3}), and then the equivalence of (\ref{theo local cont 3}) to~(\ref{theo local cont 6}).

The implication from (\ref{theo local cont 1}) to~(\ref{theo local cont 2}) is clear.

To show that(\ref{theo local cont 2}) implies~(\ref{theo local cont 2a}), let $S \in {\mc S}_Y \!\left( f(x) \right)$. Then we have $f^{-1} \left[ S \right] \in \tp_X$ by assumption. We may choose $B \in {\mc E}$ such that $B \subset f^{-1} \left[ S \right]$. It follows that
\[
f \left[ B \right] \, \subset \, f \left[ f^{-1} \left[ S \right] \right] \, \subset \, S
\]

To see that (\ref{theo local cont 2a}) implies (\ref{theo local cont 3}), let $U \in {\mc N} \!\left\{ f(x) \right\}$. There are $S_i \in {\mc S}_Y$ ($i \in I$), where $I$ is a finite index set, such that $f(x) \in S_i$ for every $i \in I$ and $S \subset U$ where $S = \bigcap_{i \in I} S_i$. For every $i \in I$, we may choose $B_i \in {\mc E}$ such that $f \left[ B_i \right] \subset S_i$ by assumption. We define $B = \bigcap_{i \in I} B_i$. It follows that $x \in B \subset f^{-1} \left[ U \right]$.

To show that (\ref{theo local cont 3}) implies (\ref{theo local cont 1}), let $U \in \tp_Y$ with $f(x) \in U$. Further let $z \in f^{-1} \left[ U \right]$. By Lemma~\ref{lemm characterization neighborhood}~(\ref{lemm characterization neighborhood 1}) we have $U \in {\mc N} \!\left\{ f(z) \right\}$, and therefore $f^{-1} \left[ U \right] \in {\mc N} \!\left\{ z \right\}$ by assumption. Thus $f^{-1} \left[ U \right] \in \tp_X$ by Lemma~\ref{lemm characterization neighborhood}~(\ref{lemm characterization neighborhood 1}).

To show that (\ref{theo local cont 3}) implies (\ref{theo local cont 6}), let $(x_n : n \in D)$ be a net in~$X$ with $x_n \rightarrow x$, and let $U \in {\mc N} \!\left\{ f(x) \right\}$. By assumption we have $f^{-1} \left[ U \right] \in {\mc N} \!\left\{ x \right\}$. There is $n \in D$ such that $x_k \in f^{-1} \left[ U \right]$ for $k \in D$, $k \geq n$. It follows that $f(x_k) \in U$ for $k \in D$, $k \geq n$.

To prove that (\ref{theo local cont 6}) implies (\ref{theo local cont 5}), let $\mc B$ be a filter base on~$X$ with ${\mc B} \rightarrow x$. Further let ${\mc F} = \Phi ({\mc B})$, $(x_n : n \in D)$ the net generated by~$\mc F$, and $U \in {\mc N} \!\left\{ f(x) \right\}$. It follows that $x_n \rightarrow x$, and thus $f(x_n) \rightarrow f(x)$ by assumption. Let $(y_m : m \in E)$ be the net generated by~$f \, \llbracket {\mc F} \, \rrbracket \,$. Then we have $\lim_n f(x_n) \subset \lim_m y_m$\,.
\\

\hspace{0.05\textwidth}
\parbox{0.95\textwidth}
{[Let $r \in \lim_n f(x_n)$ and $U \in {\mc N} \!\left\{ r \right\}$. We may choose $(z,F) = n \in D$ such that $f(x_k) \in U$ for $k \in D$, $k \geq n$. Let $e = \left( f(z), f \left[ F \right] \right)$. We clearly have $e \in E$. Now let $(v,V) = g \in E$ with $g \geq e$. We define $C = F \cap f^{-1} \left[ V \right]$. There is $A \in {\mc F}$ such that $f \left[ A \right] = V$. Since $F \cap A \in {\mc F}$ and $A \subset f^{-1} \left[ V \right]$, we have $C \in {\mc F}$. Further we have $f \left[ C \right] = V$ (exercise). Hence there is $c \in C$ with $f(c) = v$. Thus we have $g = \left( f(c), f \left[ C \right] \right)$. Since $(c,C) \geq (z,F)$, we have $y_g = f(c) = f \left( x_{(c,C)} \right) \in U$.]}
\\

It follows that $f \, \llbracket {\mc F} \, \rrbracket \rightarrow f(x)$, and thus $f \, \llbracket {\mc B} \, \rrbracket \rightarrow f(x)$ by Lemma~\ref{lemm image filter base id}.

To prove that (\ref{theo local cont 5}) implies (\ref{theo local cont 4}), note that ${\mc N} \!\left\{ x \right\}$ is a filter base on~$X$ with \mbox{${\mc N} \!\left\{ x \right\} \rightarrow x$}. Thus we have $f \, \llbracket \, {\mc N} \!\left\{ x \right\} \rrbracket \rightarrow f(x)$ by assumption. Let $U \in {\mc N} \!\left\{ f(x) \right\}$. It follows that there is $V \in {\mc N} \!\left\{ x \right\}$ such that $f \left[ V \right] \subset U$.

To see that (\ref{theo local cont 4}) implies (\ref{theo local cont 3}), let $U \in {\mc N} \!\left\{ f(x) \right\}$. By assumption there is $V \in {\mc N} \!\left\{ x \right\}$ such that $f \left[ V \right] \subset U$. It follows that $V \subset f^{-1} \left[ U \right]$, and thus $f^{-1} \left[ U \right] \in {\mc N} \!\left\{ x \right\}$.
\eproof

\bdefi
\index{Continuous!topological space}
\index{Homeomorphism}
Given two topological spaces $(X,\tp_X)$, $(Y,\tp_Y)$ and a map $f: X \longrightarrow Y$, $f$~is called $\tp_X\,$-$\tp_Y\,$-{\bf continuous}, or short {\bf continuous}, if $f$ is continuous in~$x$ for every $x \in X$. If $f$ is bijective and if both $f$ and $f^{-1}$ are continuous, then $f$ is called a $\tp_X\,$-$\tp_Y\,$-{\bf homeomorphism}, or short {\bf homeomorphism}.
\edefi

\btheo
\label{theo global cont}
Let $(X,\tp_X)$ and $(Y,\tp_Y)$ be two topological spaces, ${\mc C}_X$ and ${\mc C}_Y$ the systems of all closed subsets of $X$ and~$Y$, respectively, ${\mc S}_X$ a subbase for~$\tp_X$, and $f : X \longrightarrow Y$ a map, the following statements are equivalent:
\benum
\item \label{theo global cont 1} $f$ is continuous.
\item \label{theo global cont 2} $f^{-1} \, \llbracket {\mc S}_Y \rrbracket \, \subset \, \tp_X$
\item \label{theo global cont 3} $\forall x \in X \quad f^{-1} \, \llbracket \, {\mc N} \!\left\{ f(x) \right\} \rrbracket \, \subset \, {\mc N} \!\left\{ x \right\}$
\item \label{theo global cont 4} $\forall x \in X \quad {\mc N} \!\left\{ f(x) \right\} \, \subset_{\stackrel{}{\Phi}} f \, \llbracket \, {\mc N} \!\left\{ x \right\} \rrbracket$\,, \, i.e.\ $f$ is ${\mc N} \!\left\{ x \right\}$-${\mc N} \!\left\{ f(x) \right\}$-continuous in~$x$ for every $x \in X$.
\item \label{theo global cont 5} For every $x \in X$ and every filter base $\mc B$ on $X$, ${\mc B} \rightarrow x$ implies $f \, \llbracket {\mc B} \, \rrbracket \rightarrow f(x)$.
\item \label{theo global cont 6} For every $x \in X$ and every net $(x_n)$ in~$X$, $x_n \rightarrow x$ implies $f(x_n) \rightarrow f(x)$.
\item \label{theo global cont 7} $f^{-1} \, \llbracket \tp_Y \rrbracket \, \subset \, \tp_X$
\item \label{theo global cont 8} $f^{-1} \, \llbracket \, {\mc C}_Y \rrbracket \, \subset \, {\mc C}_X$
\eenum

\etheo

\bproof
The equivalence of (\ref{theo global cont 1}) to (\ref{theo global cont 6}) follows by Theorem~\ref{theo local cont}.

The equivalence of (\ref{theo global cont 1}) and (\ref{theo global cont 7}) is obvious.

Finally, the equivalence of (\ref{theo global cont 7}) and (\ref{theo global cont 8}) follows by complementation.
\eproof

Further equivalent statements are proven in Theorem~\ref{theo cont closure int} below.

\midvspace

\brema
Given two topologies $\tp_1$ and $\tp_2$ on a set~$X$, $\tp_1$ is finer than $\tp_2$ iff $\mathrm{id}_X$ is $\tp_1\,$-$\tp_2\,$-continuous.
\erema

In the case of a first countable domain space, continuity may be characterized through the convergence of sequences as follows.

\midvspace

\blemm
\label{lemm continuity sequence}
Let $(X,\tp_X)$ and $(Y,\tp_Y)$ be two topological spaces where $\tp_X$ is first countable, $x \in X$, and $f : X \longrightarrow Y$ a map. Then the following statements are equivalent:

\benum
\item \label{lemm continuity sequence 1} $f$ is continuous in~$x$.
\item \label{lemm continuity sequence 2} For every sequence $(x_n)$ in $X$, $x_n \rightarrow x$ implies $f(x_n) \rightarrow f(x)$.
\eenum

\elemm

\bproof
(\ref{lemm continuity sequence 1}) clearly implies (\ref{lemm continuity sequence 2}).

We prove that (\ref{lemm continuity sequence 2}) implies~(\ref{lemm continuity sequence 1}). By Lemma~\ref{lemm countable int base} we may choose a neighborhood base ${\mc B} = \left\{ B_n : n \in \naturalnumbers \right\}$ of~$x$ such that all $B_n$ ($n \in \naturalnumbers$) are open and $B_n \subset B_m$ for every $m, n \in \naturalnumbers$ with~$m < n$. Assume that (\ref{lemm continuity sequence 2}) is true and that $f$ is not continuous in~$x$. Then there is $U \in {\mc N} \!\left\{ f(x) \right\}$ such that $f^{-1} \left[ U \right] \notin {\mc N} \!\left\{ x \right\}$ by Theorem~\ref{theo local cont}. We may choose a sequence~$(x_n)$ such that $x_n \in B_n \!\setminus\! f^{-1} \left[ U \right]$ for every~$n \in \naturalnumbers$. It follows that $x_n \rightarrow x$, and thus $f(x_n) \rightarrow f(x)$ by assumption. Therefore there is $m \in \naturalnumbers$ such that $f(x_n) \in U$ for $n \in \naturalnumbers$, $n \geq m$. Hence $x_n \in f^{-1} \left[ U \right]$ for $n \in \naturalnumbers$ with $n \geq m$, which is a contradiction.
\eproof

\bdefi
\index{Homeomorphic}
Two topological spaces $(X,\tp_X)$ and $(Y,\tp_Y)$ are called {\bf homeomorphic} if there exists a homeomorphism $f : X \longrightarrow Y$.
\edefi

\blemm
\label{lemm homeo}
Let $(X,\tp_X)$, $(Y,\tp_Y)$ be homeomorphic topological spaces and let $f : X \longrightarrow Y$ be a homeomorphism. Then $f \, \llbracket \tp_X \rrbracket = \tp_Y$ and $f^{-1} \, \llbracket \tp_Y \rrbracket = \tp_X$.
\elemm

\bproof
We have $f^{-1} \, \llbracket \tp_Y \rrbracket \subset \tp_X$ by the continuity of~$f$, and $f \, \llbracket \tp_X \rrbracket \subset \tp_Y$ by the continuity of~$f^{-1}$. We define the functions
\begin{center}
\begin{tabular}{ll}
$F : \tp_X \longrightarrow \tp_Y$\,, & $F(A) = f \left[ A \right]$\;;\\[.8em]
$G : \tp_Y \longrightarrow \tp_X$\,, & $G(B) = f^{-1} \left[ B \right]$
\end{tabular}
\end{center}

It follows that
\[
G \circ F (A) = f^{-1} \left[ f \left[ A \right] \right] = A \,, \quad \quad F \circ G (B) = f \left[ f^{-1} \left[ B \right] \right] = B
\]

for every $A \in \tp_X$ and $B \in \tp_Y$. Therefore $F$ and $G$ are bijective. Thus
\[
f \, \llbracket \tp_X \rrbracket = F \left[ \tp_X \right] = \tp_Y \,, \quad \quad f^{-1} \, \llbracket \tp_Y \rrbracket = G \left[ \tp_Y \right] = \tp_X
\]

\eproof

\bexam
Let $a$, $b$ be two sets, and $X = \left\{ a, b \right\}$. Then the systems
\[
\tp_a = \big\{ \O, X, \left\{ a \right\} \!\big\}\,, \quad \quad \tp_b = \big\{ \O, X, \left\{ b \right\} \!\big\}
\]

are topologies on~$X$. The topological spaces $(X,\tp_a)$ and $(X,\tp_b)$ are homeomorphic, and the function
\[
f : X \longrightarrow X\,, \quad f(a) = b, \quad f(b) = a
\]

is a $\tp_a$\,-$\tp_b$\,-homeomorphism. If $a \neq b$, then $\tp_a \neq \tp_b$\,.
\eexam

Lemma~\ref{lemm homeo} says that two homeomorphic topological spaces are essentially the same, that is they have all properties in common that are related to their topologies. For example if a topological space is first or second countable, also its homeomorphic counterpart is first or second countable, respectively. In some cases a subspace of a topological space whose properties are known and the topological space to be analysed are homeomorphic. Then the space in question automatically has those properties that are inherited by the subspace. The following Definition is suitable for such cases.

\midvspace

\bdefi
\index{Embedding}
Let $(X,\tp_X)$ and $(Y,\tp_Y)$ be two topological spaces and $f : X \longrightarrow Y$ a map. Then $f$ is called {\bf embedding of} $X$ {\bf in}~$Y$ if the map $g : X \longrightarrow f \left[ X \right]$, $g(x) = f(x)$, is a homeomorphism.
\edefi

We now introduce the third type of continuity, both in a point and globally, viz.\ for pseudo-metric spaces.

\midvspace

\bdefi
\index{Continuous}
Given two pseudo-metric spaces $(X,d_X)$, $(Y,d_Y)$, a map $f : X \longrightarrow Y$, and a point $x \in X$, $f$~is called $d_X$-$d_Y$-{\bf continuous in}~$x$, or short {\bf continuous in}~$x$, if
\[
\forall \varepsilon \in \, \left] 0, \infty \right[ \quad \exists \,\delta \in \, \left] 0, \infty \right[ \quad \forall y \in X \quad d_X(x,y) < \delta \; \Longrightarrow \; d_Y \big( f(x), f(y) \big) < \varepsilon
\]

If $f$ is continuous in~$x$ for every $x \in X$, then $f$ is called $d_X$-$d_Y$-{\bf continuous}, or short {\bf continuous}.
\edefi

The following Lemma states that the definitions of continuity for pseudo-metric spaces and for topological spaces are in agreement with the generation of a topology by a pseudo-metric.

\midvspace

\blemm
\label{lemm cont pseudo-m top}
Let $(X,d_X)$ and $(Y,d_Y)$ be pseudo-metric spaces, respectively, $f : X \longrightarrow Y$ a map, and $x \in X$. $f$~is $d_X$-$d_Y$-continuous in~$x$ iff $f$ is $\tau(d_X)$-$\tau(d_Y)$-continuous in~$x$. Further, $f$ is $d_X$-$d_Y$-continuous iff $f$ is $\tau(d_X)$-$\tau(d_Y)$-continuous.
\elemm

\bproof
To prove the first claim, it is enough to show that $f$ is $d_X$-$d_Y$-continuous in $x$ iff $f$ is ${\mc N} \!\left\{ x \right\}$-${\mc N} \!\left\{ f(x) \right\}$-continuous in~$x$ by Theorem~\ref{theo local cont}.

First assume that $f$ is $d_X$-$d_Y$-continuous in~$x$, and let $U \in {\mc N} \!\left\{ f(x) \right\}$. By Lemma and Definition~\ref{lede pseudo-m topology} there is $r \in \left] 0, \infty \right[\,$ such that $B \subset U$ where $B$ is the open sphere about~$f(x)$ with $d_Y$-radius~$r$. It follows that there is $s \in \left] 0, \infty \right[\,$ such that $f \left[ A \right] \subset B$ where $A$ is the open sphere about~$x$ with $d_X$-radius~$s$. Since $A \in {\mc N} \!\left\{ x \right\}$, this shows that $f$ is ${\mc N} \!\left\{ x \right\}$-${\mc N} \!\left\{ f(x) \right\}$-continuous.

Now assume that $f$ is ${\mc N} \!\left\{ x \right\}$-${\mc N} \!\left\{ f(x) \right\}$-continuous. Let $r \in \left] 0, \infty \right[\,$ and $B$ the open sphere about~$f(x)$ with $d_Y$-radius~$r$. Then we have $B \in {\mc N} \!\left\{ f(x) \right\}$. There is \mbox{$C \in {\mc N} \!\left\{ x \right\}$} such that $f \left[ C \right] \subset B$ by assumption. By Lemma and Definition~\ref{lede pseudo-m topology} there is $s \in \left] 0, \infty \right[\,$ such that $A \subset C$ where $A$ is the open sphere about~$x$ with $d_Y$-radius~$r$. It follows that $f$ is $d_X$-$d_Y$-continuous in~$x$.

The second claim follows by the first one.
\eproof

Both with respect to filters and with respect to topologies, the composition of two continuous functions is continuous.

\midvspace

\blemm
\label{lemm cont comp}
Let $(X_i,\tp_i)$ ($i \in \left\{ 1, 2, 3 \right\}$) be topological spaces, $f : X_1 \longrightarrow X_2$ and $g : X_2 \longrightarrow X_3$ maps, $x_1 \in X_1$ a point, $x_2 = f(x_1)$, and $x_3 = g(x_2)$. For every $i \in \left\{ 1, 2, 3 \right\}$ let $\mc F_i$ be a filter on~$X_i$ such that $x_i$ is a cluster point of~$\mc F_i$\,. If $f$ is ${\mc F}_1\,$-${\mc F}_2\,$-continuous in $x_1$ and $g$ is ${\mc F}_2\,$-${\mc F}_3\,$-continuous in~$x_2$\,, then $g \circ f$ is ${\mc F}_1\,$-${\mc F}_3\,$-continuous in~$x_1$\,.
\elemm

\bproof
Exercise.
\eproof

\blemm
\label{composition filter cont}
Let $(X_i,\tp_i)$ ($i \in \left\{ 1, 2, 3 \right\}$) be topological spaces, \mbox{$f : X_1 \longrightarrow X_2$} and \mbox{$g : X_2 \longrightarrow X_3$} maps, and $x \in X_1$\,. If $f$ and $g$ are continuous, then $g \circ f$ is continuous. If $f$ is continuous in~$x$, and $g$ is continuous in~$f(x)$, then $g \circ f$ is continuous in~$x$.
\elemm

\bproof
To see the second claim, assume that $f$ is $\tp_1\,$-$\tp_2\,$-continuous in~$x$ and $g$ is $\tp_2\,$-$\tp_3\,$-continuous in~$f(x)$. Let $A \in \tp_3(g(f(x))$. It follows that $g^{-1} \left[ A \right] \in \tp_2(f(x))$ by the continuity of $g$, and
\[
(g \circ f)^{-1} \left[ A \right] = f^{-1} \left[ g^{-1} \left[ A \right] \right] \in \tp_1
\]

Now the first claim clearly follows.
\eproof

We adapt the following convention.

\midvspace

\bdefi
Let $f \subset X \!\times Y$ be a function. If $X = {\mathbb R}$ or $Y = {\mathbb R}$, then continuity of $f$ refers to the standard topology on~${\mathbb R}$ unless otherwise specified.
\edefi

\bdefi
\label{def distance}
\index{Distance}
Given a pseudo-metric space $(X,d)$, the function $\mathrm{dist}_d$, or short $\mathrm{dist}$, defined by
\begin{center}
\begin{tabular}{l}
$\mathrm{dist}_d : \big({\mc P}(X) \!\setminus\! \left\{ \O \right\} \!\big) \times \big({\mc P}(X) \!\setminus\! \left\{ \O \right\} \!\big) \longrightarrow {\mathbb R}$\\[.8em]
$\mathrm{dist}_d(A,B) = \inf \big\{ d(x,y) \, : \, x \in A,\, y \in B \big\}$
\end{tabular}
\end{center}

is called {\bf distance}. If one of the arguments is a singleton, we also write $\mathrm{dist}_d(A,x)$ and $\mathrm{dist}_d(x,B)$ instead of $\mathrm{dist}_d(A,\left\{ x \right\})$ and $\mathrm{dist}_d(\left\{ x \right\},B)$, respectively, where $x \in X$.
\edefi

Notice that the range of the distance is clearly a subset of~${\mathbb R}_+$.

\midvspace

\blemm
\label{lemm dist cont}
Given a pseudo-metric space $(X,d)$ and a set $A \subset X$ where $A \neq \O$, the function \mbox{$f : X \longrightarrow {\mathbb R}$}, $f(x) = \mathrm{dist}(A,x)$ is continuous.
\elemm

\bproof
Let $x \in X$ and $U \in {\mc N} \!\left\{ f(x) \right\}$. We may choose $\varepsilon \in {\mathbb R}$ such that
\[
\, \left] f(x) - \varepsilon,\, f(x) + \varepsilon \right[\, \subset U
\]

Let $u, v \in X$. Then we have
\[
d(x,v) \leq d(x,u) + d(u,v), \quad \quad d(u,v) \leq d(u,x) + d(x,v)
\]

and hence
\[
f(x) \leq d(x,u) + f(u), \quad \quad f(u) \leq d(x,u) + f(x)
\]

by Lemmas~\ref{lemm inf const} and~\ref{lemm two func inf sup}. It follows that $|f(x) - f(u)| \leq d(x,u)$. Thus $f \left[ B(x,\varepsilon)\right] \subset U$.
\eproof

In the case of interval topologies, continuity can by characterized in the following way.

\midvspace

\brema
\label{defi order cont gen}
Let $X$ and $Y$ be two sets, ${\mc R}$ and ${\mc S}$ systems of pre-orderings on $X$ and $Y$, respectively, $f : X  \longrightarrow Y$ a function, and $x \in X$. The following statements are equivalent by Theorem~\ref{theo local cont}~(\ref{theo local cont 2a}):

\benum
\item $f$ is continuous in~$x$.
\item For every $S \in {\mc S}$ and $y \in Y$ with $(y,f(x)) \in S$ there is a finite index set~$K$ and, for every $k \in K$, a pre-ordering $R_k \in {\mc R}$, a point $x_k \in X$, and an interval with either $I_k = \,\left] -\infty, x_k \right[_{\,R(k)}\,$ or $I_k = \,\left] x_k, \infty \right[_{\,R(k)}\,$, such that $x \in \bigcap_{k \in K} I_k$ and $f\left[ \bigcap_{k \in K} I_k\right] \subset \,\left] y, \infty \right[_{\,S}\,$.

Moreover, for every $S \in {\mc S}$ and $y \in Y$ with $(f(x),y) \in S$ there is a finite index set~$K$ and, for every $k \in K$, a pre-ordering $R_k \in {\mc R}$, a point $x_k \in X$, and an interval with either $I_k = \,\left] -\infty, x_k \right[_{\,R(k)}\,$ or $I_k = \,\left] x_k, \infty \right[_{\,R(k)}\,$, such that $x \in \bigcap_{k \in K} I_k$ and $f\left[ \bigcap_{k \in K} I_k \right] \subset \,\left] -\infty, y \right[_{\,S}\,$.
\eenum

\erema

Finally, we introduce the notions of open and closed maps in this Section, which are, for example, relevant in the proof of Theorem~\ref{theo pseudo-metric metric} where a metric space is generated from a pseudo-metric space.

\midvspace

\bdefi
\index{Open!function}
\index{Function!open}
\index{Closed!function}
\index{Function!closed}
Let $(X,\tp_X)$ and $(Y,\tp_Y)$ be two topological spaces, ${\mc C}_X$ and ${\mc C}_Y$ the systems of all $\tp_X\,$-closed and all $\tp_Y\,$-closed sets, respectively, and $f : X \longrightarrow Y$ a map. If $f \, \llbracket \tp_X \rrbracket \, \subset \, \tp_Y$\,, then $f$~is called $\tp_X\,$-$\tp_Y\,$-{\bf open}, or short {\bf open}. If $f \, \llbracket {\mc C}_X \rrbracket \, \subset \, {\mc C}_Y$\,, then $f$~is called $\tp_X\,$-$\tp_Y\,$-{\bf closed}, or short {\bf closed}.
\edefi

\blemm
\label{lemm charac open}
Let $(X,\tp_X)$, $(Y,\tp_Y)$ be two topological spaces, $\mc B$ a base for $\tp_X$, and $f : X \longrightarrow Y$ a map. Then the following statements are equivalent:

\benum
\item \label{lemm charac opem 1} $f$ is open.
\item \label{lemm charac opem 2} $f \, \llbracket {\mc B} \, \rrbracket \, \subset \, \tp_Y$
\eenum

\elemm

\bproof
Note that (\ref{lemm charac opem 2}) implies (\ref{lemm charac opem 1}) by Lemma~\ref{lemm func relations}~(\ref{lemm func relations 1}).
\eproof

\blemm
\label{isom cont open}
With definitions as in Lemma~\ref{lemm cont pseudo-m top}, if $f$ is an isometry, then $f$ is continuous and open.
\elemm

\bproof
Assume that $f$ is an isometry. Let $x \in X$, $y \in Y$, and $r \in \left] 0, \infty \right[\,$. We have $f \left[ B(x,r) \right] = B(f(x),r)$. Thus $f$ is open by Lemma~\ref{lemm charac open}. Moreover, we have $f^{-1} \left[ B(y,r) \right] = \bigcup \left\{ B(z,r) \, : \, z \in f^{-1} \left\{ y \right\} \right\}$. Hence $f$ is continuous by Theorem~\ref{theo global cont}~(\ref{theo global cont 2}).
\eproof

\section{Closure, interior, derived set, boundary}
\label{closure}

We now investigate, for a given topological space $(X,\tp)$ and an arbitrary subset $A \subset X$, points whose neighborhood system has particular properties with respect to the set~$A$. Specifically, we introduce the notions of interior points and their ensemble, called the interior of the set, accumulation points and their ensemble, the derived set, boundary points and their ensemble, called the boundary of the set, and finally the closure of the set.

\midvspace

\bdefi
\label{defi interior}
\index{Interior point}
\index{Interior}
\index{Set!interior}
Let $\xi = (X,\tp)$ be a topological space and $A \subset X$. A point $x \in X$ is called {\bf interior point of}~$A$ if $A \in {\mc N} \!\left\{ x \right\}$. The set of all interior points of~$A$ is called the {\bf interior of}~$A$ and is denoted by $\mathrm{int}_{\xi}(A)$ or $\mathrm{int}_{\xi} A$. If the set~$X$ is evident from the context, we also write $\mathrm{int}_{\tp}(A)$ or $\mathrm{int}_{\tp} A$. If the topological space is evident from the context, we also write $\mathrm{int} (A)$, $\mathrm{int} \, A$, or~$A^{\circ}$.
\edefi

The following is a characterization of the interior of a set.
 
\midvspace

\blemm
\label{lemm int union}
Let $(X,\tp)$ be a topological space, $A \subset X$, and ${\mc A} = \left\{ U \in \tp \, : \, U \subset A \right\}$. We have
\[
A^{\circ} \, = \, {\textstyle \bigcup} \, {\mc A} \, = \, \sup {\mc A} \, \in \, \tp
\]

where the supremum is with respect to the ordering $\subset$ on~${\mc P}(X)$. $A^{\circ}$~is the largest open set contained in~$A$, i.e.\ it is the maximum of~$\mc A$.
\elemm

\bproof
Exercise.
\eproof

\bdefi
\label{defi accum closure}
\index{Accumulation point}
\index{Set!accumulation point}
\index{Derived set}
\index{Set!derived set}
\index{Closure}
\index{Set!closure}
Let $\xi = (X,\tp)$ be a topological space and $A \subset X$. A point $x \in X$ is called {\bf accumulation point of}~$A$ if $(A \cap U) \!\setminus\! \left\{ x \right\} \neq \O$ for every $U \in {\mc N} \!\left\{ x \right\}$. The set of all accumulation points of~$A$ is called the {\bf derived set of}~$A$ and is denoted by $\mathrm{der}_{\xi}(A)$ or $\mathrm{der}_{\xi} A$. If the set~$X$ is evident from the context, we also write $\mathrm{der}_{\tp}(A)$ or $\mathrm{der}_{\tp} A$. If the topological space is evident from the context, we also write $\mathrm{der}(A)$, $\mathrm{der} \, A$, or~$A^d$.

The union $A \cup A^d$ is called the {\bf closure of}~$A$ and is denoted by $\mathrm{cl}_{\xi}(A)$ or $\mathrm{cl}_{\xi} A$. If the set~$X$ is evident from the context, we also write $\mathrm{cl}_{\tp}(A)$ or $\mathrm{cl}_{\tp} A$. If the topological space is evident from the context, we also write $\mathrm{cl}(A)$, $\mathrm{cl} \, A$, or~$\overline{A}$.
\edefi

The following result is a convenient characterization of the closure of a set.

\midvspace

\blemm
\label{lemm closure intersec}
Let $(X,\tp)$ be a topological space, $\mc C$~the system of all closed sets, $A \subset X$, and ${\mc A} = \left\{ B \in {\mc C} \, : \, A \subset B \right\}$. We have
\[
\overline{A} \, = \, {\textstyle \bigcap} \, {\mc A} \, = \, \inf {\mc A} \, \in \, {\mc C}
\]

where the infimum is with respect to the ordering $\subset$ on~${\mc P}(X)$. $\overline{A}$~is the smallest closed set containing~$A$, i.e.\ it is the minimum of~$\mc A$.
\elemm

\bproof
To show the first equation, assume that $x \notin \bigcap {\mc A}$. Then there is $B \in {\mc C}$ with $A \subset B$ and $x \notin B$. It follows that $A \cap B^c = \O$ and $x \in B^c$. Hence $x \notin A^d$. It is also clear that $x \notin A$. Therefore we have $x \notin \overline{A}$.

Conversely, assume that $x \notin \overline{A}$. Then there is $U \in {\mc N} \! \left\{ x \right\}$ such that $A \cap U = \O$. Hence there is $V \in \tp$ such that $x \in V \subset U$. We have $A \cap V = \O$. It follows that $A \subset V^c$, and hence $x \notin \bigcap {\mc A}$.

The second equation follows by Example~\ref{exam sup inf}.

To see the last claim, note that $\overline{A}$ is a lower bound and also a member of~$\mc A$.
\eproof

Those points that belong to the closure of a set may be characterized in terms of the convergence of nets and filters. The following result is applied frequently in the sequel.

\midvspace

\btheo
\label{theo closure}
Let $(X,\tp)$ be a topological space, $A \subset X$, and $x \in X$. The following statements are equivalent:

\benum
\item \label{theo closure 1} $x \in \overline{A}$
\item \label{theo closure 1a} $\forall \, U \in {\mc N} \!\left\{ x \right\} \quad U \cap A \neq \O$
\item \label{theo closure 2} There is a net $(x_n)$ in $A$ such that $x_n \rightarrow x$.
\item \label{theo closure 3} There is a filter $\mc F$ on~$X$ such that $A \in {\mc F}$ and ${\mc F} \rightarrow x$.
\item \label{theo closure 4} There is a filter base $\mc B$ on~$X$ such that ${\mc B} \subset {\mc P}(A)$ and ${\mc B} \rightarrow x$.
\item \label{theo closure 5} $x \in \mathrm{adh} \; {\mc F}$ where ${\mc F} = \left\{ F \subset X \, : \, A \subset F \right\}$
\eenum

\etheo

\bproof
The equivalence of~(\ref{theo closure 1}) and~(\ref{theo closure 1a}) is a direct consequence of Definition~\ref{defi accum closure}.

We now show the equivalence of~(\ref{theo closure 1}), and (\ref{theo closure 2}) to~(\ref{theo closure 5}).

To show that (\ref{theo closure 1}) implies (\ref{theo closure 2}), let $x \in \overline{A}$. If $x \in A$, then we may choose a constant net $(x_n : n \in D)$ in~$A$ with $x_n = x$ ($n \in D$). If $x \in A^d$, we may choose for each $U \in {\mc N} \!\left\{ x \right\}$ a point $x_U \in A \cap U$. Since $\left({\mc N} \!\left\{ x \right\}, \supset \right)$ is a directed space, $\left( x_U \, : \, U \in {\mc N} \!\left\{ x \right\} \right)$ is a net. Moreover, $x_U \rightarrow x$.

To see that (\ref{theo closure 2}) implies (\ref{theo closure 3}), let $(x_n)$ be a net in~$X$ with values in~$A$ such that $x_n \rightarrow x$. Further let $\mc F$ be the filter on~$X$ generated by~$(x_n)$. Then we clearly have $A \in {\mc F}$ and ${\mc F} \rightarrow x$.

To prove that (\ref{theo closure 3}) implies (\ref{theo closure 4}), let $\mc F$ be a filter on~$X$ such that $A \in {\mc F}$ and ${\mc F} \rightarrow x$. Then ${\mc B} = \left\{ F \cap A \, : \, F \in {\mc F} \right\}$ is a filter base on~$X$. Moreover, we have ${\mc B} \subset {\mc P}(A)$ and ${\mc B} \rightarrow x$.

To show that (\ref{theo closure 4}) implies (\ref{theo closure 5}), let $\mc B$ be a filter base on~$X$ such that ${\mc B} \subset {\mc P}(A)$ and ${\mc B} \rightarrow x$. Further let $U \in {\mc N} \!\left\{ x \right\}$ and $F \in {\mc F}$. There exists $B \in {\mc B}$ such that $B \subset U$. It follows that $B \, \subset \, U \cap A \, \subset \, U \cap F$.

Finally, to see that (\ref{theo closure 5} implies (\ref{theo closure 1}), let $x \in \left( \mathrm{adh} \; {\mc F} \right) \!\setminus\! A$ and $U \in {\mc N} \!\left\{ x \right\}$. Since $A \in {\mc F}$, we have $A \cap U \neq \O$, and thus \mbox{$(A \cap U) \!\setminus\! \left\{ x \right\} \neq \O$}. It follows that $x \in A^d$.
\eproof

The following Theorem provides a characterization of the points of the derived set.

\midvspace

\btheo
\label{theo accum}
Let $(X,\tp)$ be a topological space, $A \subset X$, and $x \in X$. Then the following statements are equivalent:

\benum
\item \label{theo accum 1} $x \in A^d$
\item \label{theo accum 2} There is a net $(x_n)$ in $A \!\setminus\! \left\{ x \right\}$ such that $x_n \rightarrow x$.
\item \label{theo accum 3} There is a filter base $\mc B$ on~$X$ such that ${\mc B} \subset {\mc P}(A \!\setminus\! \left\{ x \right\})$ and ${\mc B} \rightarrow x$.
\eenum

\etheo

\bproof
We first show that (\ref{theo accum 1}) implies (\ref{theo accum 2}). Let $x \in A^d$. It follows that $x \in \left( A \!\setminus\! \left\{ x \right\} \right)^d$.
Since
\[
\left( A \!\setminus\! \left\{ x \right\} \right)^d \, \subset \, \overline{A \!\setminus\! \left\{ x \right\}}
\]

there is a net $(x_n)$ in $A \!\setminus\! \left\{ x \right\}$ such that $x_n \rightarrow x$ by Theorem~\ref{theo closure}.

To see that (\ref{theo accum 2}) implies (\ref{theo accum 3}), let $(x_n)$ be a net in~$X$ with values in $A \!\setminus\! \left\{ x \right\}$ such that $x_n \rightarrow x$. Further let $\mc F$ be the filter on~$X$ generated by~$(x_n)$. Then
\[
{\mc B} \, = \, \big\{ (F \cap A) \!\setminus\! \left\{ x \right\} \, : \, F \in {\mc F} \big\}
\]

is a filter base on~$X$. Moreover, we have ${\mc B} \subset {\mc P}(A \!\setminus\! \left\{ x \right\})$ and ${\mc B} \rightarrow x$.

Finally, to prove that (\ref{theo accum 3}) implies (\ref{theo accum 1}), let $\mc B$ be a filter base on~$X$ such that ${\mc B} \subset {\mc P}(A \!\setminus\! \left\{ x \right\})$ and ${\mc B} \rightarrow x$. Further let $U \in {\mc N} \!\left\{ x \right\}$. There exists $B \in {\mc B}$ such that $B \subset U$. It follows that $\left(A \!\setminus\! \left\{ x \right\}\right) \cap U \neq \O$.
\eproof

In the case of first countability, sequences can be used instead of nets in Theorems~\ref{theo closure} and~\ref{theo accum}.

\midvspace

\blemm
\label{lemm first count sequ}
Let $(X,\tp)$ be a topological space that is first countable, $A \subset X$, and $x \in X$. Then the following statements hold:

\benum
\item \label{lemm first count sequ 2} $x \in \overline{A}$\, iff there is a sequence $(x_n)$ in $A$ such that $x_n \rightarrow x$.
\item \label{lemm first count sequ 1} $x \in A^d$\, iff there is a sequence $(x_n)$ in $A \!\setminus\! \left\{ x \right\}$ such that $x_n \rightarrow x$.
\eenum

\elemm

\bproof
Exercise.
\eproof

For a pseudo-metric space there is a further characterization of the closure of a set.

\midvspace

\blemm
\label{lemma closure dist zero}
Given a pseudo-metric space $(X,d)$ and a subset $A \subset X$, we have \mbox{$\overline{A} = \left\{ x \in X \, : \, \mathrm{dist}(A,x) = 0 \right\}$} where the closure is with respect to the pseudo-metric topology.
\elemm

\bproof
Exercise.
\eproof

\bdefi
\index{Boundary point}
\index{Set!boundary point}
\index{Boundary}
\index{Set!boundary}
Let $\xi = (X,\tp)$ be a topological space and $A \subset X$. A point $x \in X$ is called {\bf boundary point of}~$A$ if $A \cap U \neq \O$ and $(X \setminus\! A) \cap U \neq \O$ for every $U \in {\mc N} \!\left\{ x \right\}$. The set of all boundary points of~$A$ is called the {\bf boundary of}~$A$ and is denoted by $\mathrm{bound}_{\xi}(A)$ or $\mathrm{bound}_{\xi} A$. If the set~$X$ is evident from the context, we also write $\mathrm{bound}_{\tp}(A)$ or $\mathrm{bound}_{\tp} A$. If the topological space is evident from the context, we also write $\mathrm{bound}(A)$ or $\mathrm{bound} \, A$, or~$\partial A$.
\edefi

We have the following characterizations of the boundary.

\midvspace

\blemm
\label{lemm bound}
Let $(X,\tp)$ be a topological space and $A \subset X$. The following statements hold:

\benum
\item \label{lemm bound 1} $x \in \partial A$\, iff there is a net $(x_n)$ in $A$ such that $x_n \rightarrow x$ and a net $(y_m)$ in $A^c$ such that $y_m \rightarrow x$.
\item \label{lemm bound 2} $\partial A = \overline{A} \cap \overline{A^c}$
\eenum

\elemm

\bproof
Exercise.
\eproof

The next Lemma contains various results involving the interior, closure, derived set, and boundary of a set.

\midvspace

\blemm
\label{lemm int closure}
Let $(X,\tp)$ be a topological space, $\mc C$~the system of all closed sets, and $A \subset X$. The following statements hold:

\benum
\item \label{lemm int closure 1} $A^{\circ} \subset A$
\item \label{lemm int closure 3} $\partial A \in {\mc C}$
\item \label{lemm int closure 4} $A^{\circ} \cap \partial A \, = \, \O$
\item \label{lemm int closure 5} $\overline{A} \, = \, A^{\circ} \cup \partial A$
\item \label{lemm int closure 6} $A \in {\mc C} \quad \Longleftrightarrow \quad \overline{A} = A \quad \Longleftrightarrow \quad A^d \subset A \quad \Longleftrightarrow \quad \partial A \subset A$
\item \label{lemm int closure 7} $A \in \tp \quad \Longleftrightarrow \quad A^{\circ} = A$
\item \label{lemm int closure 7a} $\overline{A^c} = (A^{\circ})^c$
\item \label{lemm int closure 8} $(A^{\circ})^{\circ} = A^{\circ}$
\item \label{lemm int closure 9} $\overline{\overline{A}} = \overline{A}$
\item \label{lemm int closure 10} $\partial (\partial A) \subset \partial A$
\eenum
\elemm

\bproof
(\ref{lemm int closure 1}) follows by Definition~\ref{defi interior}.

(\ref{lemm int closure 3}) follows by Lemma~\ref{lemm bound}~(\ref{lemm bound 2}).

(\ref{lemm int closure 4}) is obvious.

In order to see~(\ref{lemm int closure 5}), note that clearly $A^{\circ} \subset \overline{A}$. Moreover, we have $\partial A \subset \overline{A}$ by Lemma~\ref{lemm bound}~(\ref{lemm bound 2}). Conversely, we have $A \subset A^{\circ} \cup \partial A$. Moreover, let $x \in A^d \setminus A$ if such a point exists. Then for every $U \in {\mc N} \!\left\{ x \right\}$ we have $U \cap A \neq \O$ by definition and $U \cap A^c \neq \O$ since $x \in A^c$. Hence $x \in \partial A$.

The first equivalence of~(\ref{lemm int closure 6}) follows by Lemma~\ref{lemm closure intersec}, the second by definition, and the third by~(\ref{lemm int closure 5}).

(\ref{lemm int closure 7}) is obvious.

To show (\ref{lemm int closure 7a}), notice that
\begin{eqnarray*}
(A^{\circ})^c \!\! & = & \!\! \Big( \bigcup \big\{ U \in \tp \, : \, U \subset A \big\} \Big)^c\\[.2em]
 \!\! & = & \!\! \bigcap \big\{ U^c \, : \, U \in \tp,\, U \subset A \big\}\\[.2em]
 \!\! & = & \!\! \bigcap \big\{ C \in {\mc C} \, : \, A^c \subset C \big\} \, = \; \overline{A^c}
\end{eqnarray*}

(\ref{lemm int closure 8}) is a consequence of Lemma~(\ref{lemm int union}) and~(\ref{lemm int closure 7}).

(\ref{lemm int closure 9}) follows by Lemma~(\ref{lemm closure intersec}) and~(\ref{lemm int closure 6}).

(\ref{lemm int closure 10}) follows by (\ref{lemm int closure 3}) and~(\ref{lemm int closure 6}).
\eproof

The interplay of intersection and union with interior, closure, and derived set is now investigated.

\midvspace

\blemm
\label{lemma system closure rules}
Let $(X,\tp)$ be a topological space, $A, B \subset X$, and $A_i \subset X$ ($i \in I$) where $I$ is an index set. The following statements hold:

\benum
\item \label{lemma system closure rules 1} $(A \cap B)^{\circ} \, = \, A^{\circ} \cap B^{\circ}$
\item \label{lemma system closure rules 2} $\big( \bigcap_{i \in I} A_i \big)^{\circ} \, \subset \; \bigcap_{i \in I} A_i^{\circ}$
\item \label{lemma system closure rules 3} $\overline{ \bigcap_{i \in I} A_i } \; \subset \; \bigcap_{i \in I} \overline{A_i}$
\item \label{lemma system closure rules 4} $\big( \bigcap_{i \in I} A_i \big)^d \, \subset \; \bigcap_{i \in I} A_i^d$
\item \label{lemma system closure rules 5} $\bigcup_{i \in I} A_i^{\circ} \; \subset \; \big( \bigcup_{i \in I} A_i \big)^{\circ}$
\item \label{lemma system closure rules 6} $\overline{A \cup B} \, = \, \overline{A} \cup \overline{B}$
\item \label{lemma system closure rules 7} $\bigcup_{i \in I} \overline{A_i} \; \subset \; \overline{\bigcup_{i \in I} A_i}$
\item \label{lemma system closure rules 8} $(A \cup B)^d \, = \, A^d \cup B^d$
\item \label{lemma system closure rules 9} $\bigcup_{i \in I} A_i^d \; \subset \; \big( \bigcup_{i \in I} A_i \big)^d$
\eenum

\elemm

\bproof
(\ref{lemma system closure rules 1}), (\ref{lemma system closure rules 2}), (\ref{lemma system closure rules 5}), and (\ref{lemma system closure rules 9}) hold by definition.

To see (\ref{lemma system closure rules 3}) and (\ref{lemma system closure rules 4}), we define $B = \bigcap {\mc A}$. For every $A \in {\mc A}$ we have $B \subset A$ and therefore $\overline{B} \subset \overline{A}$ by Theorem~\ref{theo closure}, and $B^d \subset A^d$ by Theorem~\ref{theo accum}.

In order to see (\ref{lemma system closure rules 6}), note that
\begin{eqnarray*}
\overline{A \cup B} \!\! & = & \!\! \overline{\left( A^c \cap B^c \right)^c} \, = \, \left(\left( A^c \cap B^c \right)^{\circ}\right)^c \, = \, \left( \left( A^c \right)^{\circ} \cap \left( B^c \right)^{\circ} \right)^c \\[.2em]
 & = & \!\! \left( \left( A^c \right)^{\circ} \right)^c \cup  \left( \left( B^c \right)^{\circ} \right)^c \, = \, \overline{A} \cup \overline{B}
\end{eqnarray*}

by (\ref{lemma system closure rules 1}) and~Lemma~\ref{lemm int closure}~(\ref{lemm int closure 7a}).

(\ref{lemma system closure rules 7}) follows by Theorem~\ref{theo closure}.

To show (\ref{lemma system closure rules 8}) notice that $(A \cup B)^d \supset A^d \cup B^d$ clearly holds. Conversely, let $x \in (A \cup B)^d $. Assume there are $U, V \in {\mc N}^{\, \mathrm{open}}_{\phantom{1}} \!\left\{ x \right\}$ such that $(U \cap A) \!\setminus\! \left\{ x \right\} = \O$ and $(V \cap B) \!\setminus\! \left\{ x \right\} = \O$. Thus we have $\big((U \cap V) \cap (A \cup B)\big) \setminus \left\{ x \right\} = \O$, which is a contradiction.
\eproof

The following Remarks and Examples demonstrate that stricter statements than those in Lemma~\ref{lemma system closure rules} are generally not possible to achieve.

\midvspace

\bexam
In order to show that equality does generally not hold in Lemma~\ref{lemma system closure rules}~(\ref{lemma system closure rules 2}), consider the standard ordering~$<$ on~$\mathbb R$ and the system of proper intervals \mbox{${\mc A} = \left\{ \; ] {-n}^{-1}, n^{-1} [\; \, : \, n \in \naturalnumbers,\, n > 0 \right\}$}.
\eexam

\bexam
Generally equality does not hold in~(\ref{lemma system closure rules 3}) and~(\ref{lemma system closure rules 4}) even if $I$ is finite. To see this consider the standard ordering~$<$ on~$\mathbb R$ and the system ${\mc A} = \left\{ \;\left] -\infty, 0 \right[,\, \left] 0, \infty \right[\; \right\}$.
\eexam

\bexam
Generally equality does not hold in~(\ref{lemma system closure rules 5}) even if $I$ is finite. To see this consider the standard ordering~$<$ on~$\mathbb R$ and the system ${\mc A} = \left\{ \, \left[ 0, 1\right],\, \left[ 1, 2 \right] \, \right\}$.
\eexam

\bexam
To see that equality does generally not hold in (\ref{lemma system closure rules 7}) and~(\ref{lemma system closure rules 9}) consider the standard ordering~$<$ on~$\mathbb R$ and the system ${\mc A} = \left\{ \left[ 0, 1 - n^{-1} \right] \, : \, n \in \naturalnumbers,\, n > 0 \right\}$.
\eexam

A condition under which the reverse of Lemma~\ref{lemma system closure rules}~(\ref{lemma system closure rules 7}) is true is provided in Lemma~\ref{lemma sigma locally finite closure} below.

It is possible to characterize a topology on a given set~$X$ by expressing what it means to form the closure~$\overline{A}$ of~$A$ for every $A \subset X$. This characterization is provided in Theorem~\ref{theo closure top} below. We begin with listing the relevant properties of the function that maps a set~$A$ on its closure with respect to a given topology.

\midvspace

\bdefi
\label{defi closure operator}
\index{Closure operator}
Given a set $X$, a function $f : {\mc P}(X) \longrightarrow {\mc P}(X)$ is called {\bf closure operator on}~$X$ if it has the following properties:

\benum
\item \label{defi closure operator 1} $f(\O) = \O$
\item \label{defi closure operator 2} $A \subset f(A)$
\item \label{defi closure operator 3} $f(A \cup B) = f(A) \cup f(B)$
\item \label{defi closure operator 4} $f ( f(A)) = f(A)$
\eenum

\edefi

The properties (\ref{defi closure operator 1}) to (\ref{defi closure operator 4}) in Definition~\ref{defi closure operator} are called Kuratowski axioms in the literature, cf.~\cite{Gaal}.

\midvspace

\btheo
\label{theo closure top}
Let $X$ be a set, $f$ a closure operator on~$X$, and ${\mc C} = \left\{ A \subset X \, : \, f(A) = A \right\}$. The system $\tp = \left\{ A^c \, : \, A \in {\mc C} \right\}$ is the unique topology on~$X$ such that $\mathrm{cl}_{\tp} A = f(A)$ for every $A \subset X$. Moreover, we have $\mathrm{ran} \, f = {\mc C}$.
\etheo

\bproof
$\mc C$~has properties (\ref{systemclosed1}) to (\ref{systemclosed3}) in Lemma~\ref{systemclosed}.
\\

\hspace{0.05\textwidth}
\parbox{0.95\textwidth}
{[We have $\O \in {\mc C}$ because of property~(\ref{defi closure operator 1}) in Definition~\ref{defi closure operator}, and $X \in {\mc C}$ because of property~(\ref{defi closure operator 2}). Further, property~(\ref{defi closure operator 3}) implies that $A \cup B \in {\mc C}$ for every $A, B \in {\mc C}$. Finally let ${\mc A} \subset {\mc C}$ where ${\mc A} \neq \O$. We clearly have $\bigcap {\mc A} \subset f \big( \bigcap {\mc A} \big)$ by property~(\ref{defi closure operator 2}). Moreover, for every $A \in {\mc A}$ we have $\bigcap {\mc A} \subset A$, and therefore $f \big( \bigcap {\mc A} \big) \subset f(A) = A$, and thus $f \big( \bigcap {\mc A} \big) \subset \bigcap {\mc A}$. Hence we obtain $f \big( \bigcap {\mc A} \big) = \bigcap {\mc A}$, and therefore $\bigcap {\mc A} \in {\mc C}$.]}
\\

Therefore $\tp$ is a topology on~$X$ by Lemma~\ref{lemm closed sets}, and $\mc C$ is the system of all $\tp$-closed sets. 

Now let $A \subset X$. We have $f(A) \in {\mc C}$ by Definition~\ref{defi closure operator}~(\ref{defi closure operator 4}). It follows that
\[
\mathrm{cl}_{\tp} A \, = \, \bigcap \big\{ B \in {\mc C} \, : \, A \subset B \big\} \, \subset \, f(A)
\]

by Lemma~\ref{lemm closure intersec} and Definition~\ref{defi closure operator}~(\ref{defi closure operator 2}). Conversely, for every $B \in {\mc C}$, $A \subset B$ implies $f(A) \subset f(B) = B$. Hence we have $f(A) \subset \mathrm{cl}_{\tp} A$. Thus we obtain $\mathrm{cl}_{\tp} A = f(A)$.

To see the uniqueness, note that for a topology on~$X$ with the stated property, the system of all closed sets is precisely ${\mc C}$ by Lemma~\ref{lemm int closure}~(\ref{lemm int closure 6}).

The last claim clearly holds.
\eproof

Similarly to the closure operator we introduce the notion of interior operator.

\midvspace

\bdefi
\label{defi interior operator}
\index{Interior operator}
Given a set $X$, a function $f : {\mc P}(X) \longrightarrow {\mc P}(X)$ is called {\bf interior operator on}~$X$ if it has the following properties:

\benum
\item \label{defi interior operator 1} $f(X) = X$
\item \label{defi interior operator 2} $f(A) \subset A$
\item \label{defi interior operator 4} $f(A \cap B) = f(A) \cap f(B)$
\item \label{defi interior operator 3} $f(f(A)) = f(A)$
\eenum

\edefi

There is a duality between closure and interior operators as follows.

\midvspace

\blemm
\label{lemm interior closure operator}
Let $X$ be a set, $f : {\mc P}(X) \longrightarrow {\mc P}(X)$ a map, and $g : {\mc P}(X) \longrightarrow {\mc P}(X)$ the map defined by $g(A) = \big(f(A^c)\big)^c$. We have $f(A) = \big(g(A^c)\big)^c$ for every $A \subset X$. Moreover, $f$ is a closure operator on~$X$ iff $g$ is an interior operator on~$X$.
\elemm

\bproof
Exercise.
\eproof

It follows that a given interior operator on~$X$ defines a certain topology on~$X$.

\midvspace

\btheo
Let $X$ be a set, $f : {\mc P}(X) \longrightarrow {\mc P}(X)$ an interior operator on~$X$, and $\tp = \left\{ U \subset X \, : \, f(U) = U \right\}$. The system $\tp$ is the unique topology on~$X$ such that $f(A) = \mathrm{int}_{\tp} A$ for every $A \subset X$. Moreover, we have $\mathrm{ran} \, f = \tp$.
\etheo

\bproof
Let $g$ be the closure operator on~$X$ defined by $g(A) = \big(f(A^c)\big)^c$ for every $A \subset X$ (cf.\ Lemma~\ref{lemm interior closure operator}), and $\tp_g = 
\left\{ A^c \, : \, A \subset X,\, g(A) = A \right\}$. Then $\tp_g$ is the unique topology on~$X$ such that $\mathrm{cl}_{\tp(g)} A = g(A)$ for every $A \subset X$ by Theorem~\ref{theo closure top}. We have
\[
\tp_g \, = \, \big\{ U \subset X \, : \, g(U^c) = U^c \big\} \, = \, \tp
\]

which shows that $\tp$ is a topology on~$X$.

Moreover, for every $A \subset X$ we have
\[
f(A) \, = \, \big(g(A^c)\big)^c = \, \big( \mathrm{cl}_{\tp} \left( A^c \right) \!\big)^c = \, \mathrm{int}_{\tp} A
\]

by Lemma~\ref{lemm int closure}~(\ref{lemm int closure 7a}). Since the range of~$g$ is the system of all closed sets, the range of~$f$ is~$\tp$.

To see the uniqueness, notice that a topology on~$X$ with the stated property is equal to~$\tp$ by Lemma~\ref{lemm int closure}~(\ref{lemm int closure 7}).
\eproof

Using closure and interior we obtain further characterizations of continuity of a function.

\midvspace

\btheo
\label{theo cont closure int}
Given two topological spaces $(X,\tp_X)$, $(Y,\tp_Y)$, and a map $f : X \longrightarrow Y$, the following statements are equivalent:

\benum
\item \label{theo cont closure int 1} $f$ is continuous.
\item \label{theo cont closure int 2} $\forall A \subset X \quad f \left[ \overline{A} \right] \, \subset \, \overline{f \left[ A \right]}$
\item \label{theo cont closure int 3} $\forall B \subset Y \quad \overline{f^{-1} \left[ B \right]} \, \subset \, f^{-1}\left[ \overline{B} \right]$
\item \label{theo cont closure int 4} $\forall B \subset Y \quad f^{-1} \left[ B^{\circ} \right] \, \subset \, \left( f^{-1} \left[ B \right] \right)^{\circ}$
\eenum

\etheo

\bproof
The fact that (\ref{theo cont closure int 1}) implies (\ref{theo cont closure int 2}) follows by Theorems~\ref{theo local cont} and~\ref{theo closure}.

To prove that (\ref{theo cont closure int 2}) implies (\ref{theo cont closure int 3}), let $B \subset Y$ and $A = f^{-1} \left[ B \right]$. Then we have $f \left[ A \right] \subset B$, and thus $f \left[ \overline{A} \, \right] \subset \overline{f \left[ A \right]} \subset \overline{B}$. It follows that
\[
\overline{f^{-1} \left[ B \right]} \, = \, \overline{A} \, \subset \, f^{-1} \!\left[ f \left[ \overline{A} \, \right] \right] \, \subset \, f^{-1} \left[ \overline{B} \, \right]
\]

To see that (\ref{theo cont closure int 3}) implies (\ref{theo cont closure int 4}), let $B \subset Y$. We have
\begin{center}
\begin{tabular}{l}
$f^{-1} \left[ B^{\circ} \right] \, = \, f^{-1} \left[ \left( \overline{B^c} \, \right)^c \right] \, = \, \left( f^{-1} \left[ \overline{B^c} \, \right] \right)^c \, \subset \, \left(\overline{f^{-1} \left[ B^c \right]} \, \right)^c$\\[.9em]
\quad \quad \quad $ = \, \left(\overline{\left(f^{-1} \left[ B \right]\right)^c} \, \right)^c \, = \, \left(f^{-1} \left[ B \right] \right)^{\circ}$
\end{tabular}
\end{center}

To see that (\ref{theo cont closure int 4}) implies (\ref{theo cont closure int 1}), let $x \in X$ and $U \in {\mc N} \!\left\{ f(x) \right\}$. Then we have \mbox{$f(x) \in U^{\circ}$}. It follows that
\[
x \, \in \, f^{-1} \left[ U^{\circ} \right] \, \subset \, \left(f^{-1} \left[ U \right]\right)^{\circ} \, \subset \, f^{-1} \left[ U \right]
\]

and thus \mbox{$f^{-1} \left[ U \right] \in {\mc N} \!\left\{ x \right\}$}.
\eproof

In Section~\ref{seqconv} we have investigated the convergence of sequences with respect to comparable topologies on the same set, see Lemmas~\ref{lemm finer seq stronger} and~\ref{lemm seq stronger finer}. In particular, we have found that if a topology $\tp_1$ is finer than a topology $\tp_2$, then $\tp_1$ is also sequentially stronger than~$\tp_2$. As shown now the converse implication holds if the topologies are countable. We also establish similar results for the convergence properties of filters and nets. In the proofs we follow~\cite{Wilansky}.

\midvspace

\bprop
\label{prop equiv top fin closures}
Let $X$ be a set and $\tp_1$, $\tp_2$ two topologies on~$X$. $\tp_1$~is finer than $\tp_2$ iff \mbox{$\mathrm{cl}_{\tp(1)} A \, \subset \, \mathrm{cl}_{\tp(2)} A$}.
\eprop

\bproof
For $i \in \left\{ 1, 2 \right\}$ let ${\mc C}_i$ be the system of all $\tp_i$-closed sets.

First assume that $\tp_2 \subset \tp_1$. It follows that ${\mc C}_2 \subset {\mc C}_1$ and therefore
\begin{eqnarray*}
\mathrm{cl}_{\tp(1)} A \!\! & = & \! \bigcap \big\{ B \in {\mc C}_1 \, : \, A \subset B \big\}\\[.2em]
 & \subset & \! \bigcap \big\{ B \in {\mc C}_2 \, : \, A \subset B \big\} \; = \; \mathrm{cl}_{\tp(2)} A
\end{eqnarray*}

by Lemma~\ref{lemm closure intersec}.

Now assume that $\mathrm{cl}_{\tp(1)} A \subset \mathrm{cl}_{\tp(2)} A$ for every $A \subset X$. Let $A \in \tp_2$. Then we have $A^c \in {\mc C }_2$. Moreover, by assumption we have $\mathrm{cl}_{\tp(1)} A^c \subset \mathrm{cl}_{\tp(2)} A^c = A^c$. It follows that $\mathrm{cl}_{\tp(1)} A^c = A^c$. Thus we have $A^c \in {\mc C}_1$, and hence $A \in \tp_1$.
\eproof

\blemm
\label{lemm seq stronger finer}
Let $X$ be a set, $\tp_1$ and $\tp_2$ two topologies on $X$ where $\tp_1$ is first countable. $\tp_1$~is finer than $\tp_2$ iff $\tp_1$ is sequentially stronger than~$\tp_2$.
\elemm

\bproof
First assume that $\tp_1$~is sequentially stronger than~$\tp_2$. Then we have $\mathrm{cl}_{\tp(1)} A \subset \mathrm{cl}_{\tp(2)} A$ by Lemma~\ref{lemm first count sequ}~(\ref{lemm first count sequ 1}). It follows that $\tp_1$ is finer than $\tp_2$ by Proposition~\ref{prop equiv top fin closures}.

The converse follows by Lemma~\ref{lemm finer seq stronger}.
\eproof

\btheo
\label{theo equiv conver}
Let $X$ be a set, and $\tp_1$ and $\tp_2$ two topologies on~$X$. Then the following statements are equivalent:

\benum
\item \label{theo equiv conver 1} $\tp_1$ is finer than $\tp_2$.
\item \label{theo equiv conver 2} For every net $(x_n)$ in $X$ and every $x \in X$, $x_n \rightarrow x$ with respect to $\tp_1$ implies $x_n \rightarrow x$ with respect to $\tp_2$.
\item \label{theo equiv conver 3} For every filter $\mc F$ in $X$ and every $x \in X$, ${\mc F} \rightarrow x$ with respect to $\tp_1$ implies ${\mc F} \rightarrow x$ with respect to $\tp_2$.
\eenum

\etheo

\bproof
(\ref{theo equiv conver 1}) implies~(\ref{theo equiv conver 2}) by Lemma~\ref{base finer neighborhood}.

To see that (\ref{theo equiv conver 2}) implies~(\ref{theo equiv conver 1}), assume that (\ref{theo equiv conver 2}) holds. Let $A \subset X$. We have $\mathrm{cl}_{\tp(1)} A \subset \mathrm{cl}_{\tp(2)} A$.
\\

\hspace{0.05\textwidth}
\parbox{0.95\textwidth}
{[Let $x \in \mathrm{cl}_{\tp(1)} A$. We may choose a net $(x_n)$ in~$A$ such that $x_n \rightarrow x$ with respect to~$\tp_1$ by Theorem~\ref{theo closure}. Then we have $x_n \rightarrow x$ with respect to~$\tp_2$ by assumption. It follows that $x \in \mathrm{cl}_{\tp(2)} A$.]}
\\

Now (\ref{theo equiv conver 1}) follows by Proposition~\ref{prop equiv top fin closures}.

The equivalence of~(\ref{theo equiv conver 2}) and~(\ref{theo equiv conver 3}) follows by Lemma~\ref{lemm filter gen conv} and by Lemma and Definition~\ref{lede net gen}.
\eproof

\section{Separability}
\label{separability}

In this Section we provide a brief discussion of separability.

\midvspace

\bdefi
\index{Dense}
\index{Set!dense}
Given a topological space $(X,\tp)$ and $A \subset X$, $A$~is called {\bf dense in}~$X$ if $\overline{A} = X$.
\edefi

Order dense sets as introduced in Definition~\ref{defi order dense} in the context of ordered spaces are related to dense sets as follows.

\midvspace

\blemm
\label{lemm dense interval}
Let $(X,\prec)$ be a pre-ordered space where $\prec$ has full field and the interval intersection property. Further let $A \subset X$ be an order dense subset. Then $A$ is dense in~$X$ with respect to the interval topology.
\elemm

\bproof
We define
\begin{eqnarray*}
{\mc S} \!\! & = & \!\! \big\{ \left] -\infty, x \right[\, , \;\left] x, \infty \right[\; : \, x \in X \big\},\\[.2em]
{\mc A} \!\! & = & \!\! \big\{ \,] x , y [\;\, : \, x, y \in X,\; x \prec y \big\},\\[.2em]
{\mc B} \!\! & = &  \!\! {\mc S} \cup {\mc A} \cup \left\{ \O \right\}
\end{eqnarray*}

The system $\mc B$ is a base for the interval topology by Lemma~\ref{lemm interval top base}. Let $x \in X$ and $U \in {\mc N} \!\left\{ x \right\}$. Then there is $B \in {\mc B}$ such that $x \in B \subset U$. Moreover we have $B \cap A \neq \O$ since $A$ is order dense. Thus $x \in \overline{A}$ by Theorem~\ref{theo closure}.
\eproof

\bdefi
\index{Separable}
\index{Topological space!separable}
A topological space is called {\bf separable} if it contains a countable subset that is dense in~$X$.
\edefi

\bexam
Let $<$ be the ordering on~${\mathbb R}_+$ as defined in Lemma and Definition~\ref{lede pos real num}. It has full field and the interval intersection property, cf.\ Remark~\ref{rema r+ interval rel}. Moreover ${\mathbb D}_+$ is order dense in~${\mathbb R}_+$ by Lemma~\ref{lemm D0 dense}. Thus ${\mathbb D}_+$ is dense in~${\mathbb R}_+$ with respect to the standard topology by Lemma~\ref{lemm dense interval}. Since ${\mathbb D}_+$ is countable by Corollary~\ref{coro dyadic countable}, $({\mathbb R}_+,<)$ is separable. Similarly it can be shown that $({\mathbb R},<)$ is separable where $<$ is the standard ordering.
\eexam

We now briefly examine how the properties of first and second countability of a topological space are related to separability.

\midvspace

\blemm
\label{lemma sec count separable}
A topological space $(X,\tp)$ that is second countable is separable.
\elemm

\bproof
Let $\mc B$ be a countable base for~$\tp$. For each $B \in {\mc B}$ with $B \neq \O$ we may choose a point $x_B \in B$. We define $A = \left\{ x_B \, : \, B \in {\mc B} \right\}$. Let $x \in X$ and $U \in {\mc N} \!\left\{ x \right\}$. Then $U \cap A \neq \O$. It follows that $x \in \overline{A}$ by Theorem~\ref{theo closure}.
\eproof

Notice that a separable space need not even be first countable as shown in the following Example that we take from~\cite{Seebach}.

\midvspace

\bexam
Let $X$ be an infinite set that is not countable and $\tp_{\mathrm{cf}}$ the cofinite topology, cf.\ Lemma and Definition~\ref{lede examples topologies}~(\ref{coftop}). 

Let $A \subset X$ be countable and infinite. Then $U \cap A \neq \O$ for every $U \in {\mc N} \! \left\{ x \right\}$ and every $x \in X$. Hence $(X,\tp_{\mathrm{cf}})$ is separable by Theorem~\ref{theo closure}.

Now assume that $(X,\tp_{\mathrm{cf}})$ is first countable. Then there is a point $x \in X$ and a countable neighborhood base~$\mc B$ of~$x$. We define $D = \bigcap {\mc B}$. We have $D = \left\{ x \right\}$ by the definition of~$\tp_{\mathrm{cf}}$, and therefore $X \!\setminus\! \left\{ x \right\} = \bigcup \left\{ B^c \, : \, B \in {\mc B} \right\}$. Since $B^c$ is finite for every $B \in {\mc B}$, it follows that 
$X \!\setminus\! \left\{ x \right\}$ is countable by Lemma~\ref{lemm count times count}, which is a contradiction.
\eexam

Finally the following results holds for pseudo-metrizable spaces.

\midvspace

\btheo
A pseudo-metrizable topological space is separable iff it is second countable.
\etheo

\bproof
Let $\xi = (X,\tp)$ be a pseudo-metrizable topological space and $d$ a pseudo-metric on~$X$ that generates~$\tp$. Assume that $\xi$ is separable. Let $A$ be a countable set that is dense in~$X$ and $R = {\mathbb D}_+ \!\!\setminus\! \left\{ 0 \right\}$. Then the system
\[
{\mc B} = \big\{ B(x,r) \, : \, x \in A,\, r \in R \big\}
\]

is a countable base for~$\tp$ by Lemma~\ref{lemm countable sets}.\\

\hspace{0.05\textwidth}
\parbox{0.95\textwidth}
{[Let $U \in \tp$ and $x \in U$. We may choose $y \in X$ and $r \in R$ such that $x \in B(y,r) \subset U$ by the definition of the pseudo-metric topology. We define $\varepsilon = r - d(x,y)$. $\tp$~is first countable by Theorem~\ref{theo pseudo-m first countable}. Hence there is a sequence $(x_n)$ in~$A$ such that $x_n \rightarrow x$ by Lemma~\ref{lemm first count sequ}~(\ref{lemm first count sequ 2}). We may choose $m \in \naturalnumbers$ such that $d(x_m,x) < \varepsilon / 2$. It follows that $x \in B$ where $B = B(x_m, \varepsilon / 2)$. For every $z \in B$ we have
\[
d(y,z) \, \leq \, d(y,x) + d(x,x_m) + d(x_m,z) \, < \, d(x,y) + \varepsilon \, = \, r
\]

and hence $B \subset B(y,r)$.]}
\\

The converse is true by Lemma~\ref{lemma sec count separable}.
\eproof


\chapter{Generated topologies}
\label{generated topologies}
\setcounter{equation}{0}

\pagebreak

In this Chapter we explore how a topology on a given set can be defined by means of already specified topologies on---generally different---sets and functions between the old and the new sets. Obviously there are two possibilities: the new topology may be generated on the domain of the functions by given topologies on their ranges (so-called inverse image topology, see Section~\ref{inverse image topology}) or generated on their range by given topologies on their domains (so-called direct image topology, see Section~\ref{direct image topology}).  In Sections~\ref{topological subspace} and~\ref{product topology} we analyse two particular cases of inverse image topologies in detail: subspace topologies and product topologies. Section~\ref{quotient topology} is devoted to quotient topologies which are special cases of direct image topologies.

\midvspace

\section{Inverse image topology}
\label{inverse image topology}

\blede
\label{definitialtop}
\index{Inverse image topology}
\index{Topology!inverse image}
\index{Topology!generated}
Given a set $X$, topological spaces $(Y_i, \tp_i)$ $(i \in I)$, where $I$ is an index set, and functions $f_i : X \longrightarrow Y_i$ $(i \in I)$, the system ${\mc S} = \bigcup_{i \in I} f_i^{-1} \, \llbracket \tp_i \, \rrbracket \,$ is a topological subbase on~$X$. The topology on $X$ that is generated by $\mc S$ is called {\bf inverse image topology} or the {\bf topology generated by} $F = \left\{ (f_i, \tp_i) \, : \, i \in I \right\}$ and denoted by~$\tau (F)$. It is the coarsest topology $\tp$ on~$X$ such that $f_i$ is $\tp$-$\tp_i \,$-continuous for every $i \in I$. We also use the notation $x_i = f_i(x)$ for $x \in X$ and $i \in I$.
\elede

\bproof
${\mc S}$ is a topological subbase on~$X$ by Lemma~\ref{lemm crit subbase}. We denote by $\mathscr{A}$ the set of all topologies $\tp$ on~$X$ such that $f_i$ is $\tp$-$\tp_i \,$-continuous for every $i \in I$. Then we have $\tau(F) \in \mathscr{A}$ by Theorem~\ref{theo global cont}. Moreover, for every $\tp \in \mathscr{A}$ we have ${\mc S} \subset \tp$. Hence $\tau(F)$~is the coarsest member of~$\mathscr{A}$ by Lemma~\ref{lemm crit subbase}.
\eproof

\blemm
\label{lemm initial top}
With definitions as in Lemma and Definition~\ref{definitialtop}, let ${\mc S}_i$ be a subbase for~$\tp_i$ for every $i \in I$. The system ${\mc R} = \bigcup_{i \in I} f_i^{-1} \, \llbracket {\mc S}_i \rrbracket \,$ is a subbase for~$\tau (F)$.
\elemm

\bproof
First note that ${\mc R}$ is a topological subbase on~$X$ by Lemma~\ref{lemm crit subbase}. We denote by $\tp$ the topology generated by~${\mc R}$. Now, ${\mc R} \subset {\mc S}$ implies $\tp = \Theta \Psi ({\mc R}) \subset \Theta \Psi ({\mc S}) = \tau (F)$. Moreover, for every $i \in I$, $f_i$ is $\tp$-$\tp_i \,$-continuous by Theorem~\ref{theo global cont}. Since $\tau (F)$ is the coarsest such topology, we have $\tau (F) \subset \tp$.
\eproof

\bcoro
With definitions as in Lemma and Definition~\ref{definitialtop}, let $\mc C$ be the system of all $\tau(F)$-closed sets and, for each $i \in I$, let $\mc E_i$~be a subbase for the $\tp_i$\,-\,closed sets. Then $\bigcup_{i \in I} f_i^{-1} \, \llbracket {\mc E}_i \rrbracket \,$ is a subbase for~$\mc C$.
\ecoro

\bproof
This is a consequence of Lemma~\ref{lemm initial top} by complementation.
\eproof

The following is an important special case of Lemma and Definition~\ref{definitialtop} and Lemma~\ref{lemm initial top}.

\midvspace

\bcoro
\label{coro gen top single space}
Let $X$ be a set and $I$ an index set. For each $i \in I$, let $\tp_i$ be a topology on~$X$ and ${\mc S}_i$ a subbase for~$\tp_i$. Further let $F = \left\{(\mathrm{id}_X,\tp_i) \, : \, i \in I \right\}$. $\bigcup_{i \in I} {\mc S}_i$~is a subbase for~$\tau (F)$. $\tau (F)$~is the supremum of $\left\{ \tp_i \, : \, i \in I \right\}$ in the ordered space $(\mathscr{T}(X),\subset)$, i.e.\ it is the coarsest topology on~$X$ that is finer than $\tp_i$ for every $i \in I$.
\ecoro

\bproof
Exercise.
\eproof

Another relevant special case is that of a topology generated by a single map. We have the following result.

\midvspace

\bcoro
\label{coro gen top single map}
Let $X$ be a set, $(Y,\tp)$ a topological space, $\mc B$~a base for~$\tp$, $\mc S$~a subbase for~$\tp$, $\mc C$~the system of all $\tp$-closed sets, $\mc D$~a base for~$\mc C$, $\mc E$~a subbase for~$\mc C$, $f : X \longrightarrow Y$ a map, and $\tp_X = \tau \, \big(\!\left\{ (f,\tp) \right\}\!\big)$. Then the following statements hold:

\benum
\item $\tp_X = f^{-1} \, \llbracket \tp \, \rrbracket \,$
\item $f^{-1} \, \llbracket {\mc B} \, \rrbracket \,$ is a base for~$\tp_X$.
\item $f^{-1} \, \llbracket {\mc S} \, \rrbracket \,$ is a subbase for~$\tp_X$.
\item $f^{-1} \, \llbracket {\mc C} \, \rrbracket \,$ is the system of all $\tp_X$-closed sets.
\item $f^{-1} \, \llbracket {\mc D} \, \rrbracket \,$ is a base for $f^{-1} \, \llbracket {\mc C} \, \rrbracket \,$.
\item $f^{-1} \, \llbracket {\mc E} \, \rrbracket \,$ is a subbase for $f^{-1} \, \llbracket {\mc C} \, \rrbracket \,$.
\eenum

\ecoro

\bproof
Exercise.
\eproof

\bexam
Let $m, n \in \naturalnumbers$ with $0 < m < n$, $\tp_m$~and $\tp_n$ the standard topologies on~${\mathbb R}^m$ and~${\mathbb R}^n$, respectively, and $x \in {\mathbb R}^{n - m}$. Further we define the function
\[
f : {\mathbb R}^m \longrightarrow {\mathbb R}^n\,, \quad f(z) = (z,x)
\]

where we have identified ${\mathbb R}^n$ and ${\mathbb R}^m \!\times {\mathbb R}^{n - m}$. Then we have~$\tau (\left\{(f,\tp_n)\right\}) = \tp_m$\,.
\eexam

The following Lemma shows how the generation of topologies behaves under the composition of functions.

\midvspace

\blemm
\label{induced assoc}
Let $X$ and $Y_i$ $(i \in I)$ be sets where $I$ is an index set, $(Z_j, \tp_j)$ $(j \in J_i, \, i \in I)$ topological spaces where $J_i$ ($i \in I$) are distinct index sets, and $f_i : X \longrightarrow Y_i$ $(i \in I)$ and $g_j : Y_i \longrightarrow Z_j$ $(i \in I,\, j \in J_i)$ maps. For every $i \in I$, we define the following topology on~$Y_i$\,:
\[
\tp_i \, = \, \tau \, \big(\!\left\{ (g_j, \tp_j) \, : \, j \in J_i \right\}\!\big)
\]

We have
\[
\tau \, \big(\!\left\{ (f_i, \tp_i) \, : \, i \in I \right\}\!\big) \, = \, \tau \,\big(\!\left\{ (g_j \circ f_i, \tp_j) \, : \, i \in I,\, j \in J_i \right\}\!\big)
\]

\elemm

\bproof
This is a consequence of Lemma~\ref{lemm initial top}.
\eproof

\btheo
\label{universalproperty}
Let $(X,\tp)$ be a topological space, $(Y_i, \tp_i)$ $(i \in I)$ topological spaces where $I$ is an index set, $f_i : X \longrightarrow Y_i$ $(i \in I)$ functions, and $F = \left\{(f_i, \tp_i) \, : \, i \in I \right\}$. The following statements are equivalent:

\benum
\item \label{universalproperty 1} $\tp = \tau(F)$
\item \label{universalproperty 2} For every topological space $(Z,\tp_Z)$ and every function $g : Z \longrightarrow X$, $g$~is $\tp_Z\,$-$\tp$-continuous iff $f_i \circ g$ is $\tp_Z\,$-$\tp_i\,$-continuous for every $i \in I$.
\eenum

\etheo

\bproof
To see that (\ref{universalproperty 1}) implies (\ref{universalproperty 2}), assume that $\tau (F) = \tp$. Then $f_i$ is $\tp$-$\tp_i\,$-continuous for every $i \in I$. Further let $(Z,\tp_Z)$ be a topological space and \mbox{$g : Z \longrightarrow X$} a map. Then the continuity of~$g$ implies the continuity of $f_i \circ g$ for every $i \in I$ by Lemma~\ref{composition filter cont}. Conversely, for every $i \in I$ let ${\mc S}_i$ be a subbase for~$\tp_i$\,. Then we have $g^{-1} \left[ f_i^{-1} \left[ S \right] \right] \in \tp_Z$ for every $i \in I$ and $S \in {\mc S}_i$ by the continuity of $f_i \circ g$. Since $\left\{ f_i^{-1} \left[ S \right] \, : \, i \in I,\, S \in {\mc S}_i\right\}$ is a subbase of~$\tp$, the continuity of $g$ follows.

To show that (\ref{universalproperty 2}) implies (\ref{universalproperty 1}), it is enough to show that the topology~$\tp$ is uniquely specified by property~(\ref{universalproperty 2}). Assume that $\tp_1$ and $\tp_2$ are two topologies on~$X$ such that (\ref{universalproperty 2}) is satisfied in both cases. Now let $Z = X$ and $g = \mathrm{id}_X$. Since $g$ is $\tp_m\,$-$\tp_m\,$-continuous for $m \in \left\{ 1, 2 \right\}$, it follows that $f_i$ is $\tp_m\,$-$\tp_i\,$-continuous for $m \in \left\{ 1, 2 \right\}$ and $i \in I$. Thus $g$ is $\tp_1\,$-$\tp_2\,$-continuous and $\tp_2\,$-$\tp_1\,$-continuous, and hence $\tp_1 = \tp_2$\,.
\eproof

Property (\ref{universalproperty 2}) in Theorem~\ref{universalproperty} is often called a universal property. Below we encounter other universal properties a topological or a uniform space may have.

\midvspace

\btheo
\label{theo conv induced}
With definitions as in Lemma and Definition~\ref{definitialtop}, let $(x_n : n \in D)$ be a net in~$X$, $\mc F$ a filter on~$X$, and $x \in X$. The following statements hold:

\benum
\item \label{theo conv induced 1} $x_n \rightarrow x \;\; \Longleftrightarrow \;\; \forall i \in I \quad f_i(x_n) \rightarrow f_i(x)$
\item \label{theo conv induced 2} ${\mc F} \rightarrow x \;\; \Longleftrightarrow \;\; \forall i \in I \quad f_i \, \llbracket {\mc F} \, \rrbracket \rightarrow f_i(x)$
\eenum

\etheo

\bproof
To see (\ref{theo conv induced 1}), assume that $f_i(x_n) \rightarrow f_i(x)$ for every $i \in I$ and let $U \in {\mc N} \!\left\{ x \right\}$. Then there is ${\mc G} \sqsubset {\mc S}$ such that ${\mc G} \neq \O$ and $x \in \bigcap {\mc G} \subset U$. Now let $G \in {\mc G}$. Further let $i \in I$ and $V \in \tp_i$ such that $G = f_i^{-1} \left[ V \right]$. Thus there is $N_G \in D$ such that $f_i(x_n) \in V$ for every $n \geq N_G$ by assumption. We may choose such $N_G$ for each $G \in {\mc G}$. It follows that $x_n \in U$ for $n \geq N$ where $N$ is the maximum of~$\left\{ N_G \, : \, G \in {\mc G} \right\}$. The converse follows by Theorem~\ref{theo global cont}.

To see (\ref{theo conv induced 2}), assume that $f_i \, \llbracket {\mc F} \, \rrbracket \rightarrow f_i(x)$ for every $i \in I$. Let $(y_m)$ be the net generated by~$\mc F$, $U \in {\mc N} \!\left\{ x \right\}$, and $i \in I$. There is $F \in {\mc F}$ such that $f_i \left[ F \right] \subset U$ by assumption. Now let $M = (y,F)$ with $y \in F$. It follows that $f_i(y_m) \in U$ for every $m \geq M$. Hence we have $y_m \rightarrow x$ by~(\ref{theo conv induced 1}), and thus ${\mc F} \rightarrow x$. Again the reverse implication follows by Theorem~\ref{theo global cont}.
\eproof

\section{Topological subspace}
\label{topological subspace}

In this Section a special case of inverse image topologies is analyzed, viz.\ topological subspaces.

\midvspace

\blede
\label{lede rel top}
\index{Topology!relative}
\index{Topological subspace}
\index{Subspace}
\index{Open in}
\index{Set!open in}
\index{Closed in}
\index{Set!closed in}
Given a topological space $\xi = (X,\tp)$ and a subset $A \subset X$, the system $\left\{ U \cap A \, : \, U \in \tp \right\}$ is a topology on~$A$. It is called the {\bf relative topology on}~$A$ and denoted by~$\tp \, | \, A$. The pair $(A,\tp \, | \, A)$ is called {\bf topological subspace of}~$\xi$, or short {\bf subspace of}~$\xi$, and denoted by~$\xi \, | \, A$. A member of~$\tp \, | \, A$ is also called {\bf open in}~$A$ whereas a member of~$\tp$ is also called {\bf open in}~$X$. Furthermore a $(\tp \, | \, A)$-closed set is also called {\bf closed in}~$A$, and a $\tp$-closed set is also called {\bf closed in}~$X$.
\elede

\bproof
Exercise.
\eproof

\bdefi
\index{Inclusion}
Given a set $X$ and a subset $A \subset X$, the map $j : A \longrightarrow X$, $j(x) = x$ ($x \in A)$, is called {\bf inclusion}. We also write $j : A \hookrightarrow X$.
\edefi

\midvspace

\blemm
\label{lemm top subspace}
Let $\xi = (X,\tp)$ be a topological space, $\alpha = (A,\tp_A)$ a subspace of~$\xi$, and $\mc C$ and ${\mc C}_A$ their respective systems of closed sets. Further let $\mc B$ be a base for~$\tp$, $\mc S$~a subbase for~$\tp$, $\mc D$~a base for the closed sets in~$X$, $\mc E$~a subbase for the closed sets in~$X$, $j : A \hookrightarrow X$ the inclusion, $x \in A$, and $B \subset A$. Then the following statements hold:

\benum
\item \label{lemm top subspace 1} We have $\tp_A = j^{-1} \, \llbracket \tp \, \rrbracket \,$. $\tp_A$~is generated by $\left\{ (j,\tp) \right\}$, and $j$ is $\tp_A$-$\tp$-continuous.
\item \label{lemm top subspace 2} The system $j^{-1} \, \llbracket {\mc B} \, \rrbracket$ is a base for~$\tp_A$.
\item \label{lemm top subspace 3} The system $j^{-1} \, \llbracket {\mc S} \, \rrbracket$ is a subbase for~$\tp_A$.
\item \label{lemm top subspace 4} Let $(Y,\tp_Y)$ be a topological space and $g : Y \longrightarrow A$ a map. Then $g$ is continuous iff $j \circ g$ is continuous.
\item \label{lemm top subspace 5} $\left\{ A \cap C \, : \, C \in {\mc C} \right\} \, = \, j^{-1} \, \llbracket {\mc C} \, \rrbracket \, = \, {\mc C}_A$
\item \label{lemm top subspace 6} The system $j^{-1} \, \llbracket {\mc D} \, \rrbracket$ is a base for~${\mc C}_A$.
\item \label{lemm top subspace 7} The system $j^{-1} \, \llbracket {\mc E} \, \rrbracket$ is a subbase for~${\mc C}_A$.
\item \label{lemm top subspace 8} $\big\{ A \cap U \, : \, U \in {\mc N}_{\xi} \!\left\{ x \right\} \!\big\} = \, j^{-1} \, \llbracket \, {\mc N}_{\xi} \!\left\{ x \right\} \rrbracket \, = \, {\mc N}_{\alpha} \!\left\{ x \right\}$
\item \label{lemm top subspace 9} $\mathrm{cl}_{\alpha} B \, = \, \mathrm{cl}_{\xi}(B)  \cap A$
\item \label{lemm top subspace 10} $\mathrm{int}_{\xi} B \, \subset \, \mathrm{int}_{\alpha} B$
\item \label{lemm top subspace 11} $\mathrm{bound}_{\alpha} B \, \subset \, \mathrm{bound}_{\xi} B$
\eenum

\elemm

\bproof
(\ref{lemm top subspace 1}) is obvious.

(\ref{lemm top subspace 2}), (\ref{lemm top subspace 3}), and (\ref{lemm top subspace 5}) to (\ref{lemm top subspace 7} are consequences of~(\ref{lemm top subspace 1}) and Corollary~\ref{coro gen top single map}.

(\ref{lemm top subspace 4}) follows by (\ref{lemm top subspace 1}) and Theorem~\ref{universalproperty}.

The first equation in~(\ref{lemm top subspace 8}) follows by definition of~$j$. Moreover we have 
\[
j^{-1} \, \llbracket \, {\mc N}_{\xi} \!\left\{ x \right\} \rrbracket \, \subset \, {\mc N}_{\alpha} \!\left\{ x \right\}
\]

by the continuity of~$j$. Conversely, let $U \in {\mc N}_{\alpha} \!\left\{ x \right\}$. We may choose $V \in \tp_A$ with $x \in V \subset U$, and $W \in \tp$ with $V = W \cap A$. Let $R = W \cup U$. Then $R \in {\mc N}_{\xi} \!\left\{ x \right\}$, and $j^{-1} \left[ R \right] = U$.

(\ref{lemm top subspace 9}) follows by (\ref{lemm top subspace 5}) and Lemma~\ref{lemm closure intersec} as follows:
\begin{eqnarray*}
\mathrm{cl}_{\alpha} B \!\!\! & = & \!\!\! \bigcap \big\{ C \in {\mc C}_A \, : \, B \subset C \big\}\\[.2em]
& = & \!\!\! \bigcap \big\{ D \cap A \, : \, B \subset D,\,  D \in {\mc C} \big\}\\[.2em]
& = & \!\!\! \bigcap \big\{ D \in {\mc C} \, : \, B \subset D \big\} \cap A \; = \; \mathrm{cl}_{\xi}(B) \cap A
\end{eqnarray*}

To see (\ref{lemm top subspace 10}), notice that a set $U \subset A$ that is open in~$X$ is also open in~$A$. Therefore we have:
\begin{eqnarray*}
\mathrm{int}_{\xi} B \!\!\! & = & \!\!\! \bigcup \big\{ U \in \tp \, : \, U \subset B \big\}\\[.2em]
& \subset & \!\!\! \bigcup \big\{ U \in \tp_A \, : \, U \subset B \big\} \, = \, \mathrm{int}_{\alpha} B
\end{eqnarray*}

Finally, in order to prove (\ref{lemm top subspace 11}), let $x \in \mathrm{bound}_{\alpha} B$, $U \in {\mc N}^{\, \mathrm{open}}_{\xi} \! \left\{ x \right\}$, and $V = U \cap A$. It follows that $V \in {\mc N}_{\alpha} \!\left\{ x \right\}$, and thus $V \cap B \neq \O$ and $V \cap \left( A \!\setminus\! B \right) \neq \O$. Therefore we have $U \cap B \neq \O$ and $U \cap \left( X \!\setminus\! B \right) \neq \O$. It follows that \mbox{$x \in \mathrm{bound}_{\xi} B$}.
\eproof

\brema
\label{rema open closed subspace}
Let $\xi = (X,\tp)$ be a topological space and $B \subset A \subset X$. If $B$ is open in~$X$, then $B$ is open in~$A$ by definition of the relative topology. If $B$ is closed in~$X$, then $B$ is closed in~$A$ by Lemma~\ref{lemm top subspace}~(\ref{lemm top subspace 5}). The converse of both implications is generally not true. However, it is true under additional assumptions, cf.\ Lemmas~\ref{lemm closed subspace} and~\ref{lemm open subspace}.
\erema

\brema
\label{rema iterated subspace}
Let $\xi = (X,\tp)$ be a topological space and $B \subset A \subset X$. We have $\xi \, | \, B = (\xi \, | \, A) \, | \, B$, which is a consequence of Definition~\ref{lede rel top}. It also follows by Lemmas~\ref{induced assoc} and~\ref{lemm top subspace}~(\ref{lemm top subspace 1}).
\erema

\blemm
\label{lemm standard top r r+}
Let $\tp$ and $\tp_+$ be the standard topologies on~${\mathbb R}$ and ${\mathbb R}_+$\,, respectively. We have $\tp_+ = \tp \, | \, {\mathbb R}_+$\,. Further let $A \subset {\mathbb R}_+$\,. Then $\tp_+ | A  = \tp \, | A$.
\elemm

\bproof
To see the first claim note that
\[
{\mc S} = \big\{ \left] -\infty, x \right[\, , \, \left] x, \infty \right[\, : \, x \in {\mathbb R} \big\}
\]

is a subbase for~$\tp$ and
\[
{\mc S}_+ = \big\{ \left] -\infty, x \right[\, , \, \left] x, \infty \right[\, : \, x \in {\mathbb R}_+ \big\} \cup \left\{ {\mathbb R}_+ \right\}
\]

is a subbase for~$\tp_+$\,. These two systems are related by $j^{-1} \, \llbracket {\mc S} \rrbracket = {\mc S}_+$\, where $j : {\mathbb R}_+ \hookrightarrow {\mathbb R}$ is the inclusion. The claim follows by~Theorem~\ref{lemm top subspace}~(\ref{lemm top subspace 3}).

Now the second claim follows by Remark~\ref{rema iterated subspace}.
\eproof

The following three Lemmas show that the concept of topological subspace does not lead to any complications when considering convergent nets or filters and continuous functions.

\midvspace

\blemm
\label{conv subspace net}
Let $\xi = (X,\tp)$ be a topological space, $(A,{\mc A})$ a subspace of~$\xi$, $(x_n)$ a net in~$A$, and $x \in A$. Then $x_n \rightarrow x$ with respect to~$\mc A$ iff $x_n \rightarrow x$ with respect to~$\tp$.
\elemm

\bproof
This is a special case of Theorem~\ref{theo conv induced}~(\ref{theo conv induced 1}).
\eproof

\blemm
\label{con vsubspace filter}
Let $\xi = (X,\tp)$ be a topological space, $(A,{\mc A})$ a subspace of~$\xi$, $x \in A$, and $\mc G$ a filter on~$A$. Then $\mc G$ is a filter base for a filter on~$X$, say~$\mc F$. Moreover the following statements are equivalent:

\benum
\item \label{con vsubspace filter 1} ${\mc G} \rightarrow x$ with respect to~$\mc A$
\item \label{con vsubspace filter 2} ${\mc G} \rightarrow x$\, with respect to~$\tp$, i.e.\ ${\mc G}$ is considered as filter base on~$X$
\item \label{con vsubspace filter 3} ${\mc F} \rightarrow x$ with respect to~$\tp$
\eenum

\elemm

\bproof
The equivalence of (\ref{con vsubspace filter 1}) and~(\ref{con vsubspace filter 2}) is a special case of Theorem~\ref{theo conv induced}~(\ref{theo conv induced 2}). The equivalence of (\ref{con vsubspace filter 2}) and~(\ref{con vsubspace filter 3}) follows by definition.
\eproof

\blemm
\label{continuity subspace}
Let $(X,\tp_X)$ and $(Y,\tp_Y)$ be topological spaces, $(A,\tp_A)$ a subspace of $(X,\tp_X)$, $(B,\tp_B)$ a subspace of $(Y,\tp_Y)$, and $f : X \longrightarrow Y$ a map with $f \left[ X \right] \subset B$. Further we define the map $g : X \longrightarrow B$, $g(x) = f(x)$. The following statements hold:
\benum
\item If $f$ is continuous, then $f \, | \, A$ is continuous.
\item $f$ is continuous iff $g$ is continuous.
\eenum

\elemm

\bproof
Exercise.
\eproof

The following two Lemmas state under which conditions the reverse implications of Remark~\ref{rema open closed subspace} are true.

\midvspace

\blemm
\label{lemm closed subspace}
Let $\xi = (X,\tp)$ be a topological space, $A$ a closed subset of~$X$, $\alpha = \xi \, | \, A$, and $B \subset A$. The following statements hold:

\benum
\item \label{lemm closed subspace 1} $B$ is closed in $X$ iff it is closed in $A$.
\item \label{lemm closed subspace 2} $\mathrm{cl}_{\alpha} B \, = \, \mathrm{cl}_{\xi} B$
\eenum

\elemm

\bproof
(\ref{lemm closed subspace 1}) follows by Lemma~\ref{lemm top subspace}~(\ref{lemm top subspace 5}) and Remark~\ref{rema open closed subspace}.

To see (\ref{lemm closed subspace 2}), notice that $\mathrm{cl}_{\xi} B \subset \mathrm{cl}_{\xi} A = A$. Thus we have $\mathrm{cl}_{\alpha} B = \mathrm{cl}_{\xi} B$ by Lemma~\ref{lemm top subspace}~(\ref{lemm top subspace 9}).
\eproof

\blemm
\label{lemm open subspace}
Let $\xi = (X,\tp)$ be a topological space, $A \in \tp$, $\alpha = \xi \, | \, A$, and $B \subset A$. The following statements hold:

\benum
\item \label{lemm open subspace 1} $B$ is open in $X$ iff it is open in~$A$.
\item \label{lemm open subspace 2} $\mathrm{int}_{\alpha} B \, = \, \mathrm{int}_{\xi} B$
\item \label{lemm open subspace 3} $\mathrm{bound}_{\alpha} B \, = \, \mathrm{bound}_{\xi}(B) \cap A$
\eenum

\elemm

\bproof
(\ref{lemm open subspace 1}) is obvious.

To show (\ref{lemm open subspace 2}), let $x \in \mathrm{int}_{\alpha} B$ and $U$ open in~$A$ with $x \in U \subset B$. Then $U$ is open in~$X$ by~(\ref{lemm open subspace 1}), and thus $x \in \mathrm{int}_{\xi} B$. The converse follows by Lemma~\ref{lemm top subspace}~(\ref{lemm top subspace 10}).

To show (\ref{lemm open subspace 3}), notice that $\mathrm{bound}_{\alpha} B \subset \mathrm{bound}_{\xi}(B) \cap A$ by Lemma~\ref{lemm top subspace}~(\ref{lemm top subspace 11}). Conversely, let $x \in \mathrm{bound}_{\xi}(B) \cap A$ and $U \in {\mc N}^{\, \mathrm{open}}_{\alpha} \!\left\{ x \right\}$. It follows by~(\ref{lemm open subspace 1}) that $U \in {\mc N}^{\, \mathrm{open}}_{\xi} \!\left\{ x \right\}$, and hence $U \cap B \neq \O$ and
\[
\O \, \neq \, U \cap \left( X \!\setminus\! B \right) \, = \, U \cap \left( A \!\setminus\! B \right)
\]

Thus $x \in \mathrm{bound}_{\alpha} B$.
\eproof

\blede
\label{def pseudometric subspace}
\index{Pseudo-metric!relative}
\index{Subspace}
\index{Pseudo-metric!subspace}
Given a pseudo-metric space $\xi = (X, d)$ and a subset $A \subset X$, the restriction $d \, | (A \times A)$ is a pseudo-metric on~$A$. It is called the {\bf relative pseudo-metric} and denoted by~$d \, | A$. The pseudo-metric space $(A, d \, | A)$ is called {\bf pseudo-metric subspace of}~$\xi$, or short {\bf subspace of}~$\xi$, and denoted by~$\xi \, | A$.
\elede

\bproof
This follows by Definition~\ref{metric}.
\eproof

\blemm
Given a metric space $\xi = (X,d)$ and a subset $A \subset X$, the subspace $\xi \, | \, A$ is a metric space.
\elemm

\bproof
This follows by Definition~\ref{metric}.
\eproof

The following Lemma states that the generation of a topology from a pseudo-metric commutes with the formation of a subspace.

\midvspace

\blemm
Given a pseudo-metric space $\xi = (X, d)$ and a subset $A \subset X$, we have $\tau(d) \, | A = \tau (d \, | A)$.
\elemm

\bproof
We define for every $r \in \,\left] 0, \infty \right[\,$ and $x \in X$:
\[
B(x,r) = \big\{ y \in X \, : \, d(x,y) < r \big\}
\]

Moreover, let
\[
{\mc B} = \big\{ B(x,r) \cap A \, : \, r \in \,\left] 0, \infty \right[,\, x \in X \big\}
\]

Then $\mc B$ is a base for for~$\tau (d) \, |A$ by Lemma~\ref{lemm top subspace}~(\ref{lemm top subspace 2}). Further let $d_A = d \, | A$ and
\[
{\mc B}_A = \big\{  B(x,r) \cap A \, : \, r \in \,\left] 0, \infty \right[,\, x \in A \big\}
\]

Then ${\mc B}_A$ is a base for~$\tau (d_A)$. We have ${\mc B}_A \subset {\mc B}$, and thus $\tau (d_A) \subset \tau(d) \, | A$. Conversely, let $r \in \,\left] 0, \infty \right[\,$, $x \in X$, and $y \in B(x,r) \cap A$. Since $B(x,r)$ is $\tau(d)$-open, there is $s \in \,\left] 0, \infty \right[\,$ such that $B(y,s) \subset B(x,r)$. Since
$B(y,s) \cap A \, \in \, {\mc B}_A$, it follows that ${\mc B} \subset_{\stackrel{}{\Phi}} {\mc B}_A$, and thus $\tau(d) \, | A \subset \tau (d_A)$.
\eproof

\bdefi
Let $(X,\tp)$ be a topological space and ${\mc A} \subset {\mc P}(X)$. $\mc A$~is called {\bf locally finite} if for every $x \in X$ there is a neighborhood $U \in {\mc N} \!\left\{ x \right\}$ and ${\mc B} \sqsubset {\mc A}$ such that $U \cap A = \O$ for every $A \in {\mc A} \!\setminus\! {\mc B}$. $\mc A$~is called {\bf locally discrete} if for every $x \in X$ there is a neighborhood $U \in {\mc N} \!\left\{ x \right\}$ and ${\mc B} \sqsubset {\mc A}$ such that $U \cap A = \O$ for every $A \in {\mc A} \!\setminus\! {\mc B}$ and either ${\mc B} \sim 0$ or ${\mc B} \sim 1$.
\edefi

\blemm
\label{lemma sigma locally finite closure}
Let $(X,\tp)$ be a topological space, $I$~an index set, and $A_i \subset X$ ($i \in I$). If $\left\{ A_i \, : \, i \in I \right\}$ is locally finite, then we have \, $\overline{\bigcup_{i \in I} A_i} = \bigcup_{i \in I} \overline{A_i}$\,.
\elemm

\bproof
Assume the stated condition. Let $x \in \overline{\bigcup_{i \in I} A_i}$. There exist $U \in {\mc N} \!\left\{ x \right\}$ and $J \sqsubset I$ such that $U \cap A_i = \O$ for every $i \in I \!\setminus\! J$. Hence $V \cap \bigcup_{i \in J} A_i \neq \O$ for every $V \in {\mc N} \!\left\{ x \right\}$ by Theorem~\ref{theo closure}~(\ref{theo closure 1a}). It follows that
\[
x \, \in \, \overline{\bigcup_{i \in J} A_i} \, = \, \bigcup_{i \in J} \overline{A}_i \, \subset \, \bigcup_{i \in I} \overline{A}_i
\]

by Lemma~\ref{lemma system closure rules}~(\ref{lemma system closure rules 6}). The converse follows by Lemma~\ref{lemma system closure rules}~(\ref{lemma system closure rules 7}).
\eproof

\btheo
\label{theo cont closed system}
Let $(X,\tp)$ and $(Y,\tp_Y)$ be topological spaces, $I$~an index set, closed sets $A_i \subset X$ ($i \in I)$ such that $\left\{ A_i \, : \, i \in I \right\}$ is locally finite and $\bigcup_{i \in I} A_i = X$, and $f : X \longrightarrow Y$ a function. $f$~is continuous iff $f \, | A_i$ is continuous for every $i \in I$.
\etheo

\bproof
If $f$ is continuous, then $f \, | A_i$ is continuous for every $i \in I$ by Lemma~\ref{continuity subspace}.

Conversely, assume that $f \, | A_i$ is continuous for every $i \in I$. For every $i \in I$ we define $\tp_i = \tp \, | A_i$\,.

Let $B \subset X$. Note that
\[
\mathrm{cl}_{\tp} B \, = \, \mathrm{cl}_{\tp} \left( {\textstyle \bigcup_{i \in I}} (B \cap A_i) \right) \, = \, \bigcup_{i \in I} \mathrm{cl}_{\tp} (B \cap A_i)
\]

by Lemma~\ref{lemma sigma locally finite closure}. Moreover, for every $i \in I$, we have $\mathrm{cl}_{\tp(i)} (B \cap A_i) = \mathrm{cl}_{\tp} (B \cap A_i)$ by Lemma~\ref{lemm closed subspace}~(\ref{lemm closed subspace 2}). It follows that
\begin{eqnarray*}
f \left[\, \mathrm{cl}_{\tp} B \, \right] \!\! & = & \! \bigcup_{i \in I} f \left[\, \mathrm{cl}_{\tp} (B \cap A_i) \right]\\[.2em]
 & = & \! \bigcup_{i \in I} \, (f \, | A_i) \left[\,  \mathrm{cl}_{\tp(i)} (B \cap A_i) \right]\\[.2em]
 & = & \! \bigcup_{i \in I} \mathrm{cl} \, \big( (f \, | A_i) \left[ B \cap A_i \right] \big)\\[.2em]
 & \subset & \! \mathrm{cl} \left( \, {\textstyle \bigcup_{i \in I}} \, f \left[ B \cap A_i \right] \, \right) \, = \; \mathrm{cl} \left( f \left[ B \right] \right)
\end{eqnarray*}

where the third equation is a consequence of Theorem~\ref{theo cont closure int}~(\ref{theo cont closure int 2}) and the fourth line follows by Lemma~\ref{lemma system closure rules}~(\ref{lemma system closure rules 7}). The continuity of~$f$ now follows by the same Theorem.
\eproof

\brema
\label{rema rel top int r}
Let $<$ be the standard ordering on~${\mathbb R}$, $\tp$~the standard topology on~${\mathbb R}$, $a, b \in {\mathbb R}$ with $a < b$, $A = \left[ a, b \right]$, $<_A$~the restriction of~$<$ to~$A$, $\tp_B$~the interval topology of~$(A,<_A)$, and $\tp_A = \tp \, | A$. We denote by subscript $A$ intervals with respect to the ordering~$<_A$\,. All other intervals refer to the ordering~$<$. The system
\[
{\mc S} \, = \, \big\{ \left] -\infty, x \right[\, , \;\left] x, \infty \right[\; : \, x \in {\mathbb R} \big\} \cup \left\{ \O \right\}
\]

is a subbase for~$\tp$. Further the system
\[
{\mc R} \, = \, \big\{ \left] -\infty, x \right[_{\,A}\, , \;\left] x, \infty \right[_{\,A}\, : \, x \in A \big\} \cup \left\{ A \right\}
\]

is a subbase for~$\tp_B$. We have
\begin{eqnarray*}
{\mc R} \!\!\! & = & \!\!\!\big\{ \left[ a, x \right[\; : \, x \in A \!\setminus\! \left\{ a \right\}  \!\big\} \, \cup \, \big\{ \, \left] x, b \right] \, : \, x \in A \!\setminus\! \left\{ b \right\} \!\big\} \, \cup \, \big\{ \! A, \O \big\}\\
 & = & \!\!\!\big\{ S \cap A \, : \, S \in {\mc S} \big\}
\end{eqnarray*}

It follows by Lemma~\ref{lemm top subspace}~(\ref{lemm top subspace 3}) that $\tp_A = \tp_B$\,.
\erema

\bexam
Let $\xi = ({\mathbb R},\tp)$ where $\tp$ is the standard topology, $a, b \in {\mathbb R}$ with $a < b$, $A = \, \left] a, b \right]$, and $\alpha = (A, \tp_A) = \xi \, | A$ the topological subspace. We define
\begin{eqnarray*}
{\mc S}_+ \!\!\! & = & \!\!\! \big\{ \left] x, b \right] \; : \, x \in A \!\setminus\! \left\{ b \right\} \big\}\\[.2em]
{\mc R}_+ \!\!\! & = & \!\!\! \big\{ \left] x, b \right] \; : \, x \in {\mathbb D} \cap \left( A \!\setminus\! \left\{ b \right\} \right) \big\} \cup \left\{ A \right\}\\[.2em]
{\mc B} \!\!\! & = & \!\!\! \big\{ \left] x, y \right[ \; : \, x, y \in {\mathbb R},\, x < y \big\} \cup \left\{ \O \right\}\\[.2em]
{\mc B}_A \!\!\! & = & \!\!\! \big\{ \left] x, y \right[ \; : \, x, y \in A,\, x < y \big\} \cup \left\{ \O \right\}\\[.2em]
{\mc A} \!\!\! & = & \!\!\! \big\{ \left] x, y \right[ \; : \, x, y \in {\mathbb D},\, x < y \big\} \cup \left\{ \O \right\}\\[.2em]
{\mc A}_A \!\!\! & = & \!\!\! \big\{ \left] x, y \right[ \; : \, x, y \in {\mathbb D} \cap A,\, x < y \big\} \cup \left\{ \O \right\}
\end{eqnarray*}

Each of the systems $\mc B$ and $\mc A$ is a base for~$\tp$ by Remark~\ref{rema r interval rel}. Each of the systems ${\mc B}_A \cup {\mc S}_+$ and ${\mc A}_A \cup {\mc R}_+$ is a base for~$\tp_A$.

Further let $c \in {\mathbb R}$ with $a < c < b$. Then the following statements hold:
\[
\ba{ll}

]c,b] \, \notin \, {\mc N}_{\xi} \!\left\{ b \right\}, & ]c,b] \, \in \, {\mc N}_{\alpha} \!\left\{ b \right\}, \\[.8em]
\mathrm{cl}_{\xi} \; ]a,c[ \;\; = \; [a,c]\,, & \mathrm{cl}_{\alpha} \; ]a,c[ \;\; = \;\; ]a,c]\,, \\[.8em]
\mathrm{cl}_{\xi} \; ]c,b[ \;\; = \; [c,b]\,, & \mathrm{cl}_{\alpha} \; ]c,b[ \;\; = \; [c,b]\,, \\[.8em]
\mathrm{int}_{\xi} \, [c,b] \; = \;\; ]c,b[ \,, &  \mathrm{int}_{\alpha} \, [c,b] \; = \;\; ]c,b]\,, \\[.8em]
\mathrm{bound}_{\xi}\; ]a,c[ \;\; = \; \left\{ a, c \right\}, &  \mathrm{bound}_{\alpha} \; ]a,c[ \;\;  = \; \left\{ c \right\}
\ea
\]

\eexam

\bexam
Let $\tp$ be the standard topology on~${\mathbb R}$. The set ${\mathbb D}$ is neither $\tp$-open nor $\tp$-closed by Lemma~\ref{lemm D dense}. Further let $\tp_D$ be the relative topology on~${\mathbb D}$, and $a, b \in {\mathbb R}$ with $a < b$. The set $[a,b] \cap {\mathbb D}$ is $\tp_D$-closed. If $a, b \notin {\mathbb D}$, then this set is also $\tp_D$-open.
\eexam

\section{Product topology}
\label{product topology}

In this Section we consider another important special case of inverse image topologies: product topologies.

\midvspace

\bdefi
\label{defprodtop}
\index{Product topology}
\index{Topology!product}
\index{Product space}
Let $\xi_i = (X_i,\tp_i)$ $(i \in I)$ be topological spaces where $I$ is an index set, $X = \bigtimes_{\!\! i \in I}\, X_i$\,, and $p_i : X \longrightarrow X_i$ the projections. The topology $\tp = \tau (\left\{ (p_i,\tp_i) \, : \, i \in I \right\})$ is called {\bf product topology} and denoted by~$\prod_{i \in I} \tp_i$\,. The topological space $\xi = (X,\tp)$ is called {\bf product topological space}, or short {\bf product space}, and denoted by~$\prod_{i \in I} \xi_i$\,. If $I = \sigma(n) \!\setminus\! m$ for some $m, n \in \naturalnumbers$ with $m < n$, then we also write $\prod_{k = m}^n \tp_k$ for~$\tp$, and $\prod_{k = m}^n \xi_k$ for~$\xi$. If $I = \naturalnumbers \!\setminus\! m$ for some $m \in \naturalnumbers$, then we also write $\prod_{k = m}^{\infty} \tp_k$ for~$\tp$, and $\prod_{k = m}^{\infty} \xi_k$ for~$\xi$.
\edefi

While the symbol $\bigtimes$ always denotes the Cartesian product, the symbol $\prod$ has several meanings, e.g.\ it denotes the product topology and the product space as well as product nets. More meanings of~$\prod$ are introduced below. In each occurrence we ensure that the correct interpretation is evident from the context.
 
We also remind the reader that we often use the notation $x_i$ for $p_i(x)$, where $i \in I$ and $x \in X$, as introduced in the general case of generated topologies in Lemma and Definition~\ref{definitialtop}.

By definition the projections $p_i$ are continuous maps. Furthermore the following result holds.

\midvspace

\blemm
\label{lemm projection open}
With definitions as in Definition~\ref{defprodtop} the projections $p_i : X \longrightarrow X_i$ are open maps.
\elemm

\bproof
Let $i \in I$, $U \in \tp$, and $U_i = p_i \left[ U \right]$. To show that $U_i$ is open let $r \in U_i$\,. We may choose $x \in U$ such that $p_i(x) = r$. Further let $K \sqsubset I$ and $V_k \in \tp_k$ ($k \in K$) such that $x \in \bigcap_{k \in K}\, p_k^{-1} \left[ V_k \right] \subset U$. If $i \notin K$, then we have $U_i = X_i$\,. If $i \in K$, then it follows that $r \in V_i \subset U_i$\,.
\eproof

\blemm
Let $I$ be an index set. For every $i \in I$ let $(X_i,R_i)$ be a pre-ordered space such that $R_i$ has full field. Further let $\tp_i$~($i \in I$) be the respective interval topologies and $X = \bigtimes_{\!\! i \in I}\, X_i$\,. Then, for every $i \in I$, $S_i = p_i^{-1} \left[ R_i \right]$ is a pre-ordering on~$X$ that has full field. Let ${\mc S} = \left\{ S_i \, : \, i \in I \right\}$, and $\tp$ be the $\mc S$-interval topology. Then we have $\tp = \prod_{i \in I} \tp_i$\,.
\elemm

\bproof
For every $i \in I$, $S_i$~is a pre-ordering by Example~\ref{exam product pre-orderings} and clearly has full field. We have
\begin{eqnarray*}
 & & \left\{ \;\left] -\infty, x \right[_{\, S(i)},\, \;\left] x, \infty \right[_{\, S(i)} \, : \, x \in X,\, i \in I \right\} \cup \left\{ \O \right\}\\[.2em]
& = & \left\{ p_i^{-1}\left[\;\left] -\infty, x_i \right[_{\, R(i)}\right]\!,\; p_i^{-1}\left[\;\left] x_i, \infty \right[_{\, R(i)}\right] \, : \, x \in X,\, i \in I \right\} \cup \left\{ \O \right\}\\[.2em]
& = & \left\{ p_i^{-1}\left[\;\left] -\infty, r \right[_{\, R(i)}\right]\!,\; p_i^{-1}\left[\;\left] r, \infty \right[_{\, R(i)}\right] \, : \, r \in X_i\,,\, i \in I \right\} \cup \left\{ \O \right\}
\end{eqnarray*}

The first expression is a subbase for~$\tp$, and the last is a subbase for the product topology.
\eproof

\brema
\label{rema standard prod top}
Let $\tp$, $\tp_+$\,, $\tp^n$, and $\tp^n_+$ be the standard topologies on~${\mathbb R}$, ${\mathbb R}_+$\,, ${\mathbb R}^n$, and ${\mathbb R}^n_+$\,, respectively, where $n \in \naturalnumbers$, $n > 0$. We have $\tp^n = \prod_{k = 1}^n \tp$ and $\tp^n_+ = \prod_{k = 1}^n \tp_+$\,.
\erema

Next we demonstrate that an iterated product of topological spaces is essentially a product of those topological spaces.

\midvspace

\blede
Let $\xi_j = ( X_j, \tp_j )$ $(j \in J_i, \, i \in I)$ be topological spaces where $I$ is an index set and $J_i$ $(i \in I)$ are disjoint index sets. Further let $K = \bigcup \left\{ J_i \, : \, i \in I \right\}$. Then $\prod_{i \in I} \left(\prod_{j \in J_i} \xi_j \right)$ and $\prod_{j \in K} \xi_j$ are homeomorphic.
\elede

\bproof
For every $i \in I$ let $(Y_i,\tp_i) = \prod_{j \in J_i}\, (X_j,\tp_j)$. Further let $(Y,\tp_Y) = \prod_{i \in I}\, (Y_i,\tp_i)$ and $(X,\tp_X) = \prod_{j \in K}\, (X_j,\tp_j)$. We define the map
\begin{center}
\begin{tabular}{l}
$f \; : \; X \longrightarrow Y$\;,\\[.8em]
$\Big(\big(f(h)\big)(i)\Big)(j) = h(j)$ \quad for every $i \in I$ and $j \in J_i \;$
\end{tabular}
\end{center}

Then $f$ is a bijection by Remark~\ref{rema it Cartesian product}. For every $i \in I$ the topology $\tp_i$ is generated by $\left\{ (p_{ij},\tp_j) \, : \, j \in J_i \right\}$ where \mbox{$p_{ij} : Y_i \longrightarrow X_j$} $(j \in J_i)$ are the projections. The topology $\tp_Y$ is generated by $\left\{ (p_i,\tp_i) \, : \, i \in I \right\}$ where \mbox{$p_i : Y \longrightarrow Y_i$} ($i \in I$) are the projections. Furthermore $\tp_X$ is generated by $\left\{ (q_j,\tp_j) \, : \, j \in K \right\}$ where \mbox{$q_j : X \longrightarrow X_j$} ($j \in K$) are the projections. Hence
\[
{\mc S}_X = \big\{ q_j^{-1} \left[ U \right] \, : \, j \in K,\, U \in \tp_j \big\}
\]

is a subbase for~$\tp_X$. Since $q_j \circ f^{-1} = p_{ij} \circ p_i$ for every $i \in I,\, j \in J_i$\,, the topology $\tp_Y$ is generated by $\left\{ \left(q_j \circ f^{-1}, \tp_j \right) \, : \, j \in K \right\}$ by Lemma~\ref{induced assoc}. This means that $f \, \llbracket {\mc S}_X \rrbracket$ is a subbase for~$\tp_Y$. Thus $f$ is a homeomorphism.
\eproof

The following result is a direct consequence of the more general Theorem in Section~\ref{inverse image topology}, however, it is particularly important in the case of product  spaces.

\midvspace

\bcoro
\label{coro conv product filter}
With definitions as in Definition~\ref{defprodtop} let $(x_n)$~be a net in~$X$, $\mc F$~a filter on~$X$, and $x \in X$. The following statements hold:
\benum
\item \label{coro conv product filter 1} $x_n \rightarrow x \;\; \Longleftrightarrow \;\; \forall i \in I \quad p_i(x_n) \rightarrow p_i(x)$
\item \label{coro conv product filter 2} ${\mc F} \rightarrow x \;\; \Longleftrightarrow \;\; \forall i \in I \quad p_i \, \llbracket {\mc F} \, \rrbracket \rightarrow p_i(x)$
\eenum

\ecoro

\bproof
This follows by Theorem~\ref{theo conv induced}.
\eproof

\bcoro
\label{coro conv prod net}
Let $I$ be an index set. For every $i \in I$ let $(X_i,\tp_i)$ be a topological space, $A_i \subset X_i$\,, and $(x^i_n)$ a net in~$A_i$\,. Further let $(X,\tp) = \prod_{i \in I}\, (X_i,\tp_i)$ and $p_i : X \longrightarrow X_i$ ($i \in I$) be the projections. Moreover, let $(x_r) = \prod_{i \in I}\, (x^i_n)$ and $x \in X$. Then we have
\[
x_r \rightarrow x \quad \Longleftrightarrow \quad \forall i \in I \quad x^i_n \rightarrow p_i(x)
\]

\ecoro

\bproof
This is a consequence of Corollary~\ref{coro conv product filter} and Lemma~\ref{lemm prod net conv}.
\eproof

\btheo
\label{theo product properties}
Let $I$ be an index set. For every $i \in I$ let $(X_i,\tp_i)$ be a topological space and $A_i \subset X_i$\,. Moreover let $(X,\tp) = \prod_{i \in I}\, (X_i,\tp_i)$ and \mbox{$p_i : X \longrightarrow X_i$} ($i \in I$) be the projections. Further we define $A = \bigtimes_{\!\! i \in I}\, A_i$\,. The following statements hold:

\benum
\item \label{theo product properties 1} $\overline{A} \, = \bigtimes_{\!\! i \in I}\, \overline{A_i}$
\item \label{theo product properties 2} $A$ is closed iff $A_i$ is closed for every $i \in I$.
\item \label{theo product properties 3} $\O \neq A \in \tp \quad \Longleftrightarrow \\[.2em]
{} \quad \big( \forall i \in I \quad \O \neq A_i \in \tp_i \, \big) \;\; \wedge \;\; \big( \exists K \sqsubset I \quad \forall j \in I \!\setminus\! K \quad A_j = X_j \big)$
\item \label{theo product properties 4} $A^{\circ} \subset \bigtimes_{\!\! i \in I}\, A_i^{\circ}$
\item \label{theo product properties 5} If $A^{\circ} \neq \O$, then we have:\\[.2em]
$A^{\circ} = \bigtimes_{\!\! i \in I}\, A_i^{\circ} \quad \Longleftrightarrow \quad \exists K \sqsubset I \quad \forall i \in I \!\setminus\! K \quad A_i = X_i$
\item \label{theo product properties 6} $\tp \, | A \, = \, \prod_{i \in I} \left( \tp_i \, | A_i \right)$
\eenum

\etheo

\bproof
Notice that (\ref{theo product properties 1}) follows by
Theorem~\ref{theo closure}~(\ref{theo closure 2}) and Corollaries~\ref{coro conv product filter} and~\ref{coro conv prod net}.

(\ref{theo product properties 2})~is a consequence of~(\ref{theo product properties 1}).

To show~(\ref{theo product properties 3}), assume that $\O \neq A \in \tp$. For every $i \in I$, we clearly have $A_i \neq \O$, and $A_i = p_i \left[ A \right] \in \tp_i$ since $p_i$ is open by Lemma~\ref{lemm projection open}. Let $x \in A$. We may choose $K \sqsubset I$ and, for each $i \in K$, $U_i \in \tp_i$ such that $x \in U \subset A$ where $U = \bigcap_{i \in K} p_i^{-1} \left[ U_i \right]$. Then we obtain $X_i = p_i \left[ U \right] \subset A_i$ for every $i \in I \!\setminus\! K$. The reverse implication holds by definition of the product topology.

To see~(\ref{theo product properties 4}), let $x \in A^{\circ}$. There is $U \in \tp$ such that $x \in U \subset A$. For every $i \in I$, we have $p_i(x) \in p_i \left[ U \right] \subset A_i$ and $p_i \left[ U \right] \in \tp_i$ since $p_i$ is open. Hence $p_i(x) \in  A_i^{\circ}$ for every $i \in I$.

In order to prove~(\ref{theo product properties 5}), first note that we clearly have $A \supset \bigtimes_{\!\! i \in I}\, A_i^{\circ}$. Assume that $A^{\circ} \neq \O$. It follows that $A_i^{\circ} \neq \O$ for every $i \in I$ by~(\ref{theo product properties 4}).

If there is $K \sqsubset I$ such that $A_i = X_i$ ($i \in I \!\setminus\! K$), then $\bigtimes_{\!\! i \in I}\, A_i^{\circ}$ is open by~(\ref{theo product properties 3}) and we obtain $A^{\circ} \supset \bigtimes_{\!\! i \in I}\, A_i^{\circ}$. It follows that $A^{\circ} = \bigtimes_{\!\! i \in I}\, A_i^{\circ}$. The converse follows by~(\ref{theo product properties 3}).

To show~(\ref{theo product properties 6}), let $i \in I$ and $U \in \tp_i$\,. We have
\[
p_i^{-1} \left[ U \right] \, \cap \, A \, = \, \left( p_i \, | \, A \right)^{-1} \left[ \, U \cap A_i \right]
\]

Thus there is a subbase for $\tp \, | A$ that is a subbase for $\prod_{i \in I} \left( \tp_i \, | A_i \right)$ as well.
\eproof

\bcoro
\label{coro standard top subspace}
Let $\tp$, $\tp_+$\,, $\tp^n$, and $\tp^n_+$ be the standard topologies on~${\mathbb R}$, ${\mathbb R}_+$\,, ${\mathbb R}^n$, and ${\mathbb R}^n_+$\,, respectively, where $n \in \naturalnumbers$, $n > 0$. We have $\tp^n_+ = \tp^n \, | \, {\mathbb R}^n_+$\,.
\ecoro

\bproof
This is shown by Remark~\ref{rema standard prod top}, Lemma~\ref{lemm standard top r r+}, and Theorem~\ref{theo product properties}~(\ref{theo product properties 6}) as follows:
\[
\tp^n_+ \, = \, \prod_{k = 1}^n \tp_+ \, = \, \prod_{k = 1}^n \left( \tp \, | \, {\mathbb R}_+ \right) \, = \, \tp^n \, | \, {\mathbb R}^n_+
\]

\eproof

\blemm
\label{lemmm cont prod}
With definitions as in Definition~\ref{defprodtop}, let $(Y,\tp_Y)$ be a topological space and $g : Y \longrightarrow X$ be a map. Then $g$ is continuous iff $p_i \circ g$ is continuous for every $i \in I$.
\elemm

\bproof
This follows by Theorem~\ref{universalproperty}.
\eproof

\blemm
Let $I$ be an index set and, for each $i \in I$, let $(X_i,\tp_i)$ and $(Y_i,\tp_i')$ be topological spaces and $f_i : X_i \longrightarrow Y_i$ a map. Further let $X = \bigtimes_{\!\! i \in I}\, X_i$ and $Y = \bigtimes_{\!\! i \in I}\, Y_i$\,, and the map $f : X \longrightarrow Y$ be defined by $\big(f(x)\big)_i = f_i(x_i)$ ($i \in I$). Then $f$ is continuous iff $f_i$ is continuous for every $i \in I$.
\elemm

\bproof
This follows by Theorem~\ref{theo global cont}~(\ref{theo global cont 6}) and Corollary~\ref{coro conv product filter}~(\ref{coro conv product filter 1}).
\eproof

We recall that the space of all functions from a set~$X$ to a set~$Y$, written $Y^X$, is identical to the Cartesian product $\bigtimes_{\!\! x \in X}\, Y$ with equal factors~$Y$. We often encounter subsets of functions $F \subset Y^X$. In such situations the following definition is convenient.

\midvspace

\bdefi
\label{evaluation function}
\index{Evaluation function}
Let $X$ and $Y$ be two sets, $F \subset Y^X$, and $p_x : Y^X \longrightarrow Y$ ($x \in X$) the projections, i.e.\ $p_x(f) = f(x)$ for every $x \in X$ and every $f \in Y^X$. Given $z \in X$, the restriction $q_z = p_z \, | \, F$ is called {\bf evaluation function at}~$z$.
\edefi

\brema
With definitions as in Definition~\ref{evaluation function}, we have $f(x) = q_x(f)$ for every $f \in F$ and $x \in X$. Let $j : F \hookrightarrow Y^X$. Then $q_x = p_x \circ j$ ($x \in X$).
\erema

Since the set of functions~$F$ in Lemma and Definition~\ref{evaluation function} is a subset of a Cartesian product, the concepts of relative topology and product topology may be used to define a topology on~$F$.

\midvspace

\bdefi
\label{pointwise}
\index{Pointwise convergence}
\index{Convergence!pointwise}
\index{Topology!of pointwise convergence}
Let $X$ be a set, $(Y,\tp_Y)$ a topological space, $\tp$ the product topology on~$Y^X$, and $(F,\tp_F)$ a subspace of $\left( Y^X, \tp \right)$. $\tp_F$~is called {\bf topology of pointwise convergence}.
\edefi

\brema
\label{rema pointw conv gen}
Let $X$ be a set, $(Y,\tp_Y)$ a topological space, $F \subset Y^X$, $q_x$~($x \in X$) the evaluation functions with domain~$F$, and $\tp_F$ the topology of pointwise convergence on~$F$. We have $\tp_F = \tau \left(\left\{(q_x,\tp_Y) \, : \, x \in X \right\}\right)$ by Lemma~\ref{induced assoc}. Moreover, the functions $q_x$ ($x \in X)$ are $\tp_F\,$-$\tp_Y\,$-continuous.
\erema

\bexam
${\mathbb R}^{\naturalnumbers}$ is the set of all real-valued sequences. Let $\tp$ be the topology of pointwise convergence on~${\mathbb R}^{\naturalnumbers}$ that corresponds to the standard topology on~${\mathbb R}$. The space $\left({\mathbb R}^{\naturalnumbers}, \tp \right)$ is second countable by Remark~\ref{rema r interval rel} and Lemmas~\ref{lemm count times count} and~\ref{lemm fin count}.
\eexam

\btheo
With definitions as in Remark~\ref{rema pointw conv gen}, let $(f_n)$ be a net in~$F$, $\mc F$ a filter on~$F$, and $f \in F$. The following statements hold:

\benum
\item $\big( f_n \rightarrow f \;\, \mbox{with respect to} \; \tp_F \big) \quad \Longleftrightarrow$\\[.2em]
${} \quad \big( \forall x \in X \quad f_n(x) \rightarrow f(x) \;\, \mbox{with respect to} \; \tp_Y \big)$
\item $\big( {\mc F} \rightarrow f \;\, \mbox{with respect to} \; \tp_F \big) \quad \Longleftrightarrow$\\[.2em]
${} \quad \big( \forall x \in X \quad q_x \llbracket {\mc F} \, \rrbracket \rightarrow f(x) \;\, \mbox{with respect to} \; \tp_Y \big)$
\eenum

\etheo

\bproof
This is a consequence of Theorem~\ref{theo conv induced} and Remark~\ref{rema pointw conv gen}.
\eproof

\section{Direct image topology}
\label{direct image topology}

\blede
\label{deffinaltop}
\index{Direct image topology}
\index{Topology!direct image}
\index{Topology!generated}
Given a set $X$, topological spaces $(Y_i, \tp_i)$ $(i \in I)$ where $I$ is an index set, and functions $f_i : Y_i \longrightarrow X$ $(i \in I)$, the system $\bigcap_{i \in I} \left\{B \subset X \, : \, f_i^{-1} \left[ B \right] \in \tp_i \right\}$ is a topology on~$X$. It is called {\bf direct image topology} or the topology {\bf generated by} $F = \left\{ (\tp_i, f_i) \, : \, i \in I \right\}$ and is denoted by~$\tau (F)$. It is the finest topology $\tp$ on~$X$ such that $f_i$ is $\tp_i\,$-$\tp$-continuous for every $i \in I$.
\elede

\bproof
$\tau (F)$ clearly has properties (\ref{defi top 1}) to (\ref{defi top 3}) in Definition~\ref{defi top}. Now let $\mathscr{A}$ be the set of all topologies $\tp$ on~$X$ such that $f_i$ is $\tp_i\,$-$\tp$-continuous for every $i \in I$. We clearly have $\tau(F) \in \mathscr{A}$. Moreover, for every $\tp \in \mathscr{A}$ we have $\tp \subset \tau (F)$. Hence $\tau(F)$~is the finest member of~$\mathscr{A}$.
\eproof

Notice that our convention is such that every member of~$F$ has the topology of the domain space as its left coordinate, in contrast to Lemma and Definition~\ref{definitialtop} where the members of the generating system have the topology of the range space as their right coordinate.

The following is an important special case.

\midvspace

\bcoro
\label{coro im top same set}
Let $X$ be a set, $I$~an index set, and, for each $i \in I$, $\tp_i$ a topology on~$X$. Further let $F = \left\{(\tp_i,\mathrm{id}_X) \, : \, i \in I \right\}$. We have $\tau (F) = \bigcap_{i \in I} \tp_i$\,. $\tau (F)$~is the infimum of $\left\{ \tp_i \, : \, i \in I \right\}$ in the ordered space $(\mathscr{T}(X),\subset)$, i.e.\ it is the finest topology on~$X$ that is coarser than $\tp_i$ for every $i \in I$.
\ecoro

\bproof
Exercise.
\eproof

Direct image topologies may be characterized by a universal property, similarly to the case of inverse image topologies (cf.\ Theorem~\ref{universalproperty}).

\midvspace

\btheo
\label{universal property final}
Let $(X,\tp)$ be a topological space, $I$~an index set, $(Y_i, \tp_i)$ $(i \in I)$ topological spaces, $f_i : Y_i \longrightarrow X$ $(i \in I)$ functions, and $F = \left\{ (\tp_i, f_i) \, : \, i \in I \right\}$. The following statements are equivalent:

\benum
\item \label{universal property final 1} $\tp = \tau (F)$
\item \label{universal property final 2} For every topological space $(Z,\tp_Z)$ and every function $g : X \longrightarrow Z$, $g$ is $\tp$-$\tp_Z\,$-continuous iff $g \circ f_i$ is $\tp_i\,$-$\tp_Z\,$-continuous for every $i \in I$.
\eenum

\etheo

\bproof
To see that (\ref{universal property final 1}) implies (\ref{universal property final 2}), assume that $\tau (F) = \tp$. Then $f_i$ is $\tp_i\,$-$\tp$-continuous for every $i \in I$. Further let $(Z,\tp_Z)$ be a topological space and \mbox{$g : X \longrightarrow Z$} a map. If~$g$ is continuous, then $g \circ f_i$ is continuous for every $i \in I$ by Lemma~\ref{lemm cont comp}. To show the converse let $U \in \tp_Z$. If $g \circ f_i$ is continuous for every $i \in I$, then we have $f_i^{-1} \left[ g^{-1} \left[ U \right] \right] \in \tp_i$ ($i \in I$), and therefore $g^{-1} \left[ U \right] \in \tp$. Thus $g$ is continuous.

To show that (\ref{universalproperty 2}) implies (\ref{universalproperty 1}), it is enough to show that the topology~$\tp$ is uniquely specified by property~(\ref{universalproperty 2}). Assume that $\tp_1$ and $\tp_2$ are two topologies on~$X$ such that (\ref{universalproperty 2}) is satisfied in both cases. Now let $Z = X$ and $g = 
\mathrm{id}_X$. Since $g$ is $\tp_m\,$-$\tp_m\,$-continuous for $m \in \left\{ 1, 2 \right\}$, it follows that $f_i$ is $\tp_i\,$-$\tp_m\,$-continuous for $m \in \left\{ 1, 2 \right\}$ and $i \in I$. Thus $g$ is $\tp_1\,$-$\tp_2\,$-continuous and $\tp_2\,$-$\tp_1\,$-continuous, and hence $\tp_1 = \tp_2$.
\eproof

The following result is a characterization of the direct image topology in the case of a single function.

\midvspace

\btheo
\label{theo char quotient}
Let $(X,\tp_X)$ and $(Y,\tp_Y)$ be two topological spaces and $f : X \longrightarrow Y$ a map. If $f$ is $\tp_X\,$-$\tp_Y\,$-continuous, surjective, and either $\tp_X\,$-$\tp_Y\,$-open or $\tp_X\,$-$\tp_Y\,$-closed, then we have $\tp_Y = \tau \left(\left\{ (\tp_X,f) \right\}\right)$.
\etheo

\bproof
Let $\tp = \tau \left(\left\{ (\tp_X,f) \right\}\right)$, and assume that $f$ is $\tp_X\,$-$\tp_Y\,$-continuous, surjective, and either open or closed.

Since $\tp$ is the finest topology on~$Y$ such that $f$ is $\tp_X\,$-$\tp$-continuous, we have $\tp_Y \subset \tp$.

Conversely, let $U \in \tp$. Then we have $f^{-1} \left[ U \right] \in \tp_X$. First consider the case that $f$ is $\tp_X\,$-$\tp_Y\,$-open. Then we have $U = f \left[ f^{-1} \left[ U \right] \right] \in \tp_Y$. Second, if $f$ is $\tp_X\,$-$\tp_Y\,$-closed we have
\[
U = \big( U^c \big)^c = \big( f \left[ f^{-1} \left[ U^c \right] \right] \big)^c = \Big( f \, \Big[ \big( f^{-1} \left[ U \right] \! \big)^c \, \Big] \Big)^c \in \tp_Y
\]

\eproof

\section{Quotient topology}
\label{quotient topology}

In this Section we analyse an important special case of the concept introduced in Section~\ref{direct image topology}.

\midvspace

\bdefi
\index{Quotient topology}
\index{Topology!quotient}
\index{Quotient space}
Let $(X,\tp)$ be a topological space, $R$~an equivalence relation on~$X$, and \mbox{$f : X \longrightarrow X/R$}, $f(x) = \left[ x \right]$. The topology $\tp_R = \tau \left(\left\{(\tp,f)\right\}\right)$ is called {\bf quotient topology}. The space $(X/R,\tp_R)$ is called {\bf quotient topological space}, or short {\bf quotient space}.
\edefi

\blemm
Let $X$ be a set, $(Y,\tp)$ a topological space, $f : Y \longrightarrow X$ a surjective function, and $R$ the equivalence relation on $Y$ defined by
\[
(y,z) \in R \quad \Longleftrightarrow \quad f(y) = f(z)
\]

Further let $\tp_R$ be the quotient topology on~$Y/R$ and $\tp_X = \tau \left(\left\{(\tp,f)\right\}\right)$. Then $(Y/R,\tp_R)$ and $(X,\tp_X)$ are homeomorphic.
\elemm

\bproof
We define the map $g : Y \longrightarrow Y/R$, $g(y) = \left[ y \right]$. Further let the map $h : Y/R \longrightarrow X$ be defined by $h(\left[ y \right]) = f(y)$ for every $y \in Y$. $h$~is clearly well-defined. We show that $h$ is a homeomorphism. $h$~is clearly bijective. We have $h \circ g = f$. Therefore the continuity of~$f$ implies the continuity of~$h$ by Theorem~\ref{universal property final}. Furthermore, we have $g = h^{-1} \circ f$. Hence the continuity of~$g$ implies the continuity of~$h^{-1}$ by the same Theorem. 
\eproof

\btheo
\label{theo pseudo-metric metric}
Let $(X,d)$ be a pseudo-metric space and $R = \left\{(x,y) \, : \, d(x,y) = 0 \right\}$. $R$~is an equivalence relation on~$X$. The map
\[
D : (X/R) \times (X/R) \longrightarrow {\mathbb R}_+\,, \quad D \left( \,\! \left[ x \right]\!, \,\! \left[ y \right] \,\! \right) = d(x,y)
\]

is a metric on~$X/R$. Moreover $\tau (D)$ is the quotient topology of~$\tau (d)$.
\etheo

\bproof
$R$~clearly is an equivalence relation on~$X$. The function~$D$ is well-defined because of the triangle inequality for~$d$. That $D$ is a pseudo-metric follows by the fact that $d$ is a pseudo-metric, and that $D$ is a metric is then obvious. In order to show that $D$ generates the quotient topology of~$\tau (d)$, it is enough by Theorem~\ref{theo char quotient} to show that the map $f : X \longrightarrow X/R$, $f(x) = \left[ x \right]$, is $\tau (d)$-$\tau (D)$-continuous and $\tau (d)$-$\tau (D)$-open. Since $f$ is an isometry, this follows by Lemma~\ref{isom cont open}.
\eproof


\chapter{Functions and real numbers}
\label{functions on Rn}
\setcounter{equation}{0}

\pagebreak

In this Chapter we use the various concepts introduced in Chapters~\ref{topologies} to~\ref{generated topologies} (topologies, pseudo-metrics, continuity, etc.)\ in the context of the number systems defined in Chapter~\ref{numbers ii}.

\midvspace

\bdefi
\label{defi convention top r r+}
We adopt the convention that all notions related to topologies on ${\mathbb R}$ and ${\mathbb R}^n$, and their subsets refer to the respective standard topologies or their relative topologies if not otherwise specified.
\edefi

In particular, according to this convention we refer to the standard topologies on ${\mathbb R}_+$ and~${\mathbb R}^n_+$ since these are the relative topologies by Lemma~\ref{lemm standard top r r+} and Corollary~\ref{coro standard top subspace}.

\midvspace

\blemm
\label{lemm real op cont}
The addition $+ : {\mathbb R}^2 \longrightarrow {\mathbb R}$, the absolute value $b : {\mathbb R} \longrightarrow {\mathbb R}_+$\,, the multiplication \mbox{$\cdot : {\mathbb R}^2 \longrightarrow {\mathbb R}$}, and, for every $m \in \naturalnumbers$, the exponentiation \mbox{$h_m : {\mathbb R}_+ \longrightarrow {\mathbb R}_+$}\,, $h_m (\alpha) = \alpha^m$, are continuous functions.
\elemm

\bproof
We show the continuity of each function by means of Lemma~\ref{lemm continuity sequence}. We use the bases for the respective standard topologies on~${\mathbb R}$ and~${\mathbb R}_+$ as given in Remarks~\ref{rema r interval rel} and~\ref{rema r+ interval rel}.

The continuity of~$b$ is clear.

Let $z \in {\mathbb R}^2$, $x$ and $y$ be its left and right coordinates, i.e.\ $z = (x,y)$, $(z_k)$ a sequence in~${\mathbb R}^2$ such that $z_k \rightarrow z$, and $(x_k)$ and $(y_k)$ the left and right coordinate sequences, i.e.\ $z_k = (x_k, y_k)$ for every $k \in \naturalnumbers$. It follows that $x_k \rightarrow x$ and $y_k \rightarrow y$ by Corollary~\ref{coro conv product filter}.

To see that addition is continuous in~$z$, let $u, v \in {\mathbb R}$ such that $u < x + y < v$. We define $w = \frac{1}{2} \, \mathrm{min} \left\{ x+y-u,\, v-x-y \right\}$. There is $n \in \naturalnumbers$ such that
\[
x - w < x_k < x + w, \quad \quad y - w < y_k < y + w
\]

for every $k \geq n$. Hence we have $u < x_k + y_k < v$ for $k \geq n$.

To show that multiplication is continuous in~$z$, let $u, v \in {\mathbb R}$ such that $u < xy < v$. We define $w = \frac{1}{2} \, \mathrm{min} \left\{ v - xy,\, xy - u \right\}$. For $k \in \naturalnumbers$, we have
\begin{eqnarray*}
|x_k y_k - xy| & = & |x_k y_k - x_k y + x_k y - xy|\\
 & \leq & |x_k| |y_k - y| + |x_k - x| |y|
\end{eqnarray*}

by Remark~\ref{rema prop abs val} and Lemma~\ref{lemm b metric}. We may choose $K \in {\mathbb R} \!\setminus\! \left\{ 0 \right\}$ such that $|x| < K$. There is $n \in {\mathbb N}$ such that $|x_k| < K$, $|y_k - y| < w K^{-1}$, and, if $y \neq 0$, $|x_k - x| < w |y|^{-1}$ for every $k \geq n$. It follows that, for $k \geq n$, $|x_k y_k - xy| < 2w$, and thus $u < x_k y_k < v$.

Finally, we show that $h_m$ is continuous for every $m \in \naturalnumbers$ by the Induction principle. The case $m = 0$ is clear. Assume that $h_m$ is continuous for some $m \in \naturalnumbers$. The function $f : {\mathbb R}_+ \longrightarrow {\mathbb R}_+^2$\,, $f(\alpha) = (h_m(\alpha), \alpha)$, is continuous by Lemma~\ref{lemmm cont prod} and Remark~\ref{rema standard prod top}. Since the multiplication on~${\mathbb R}$ is continuous, the multiplication on~${\mathbb R}_+$\,, $g : {\mathbb R}_+^2 \longrightarrow {\mathbb R}_+$\,, $g(\alpha, \beta) = \alpha \cdot \beta$, is continuous by Corollary~\ref{coro standard top subspace} and Lemma~\ref{continuity subspace}. Thus $h_{m+1} = g \circ f$ is continuous by Lemma~\ref{composition filter cont}.
\eproof

\bprop
\label{prop bound sequ r}
Let $x, y \in {\mathbb R}$ and $(x_n)$ be a sequence in ${\mathbb R}$ such that $x_n \rightarrow x$. The following statements hold:

\benum
\item \label{prop bound sequ r 1} $\left( \forall n \in \naturalnumbers \quad x_n \leq y \right) \quad \Longrightarrow \quad x \leq y$
\item \label{prop bound sequ r 2} $\left( \forall n \in \naturalnumbers \quad y \leq x_n \right) \quad \Longrightarrow \quad y \leq x$
\eenum

\eprop

\bproof
To see~(\ref{prop bound sequ r 1}), we define $A = \left] - \infty, y \right]$. $A$~is closed with respect to the standard topology on~${\mathbb R}$. Under the stated condition we have $x_n \in A$ for every $n \in \naturalnumbers$, and thus $x \in A$ by Lemma~\ref{lemm first count sequ}~(\ref{lemm first count sequ 2}).

The proof of~(\ref{prop bound sequ r 2}) is similar.
\eproof

\bprop
\label{prop max min r}
Let $a, b, c, w \in {\mathbb R}$ such that $\mathrm{min} \left\{ a, c \right\} \leq w \leq \mathrm{max} \left\{ a, c \right\}$. Then we have $\mathrm{min} \left\{ a, b \right\} \leq w \leq \mathrm{max} \left\{ a, b \right\}$ \, or \, $\mathrm{min} \left\{ b, c \right\} \leq w \leq \mathrm{max} \left\{ b, c \right\}$.
\eprop

\bproof
Exercise.
\eproof

\btheo[Intermediate value]
\label{theo interm value}
\index{Intermediate value theorem}
Let $x, y \in {\mathbb R}$ with $x < y$, $A = \left[ x, y \right]$, $f : A \longrightarrow {\mathbb R}$ a continuous function, $B = \mathrm{ran} \, f$, and $u = \mathrm{min} \left\{ f(x), f(y) \right\}$, $v = \mathrm{max} \left\{ f(x), f(y) \right\}$. The following statements hold:

\benum
\item \label{theo interm value 1} $u < v \;\; \Longrightarrow \;\; \left[ u, v \right] \subset B$
\item \label{theo interm value 2} If $f$ is strictly monotonic, then $\left[ u, v \right] = B$.
\item \label{theo interm value 3} We define the map $g : A \longrightarrow B$, $g(z) = f(z)$. If $g$ is strictly monotonic, then it is bijective and $g^{-1}$ is continuous. If $g$ is strictly increasing (strictly decreasing), then $g^{-1}$ is strictly increasing (strictly decreasing).
\eenum

\etheo

\bproof
To see~(\ref{theo interm value 1}), assume that $u < v$. Let $w \in \left[ u, v \right]$. We define two sequences $(x_n)$ and $(y_n)$ in~$A$ by $(x_0,y_0) = (x,y)$, and recursively for every $n \in \naturalnumbers$,

\vspace{.1in}
\begin{minipage}{\alltheoremwidth}
\[
(x_{n+1},y_{n+1}) = \left\{
\begin{array}{ll}
(x_n,z) \quad & \mathrm{if} \quad s \leq w \leq t\\[.5em]
(z,y_n) \quad & \mathrm{else}
\end{array}
\right.
\]

where
\[
z = \frac{1}{2}(x_n + y_n) \,, \quad
s = \mathrm{min} \left\{ f(x_n), f(z) \right\} \,, \quad
t = \mathrm{max} \left\{ f(x_n), f(z) \right\}
\]

\end{minipage}
\vspace{.1in}

It follows by the Induction principle that \mbox{$y_n - x_n = (y_0 - x_0) / 2^n$} for every $n \in \naturalnumbers$. Thus, for every $n \in \naturalnumbers$, we have $x_n < y_n$\,, whence $x_n < (x_n + y_n)/2 < y_n$\,, and thus $x_n \leq x_{n+1}$ and $y_{n+1} \leq y_n$\,. Hence $(x_n)$ is increasing and $(y_n)$ is decreasing by the Induction principle. Let $\tp$ be the standard topology on~${\mathbb R}$. Now regarding $(x_n)$ and $(y_n)$ as sequences in the whole of~${\mathbb R}$, they are increasing and decreasing, respectively, and bounded. Thus they are convergent with respect to~$\tp$ by Lemma~\ref{lemm increasing convergent} and Remark~\ref{rema r increasing conv}, say $x_n \rightarrow x_{\infty}$ and $y_n \rightarrow y_{\infty}$ where $x_{\infty}\,,\, y_{\infty} \in X$. Therefore $x_{\infty} = y_{\infty}$ by Lemmas~\ref{lemm real op cont} and~\ref{lemm sequ inv conv}. Since $A$ is $\tp$-closed, we have $x_{\infty} \in A$. Again regarding $(x_n)$ and $(y_n)$ as sequences in~$A$, it follows that $x_n \rightarrow x_{\infty}$ and $y_n \rightarrow x_{\infty}$ with respect to~$\tp \, | A$ by Lemma~\ref{conv subspace net}.

The continuity of $f$ implies $f(x_n) \rightarrow f(x_{\infty})$ and $f(y_n) \rightarrow f(x_{\infty})$. Moreover we have for every $n \in \naturalnumbers$
\[
\mathrm{min} \left\{ f(x_n), f(y_n) \right\} \, \leq \, w \, \leq \, \mathrm{max} \left\{ f(x_n), f(y_n) \right\}
\]

by the Induction principle and Proposition~\ref{prop max min r}. It follows that $w = f(x_{\infty})$ by Proposition~\ref{prop bound sequ r}.

To see~(\ref{theo interm value 2}), note that, if $f$ is strictly monotonic, then clearly $B \subset \left[ u, v \right]$.

To see~(\ref{theo interm value 3}), notice that the map $g$ is surjective by definition. Now assume that $g$ is strictly monotonic. Then $g$ is clearly injective. To see that $g^{-1}$ is continuous, note that the system
\[
{\mc S}_A = \big\{ \left[ x, z \right[\; : \, z \in \left] x, y \right] \big\} \cup \big\{ \left] z, y \right] \, : \, z \in \left[ x, y \right[ \big\}
\]

is a subbase for~$\tp \, | A$, and
\[
{\mc S}_B = \big\{ \left[ u, z \right[\; : \, z \in \left] u, v \right] \big\} \cup \big\{ \left] z, v \right] \, : \, z \in \left[ u, v \right[ \big\}
\]

is a subbase for~$\tp \, | \, B$, cf.\ Remark~\ref{rema rel top int r}. We have $g \, \llbracket {\mc S}_A \rrbracket \subset \, {\mc S}_B$\,, and thus $g^{-1}$ is continuous by Theorem~\ref{theo global cont}~(\ref{theo global cont 2}).
\\

\hspace{0.05\textwidth}
\parbox{0.95\textwidth}
{[If $g$ is strictly increasing, we have for every $z \in \left] x, y \right]$:
\begin{eqnarray*}
g \big[ \left[ x, z \right[\, \big] \!\!\!\! & = & \!\!\! 
g \big[ \left[ x, z \right] \setminus \left\{ z \right\} \big] = \,
g \big[ \left[ x, z \right] \big] \setminus \left\{ g(z) \right\}\\[.2em]
 & = & \!\!\!
\left[ g(x), g(z) \right] \setminus \left\{ g(z) \right\} \, = \, \left[ u, g(z) \right[
\end{eqnarray*}

where the third equation is a consequence of~(\ref{theo interm value 2}) and Lemma~\ref{continuity subspace}.

Moreover if $g$ is strictly increasing, we have for every $z \in \left[ x, y \right[\,$:
\[
g \big[ \,\left] z, y \right] \big] = \; \left] g(z), v \right]
\]

If $g$ is strictly decreasing, we have

\begin{center}
\begin{tabular}{l}
$\forall z \in \left] x, y \right] \quad g \big[ \left[ x, z \right[\, \big] = \;\left] g(z), v \right]$ ,\\[.5em]
$\forall z \in \left[ x, y \right[ \quad g \big[ \,\left] z, y \right] \big] = \,\left[ u, g(z) \right[$
\end{tabular}
\end{center}]}
\\

The last two claims are now obvious.
\eproof

\bcoro
\label{coro strictly inverse}
Let $f : {\mathbb R}_+ \longrightarrow {\mathbb R}_+$ be a map with $f(0) = 0$. If $f$ is continuous, strictly increasing, and unbounded, then $f$ is bijective and $f^{-1}$ is continuous and strictly increasing.
\ecoro

\bproof
Assume the stated conditions. $f$~is clearly injective. To see that it is surjective, let $y \in {\mathbb R}_+$\,. Since $f$ is unbounded, there is $x \in {\mathbb R}_+$ such that $f(x) > y$. It follows that $y \in \left[ 0, f(x) \right] = f \big[ \left[ 0, x \right] \big]$ by Lemma~\ref{continuity subspace}, and Theorem~\ref{theo interm value}~(\ref{theo interm value 2}). Moreover $f^{-1}$ is clearly strictly increasing. 

Finally we show that $f^{-1}$ is continuous. For every $m \in \naturalnumbers$ we define $A_m = \left[ m, m+1 \right]$ and the function
\[
f_m : A_m \longrightarrow {\mathbb R}, \quad f_m(x) = f(x)
\]

and $B_m = \mathrm{ran} \, f_m$\,. Clearly these maps are strictly increasing and continuous by Lemma~\ref{continuity subspace}. We have $B_m = \left[ f_m(m), f_{m+1}(m+1) \right]$ \, ($m \in \naturalnumbers$) by Theorem~\ref{theo interm value}~(\ref{theo interm value 2}). Further for every $m \in \naturalnumbers$ we define the map
\[
g_m : A_m \longrightarrow B_m\,, \quad g_m(x) = f_m(x)
\]

The functions $g_m$ ($m \in \naturalnumbers$) are strictly increasing. By Theorem~\ref{theo interm value}~(\ref{theo interm value 3}), for every $m \in \naturalnumbers$, $g_m$ is bijective and $g_m^{-1}$ is continuous. For every $m \in \naturalnumbers$ we define
\[
t_m : B_m \longrightarrow {\mathbb R}_+\,, \quad t_m(y) = g_m^{-1}(y)
\]

The functions $t_m$ ($m \in \naturalnumbers$) are continuous by Lemma~\ref{continuity subspace}. We have $t_m = f^{-1}| \, B_m$ ($m \in \naturalnumbers$).
\\

\hspace{0.05\textwidth}
\parbox{0.95\textwidth}
{[This is seen as follows:
\begin{eqnarray*}
(y,x) \in t_m & \Longleftrightarrow & (y,x) \in g^{-1}_m \;\; \Longleftrightarrow \;\; (x,y) \in g_m \;\; \Longleftrightarrow \;\; (x,y) \in f_m\\
 & \Longleftrightarrow & (x,y) \in f \;\; \wedge \;\; x \in A_m\\
 & \Longleftrightarrow & (y,x) \in f^{-1} \;\; \wedge \;\; y \in B_m\\
 & \Longleftrightarrow & (y,x) \in f^{-1} | \, B_m
\end{eqnarray*}
]}
\\

Thus $f^{-1}$ is continuous by Theorem~\ref{theo cont closed system}.
\eproof

\blede
\index{m-th root function}
\index{Root function}
For every $m \in \naturalnumbers$ with $m > 0$ let $h_m : {\mathbb R}_+ \longrightarrow {\mathbb R}_+$\,, $h_m (\alpha) = \alpha^m$. The map $h_m$ is bijective. $h_m^{-1}$~is called {\bf $m$-th root function}. $h_2^{-1}$~is called {\bf square root function}. $h_m^{-1}$~is continuous and strictly increasing for $m \in \naturalnumbers$, $m > 0$. We also write \,\! $\alpha^{\frac{1}{m}}$, $\alpha^{1/m}$, or \,\! $\sqrt[m]{\alpha}$ \,\! for \,\! $h_m^{-1}(\alpha)$. Moreover, we also write $\sqrt{\alpha}$ \,\! for \,\! $h_2^{-1}(\alpha)$. The value $\alpha^{1/m}$ is called {\bf $m$-th root of}~$\alpha$. The value $\alpha^{1/2}$ is called {\bf square root of}~$\alpha$.

For $\alpha, \beta \in {\mathbb R}_+$ and $m \in \naturalnumbers$ with $m > 0$, we have $\left( \alpha \, \beta \right)^{1/m} = \alpha^{1/m} \, \beta^{1/m}$.
\elede

\bproof
Let $m \in \naturalnumbers$ with $m > 0$. Then $h_m$ is continuous by Lemma~\ref{lemm real op cont}, and strictly increasing and unbounded by Lemma~\ref{lemm op increasing}. Furthermore, we have $h_m(0) = 0$. It follows by Corollary~\ref{coro strictly inverse} that $h_m$ is bijective, and that $h_m^{-1}$ is strictly increasing and continuous.

To see the last claim notice that for every $\alpha, \beta \in {\mathbb R}_+$ and every $m \in \naturalnumbers$ with $m > 0$ we have
\[
\big( \alpha^{1/m} \, \beta^{1/m} \big)^m \! = \, \alpha \, \beta
\]

by Lemma and Definition~\ref{lede exp pos reals}.
\eproof

\bdefi
\label{defi finite series}
\index{Finite series}
\index{Series!finite}
Let $n \in \naturalnumbers$ and $x \in {\mathbb R}^n$. Further let the function $S : \sigma(n) \!\setminus\! \left\{ 0 \right\} \longrightarrow {\mathbb R}$ be recursively defined by $S(1) = x_1$ and $S(k + 1) = S(k) + x_{k+1}$ for $1 \leq k \leq n-1$. For every $m \in \naturalnumbers$ with $1 \leq m \leq n$, $S(m)$~is called a {\bf finite series} and denoted by
\[
{\displaystyle \sum_{k=1}^m} x_k
\]

If $n \geq 2$, then we write ${\displaystyle \sum_{k=l}^m} x_k$ for $S(m) - S(l-1)$ where $l, m \in \naturalnumbers$ with \mbox{$2 \leq l \leq m \leq n$}.
\edefi

Note that Definition~\ref{defi finite series} is based on the Local recursion theorem~\ref{theo recursion}. The Recursion theorem for natural numbers, Theorem~\ref{theo recursive def}, does not suffice.

\midvspace

\blemm
\label{lemm calc finite series}
Let $n \in \naturalnumbers$, $x, y \in {\mathbb R}^n$, and $z \in {\mathbb R}$. Then the following statements hold:

\benum
\item \label{lemm calc finite series 1} ${\displaystyle \sum_{k=1}^n} (x_k + y_k) \, = \, \left( {\displaystyle \sum_{k=1}^n} x_k \right) + \left( {\displaystyle \sum_{k=1}^n} y_k \right)$
\item \label{lemm calc finite series 2} ${\displaystyle \sum_{k=1}^n} (z \cdot x_k) \, = \, z \cdot \left( {\displaystyle \sum_{k=1}^n} x_k \right)$
\item \label{lemm calc finite series 3} $x_k  > 0 \;\; (1 \leq k \leq n) \quad\! \Longrightarrow \quad {\displaystyle \sum_{k=1}^n} x_k \, > \, 0$
\item \label{lemm calc finite series 4} $x_k \geq 0 \;\; (1 \leq k \leq n) \quad\! \Longrightarrow \quad {\displaystyle \sum_{k=1}^n} x_k \, \geq \, 0$
\item \label{lemm calc finite series 5} $\left( x_k \geq 0 \;\; (1 \leq k \leq n) \;\; \wedge \;\; {\displaystyle \sum_{k=1}^n} x_k \, = \, 0 \right) \quad \Longrightarrow \quad x_k = 0 \;\; (1 \leq k \leq n)$
\eenum

\elemm

\bproof
Exercise.
\eproof

\midvspace

\blemm
\label{lemm rel rn}
\index{Cauchy-Schwarz inequality}
\index{Inequality!Cauchy-Schwarz}
Let $x, y \in \mathbb R^n$. The following statements hold:

\benum
\item \label{lemm rel rn 1} $\left( {\displaystyle \sum_{k=1}^n} x_k y_k \right)^2 \, \leq \, \left( {\displaystyle \sum_{k=1}^n} x_k^2 \right)\left( {\displaystyle \sum_{k=1}^n} y_k^2 \right)$ \quad (Cauchy-Schwarz inequality)
\item \label{lemm rel rn 2} $\left( {\displaystyle \sum_{k=1}^n} (x_k + y_k)^2 \right)^{1/2} \leq \, \left( {\displaystyle \sum_{k=1}^n} x_k^2 \right)^{1/2} \!\! + \left( {\displaystyle \sum_{k=1}^n} y_k^2 \right)^{1/2}$
\eenum

\elemm

\bproof
We first prove~(\ref{lemm rel rn 1}). We have
\begin{eqnarray*}
0 \!\! & \leq & \!\! \sum_{k=1}^n \sum_{l=1}^n \left(x_k y_l - x_l y_k\right)^2 \; = \; \sum_{k=1}^n \sum_{l=1}^n \left( x_k^2 y_l^2 + x_l^2 y_k^2 - 2 x_k x_l y_k  y_l \right)\\[.2em]
 & = & \!\! 2 \left( \sum_{k=1}^n x_k^2 \right)\left( \sum_{l=1}^n y_l^2 \right) - 2 \left( \sum_{k=1}^n x_k y_k \right) \left( \sum_{l=1}^n x_l y_l \right)
\end{eqnarray*}

In order to show~(\ref{lemm rel rn 2}), notice that
\begin{eqnarray*}
\sum_{k=1}^n (x_k + y_k)^2 \!\! & = & \!\! \sum_{k=1}^n x_k^2 \, + \, \sum_{k=1}^n y_k^2 \, + \, 2 \sum_{k=1}^n x_k y_k\\[.2em]
 & \leq & \!\! \sum_{k=1}^n x_k^2 \, + \, \sum_{k=1}^n y_k^2 \, + \, 2 \left(\sum_{k=1}^n x_k^2 \right)^{1/2} \! \left(\sum_{k=1}^n y_k^2 \right)^{1/2}\\[.2em]
 & = & \!\! \left(\left(\sum_{k=1}^n x_k^2 \right)^{1/2} \! + \, \left(\sum_{k=1}^n y_k^2 \right)^{1/2}\right)^2
\end{eqnarray*}

by the Cauchy-Schwarz inequality and the fact that the square root function is increasing.
\eproof

\blede
\index{Euclidean metric}
\index{Metric!Euclidean}
For $n \in \mathbb{N}$ with $n > 0$ the map
\[
d : {\mathbb R}^n \!\times {\mathbb R}^n \longrightarrow {\mathbb R}_+\,, \quad d(x,y) = \left( {\displaystyle \sum_{k=1}^n} |x_k - y_k|^2 \right)^{1/2}
\]

is a metric. It is called {\bf Euclidean metric}.
\elede

\bproof
Notice that $d$ satisfies the triangle inequality by Lemmas~\ref{lemm b metric} and~\ref{lemm rel rn}~(\ref{lemm rel rn 2}), and thus it is a pseudo-metric. Further $d$ is a metric by Lemma~\ref{lemm calc finite series}~(\ref{lemm calc finite series 5}).
\eproof

\blemm
\label{lemm sup int top equiv}
Let $n \in \mathbb{N}$ with $n > 0$, $d$~the maximum metric on~${\mathbb R}^n$, and $\tp$ the standard topology on~${\mathbb R}^n$. Then $\tau(d) = \tp$. 
\elemm

\bproof
We define the function

\begin{center}
\begin{tabular}{l}
$B : {\mathbb R}^n \, \times \;]0,\infty[\; \longrightarrow {\mc P}(X)$,\\[.8em]
$B(x,r) = \big\{ y \in {\mathbb R}^n \, : \, |x_k - y_k| < r \;\; (1 \leq k \leq n) \big\}$
\end{tabular}
\end{center}

The system
\[
{\mc B} = \big\{ B(x,r) \, : \, x \in {\mathbb R}^n,\; r \in \;]0,\infty[\; \big\} \cup \left\{ \O \right\}
\]

is a base for $\tau(d)$ by definition. Moreover the system
\[
{\mc A} = \big\{ \left] x, y \right[ \; : \, x, y \in \mathbb R^n,\; (x,y) \in S \big\} \cup \left\{ \O \right\}
\]

is a base for~$\tp$ where the interval refers to the ordering $S$ on~${\mathbb R}^n$ as defined in Remark~\ref{rema rn interval rel}. We clearly have ${\mc B} \subset {\mc A}$ and ${\mc A} \subset \Theta ({\mc B})$.
\eproof

\blemm
\label{lemm pseudo-m equiv top}
Let $n \in \naturalnumbers$ with $n > 0$. Let $e$ be the Euclidean metric on~${\mathbb R}^n$ and $d$ the maximum metric on~${\mathbb R}^n$. We have $\tau(d) = \tau(e)$.
\elemm

\bproof
Notice that we have $e(x,y) \leq \sqrt{n} \, d(x,y)$ and $d(x,y) \leq e(x,y)$ for every $x, y \in {\mathbb R}^n$. The claim follows by Lemma~\ref{lemm pseudo-m comp}.
\eproof

\bcoro
For every $n \in \naturalnumbers$ with $n > 0$, the topology generated by the Euclidean metric on~$\mathbb R^n$ is the standard topology.
\ecoro

\bproof
This is a consequence of Lemmas~\ref{lemm sup int top equiv} and~\ref{lemm pseudo-m equiv top}.
\eproof

\bcoro
\label{coro sequ r unique limit}
Let $(x_n)$ be a sequence in~$\mathbb R$. If $(x_n)$ is convergent, then it has a unique limit point.
\ecoro

\bproof
Let $x, y \in {\mathbb R}$ be two limit points of~$(x_n)$. For every $n \in \naturalnumbers$ there is $m \in \naturalnumbers$ such that $|x_m - x| < 1/n$ and $|x_m - y| < 1/n$ by Remark~\ref{rema equiv net r}. It follows that $|x - y| \leq |x - x_m| + |x_m - y| < 2/n$. Therefore we have $|x - y| = 0$ by Lemma~\ref{lemm sequ inv conv} and Proposition~\ref{prop bound sequ r}~(\ref{prop bound sequ r 2}), and thus $x = y$.
\eproof

\bprop
\label{prop power xm}
For every $x \in {\mathbb R}$ with $x > 1$ and every $m \in \naturalnumbers$ we have $x^m \geq m (x - 1) + 1$.
\eprop

\bproof
The claim is clear for $m = 0$ and every $x \in {\mathbb R}$ with $x > 1$. Now assume that it holds for some $m \in \naturalnumbers$ and every $x \in {\mathbb R}$ with $x > 1$. Then we have
\begin{eqnarray*}
x^{m+1} \!\!\! & = & \!\! x^m \cdot x \; \geq \; \big( m (x-1) + 1 \big) \cdot x \; > \; m (x - 1) + x\\
 & = & \!\! m (x - 1) + (x - 1) + 1 \; = \; (m + 1)(x - 1) + 1
\end{eqnarray*}

\eproof

\bprop
\label{prop seq power}
Let $x \in {\mathbb R}_+$\,. The sequence $\left( x^m : m \in \naturalnumbers \right)$ is unbounded if $x > 1$, and converges to~$0$ if $0 < x < 1$.
\eprop

\bproof
If $x > 1$, then $x^m$ is unbounded by Proposition~\ref{prop power xm} and Lemma~\ref{lemm op increasing}. If $0 < x < 1$, then $x^m \rightarrow 0$ by Lemma~\ref{lemm sequ inv conv}.
\eproof

\blemm[Finite geometric series]
\label{lemm finite geom series}
\index{Finite geometric series}
\index{Geometric series}
\index{Series!finite geometric}
Let $x \in {\mathbb R}$. We have
\[
\sum_{k=0}^m x^k \, = \, \frac{\displaystyle 1 - x^{m+1}}{\displaystyle 1 - x}
\]

\elemm

\bproof
The equality is clearly true for $m = 0$. Assuming that it is true for some $m \in \naturalnumbers$, we have
\[
\sum_{k=0}^{m+1} x^k \, = \, \frac{\displaystyle 1 - x^{m+1}}{\displaystyle 1 - x} + x^{m+1} \, = \, \frac{\displaystyle 1 - x^{m+2}}{\displaystyle 1 - x}
\]

\eproof

\bdefi
\label{defi finite series seq}
Let $(x_n)$ be a sequence in~${\mathbb R}$. We define the function $S : \mathbb{N} \longrightarrow {\mathbb R}$ recursively by $S(0) = x_0$\,, and $S(k+1) = S(k) + x_{k+1}$ ($k \in \naturalnumbers$). For every $m \in \naturalnumbers$, $S(m)$~is called a {\bf finite series} and denoted by
\[
\sum_{k=0}^m x_k
\]

Moreover, we write ${\displaystyle \sum_{k=m}^n} x_k$ for $S(n) - S(m-1)$ where $m, n \in \naturalnumbers$ and $0 < m \leq n$.
\edefi

Similarly to Definition~\ref{defi finite series}, Definition~\ref{defi finite series seq} is based on the Local recursion theorem~\ref{theo recursion}.

\midvspace

\blede
\label{lede infinite sequence}
\index{Infinite series}
\index{Series!infinite}
Let $(x_n)$ be a sequence in~${\mathbb R}$. If the sequence $\big( \sum_{k=0}^m x_k \big)_m$ is convergent, its limit point is unique and denoted by
\[
\sum_{k=0}^{\infty} x_k
\]

It is called {\bf infinite series}.
\elede

\bproof
The uniqueness follows by Corollary~\ref{coro sequ r unique limit}.
\eproof

\blemm[Geometric series]
\index{Geometric series}
\index{Series!geometric}
The sequence $\big( \sum_{k=0}^m x^k \big)_m$ has a limit point for $x \in \left[ 0, 1 \right[\;$. In this case we have
\[
\sum_{k=0}^{\infty} x^k = \frac{1}{1 - x}
\]

\elemm

\bproof
We have
\[
{\textstyle \lim_m} \sum_{k=0}^m x^k \, = \, {\textstyle \lim_m} \frac{\displaystyle 1 - x^{m+1}}{\displaystyle 1 - x} \, = \, \frac{\displaystyle 1}{\displaystyle 1 - x}
\]

by Lemmas~\ref{lemm finite geom series} and~\ref{lemm real op cont}, and Proposition~\ref{prop seq power}.
\eproof


\backmatter


\cleardoublepage
\phantomsection
\addcontentsline{toc}{chapter}{Index}
\printindex

\end{document}